\def\versiondate{V 1.1 April 1999} 
 
\input math.macros
\input Ref.macros
\newcount\constnum
\constnum=0
\newcount\consSer 
\consSer=0
\def\OSmark#1{const.#1.OS}
\def\OSlbl#1#2{\expandafter\xdef\csname \OSmark{#1}\endcsname{#2}}
\def\OSref#1{\csname\OSmark{#1}\endcsname}
\def\nco#1{\OSlbl{\the\consSer#1}{\the\constnum}%
\global\advance\constnum by 1\co{#1}\relax}
\def\co#1{\theConst{\OSref{\the\consSer#1}}}
\def\theConst#1{C_{#1}}
\def\resetConsts{\global\advance\consSer by 1\constnum=0\relax}

\input EPSfig.macros
  \ifx\LabelFigloaded\MYundefined\relax
  \else
    \message{ !!! labelfig.tex ALREADY loaded !!!}
   \fi

  \def\LabelFigloaded{\relax}


  \chardef\LabelFigCatAt\the\catcode`\@
  \catcode`\@=11

 \let\LabelFigwlog@ld\wlog
 \def\wlog#1{\relax}

 \ifx\\\MYundefined@
    \let\\\relax
 \fi


  \def\ms@g{\immediate\write16}

 \def\N@wif{\csname newif\endcsname }
 \def\Temp@ {\N@wif\ifIN@}
 \ifx\INN@\MYundefined@
    \else \let\Temp@\relax
 \fi
 \Temp@

  \def\IN@{\expandafter\INN@\expandafter}
  \long\def\INN@0#1@#2@{\long\def\NI@##1#1##2##3\ENDNI@
    {\ifx\m@rker##2\IN@false\else\IN@true\fi}%
     \expandafter\NI@#2@@#1\m@rker\ENDNI@}
  \def\m@rker{\m@@rker}
 
  \newtoks\Initialtoks@  \newtoks\Terminaltoks@
  \def\SPLIT@{\expandafter\SPLITT@\expandafter}
  \def\SPLITT@0#1@#2@{\def\TTILPS@##1#1##2@{%
     \Initialtoks@{##1}\Terminaltoks@{##2}}\expandafter\TTILPS@#2@}

 \def\Shifted@@#1#2#3{\setbox0=\hbox{#3}%
   \raise -\dp0\vbox {\kern-#2%
       \hbox {\kern#1\unhbox0\kern-#1}%
           \kern#2}}

 \newcount\gridcount
 \newbox\auxGridbox@ \newbox\hGridbox@ \newbox\vGridbox@
 \newbox\Labelbox@ \newbox\auxLabelbox@
 \newbox\Coordinatebox@
 \newtoks\Labeltoks@
 \newdimen\Wdd@ \newdimen\Htt@
 \newdimen\Wddd@ \newdimen\Httt@
 
 \def\Wr@{\immediate\write16}

 \newdimen\GL@wd
 \GL@wd=.02pt
 \def\GridLineWidth#1{\GL@wd=#1}

 \def\gobble#1{}
 \def\EdgeErr@{\Wr@{}%
      \Wr@{\string\Edges\space argument
      1, 10, 100 or 1000 please\string!}%
      }

 \newcount\Edgect@

 \def\Sweepup#1\endSweepup{}

 \def\SetEdges@{%
    \edef\Zr@@s{\expandafter\gobble\number\Edgect@\empty}%
        \count255=0\Zr@@s\relax
        \ifnum\count255=\z@\else\EdgeErr@\show\tailtest\fi
        \count255=1\Zr@@s\relax
        \ifnum\count255=\Edgect@\relax\else\EdgeErr@\show\leadtest\fi
    \EdgGl@b\edef\Zr@s{\expandafter\gobble\Zr@@s\empty}
    \ifnum\Edgect@>\@ne\relax\EdgGl@b\let\L@Dc\empty
        \else\EdgGl@b\edef\L@Dc{\string.}\fi
    \ifnum\Edgect@>\@ne\relax
        \EdgGl@b\edef\Edgescale@##1{\divide##1 by \Edgect@}%
        \else\EdgGl@b\edef\Edgescale@##1{}\fi
    }

 \def\Edges#1{\Edgect@=#1\relax
     \let\EdgGl@b\global \SetEdges@}

 \Edges{1}

 \def\hhrule{\hrule height \GL@wd\vskip-.\GL@wd}

 \def\hRule@{%
   \advance\gridcount -2%
   \vfil\hhrule\vfil
   \llap{\smash{\raise -2.5pt
     \hbox{\L@Dc\number\gridcount\Zr@s\kern2pt}}}%
   \hhrule
   }

\def\vvrule{\vrule width \GL@wd \kern-\GL@wd}

 \def\vRule@{\advance\gridcount 2%
   \hfil\vvrule\hfil
   \setbox\auxGridbox@=\vbox to 0pt
      {\vskip \Htt@\vskip 2pt
        \hbox to 0pt{\hss\L@Dc\number\gridcount\Zr@s\hss}\vss}%
      \wd\auxGridbox@=0pt \box\auxGridbox@
   \vvrule
   }

 \def\PlaceGrid@@{\gridcount=10 
  \setbox\hGridbox@=\hbox{%
        \hbox{%
             \hskip-.4pt\vrule
             \vbox to \Htt@{%
               \offinterlineskip\parindent=\z@\relax
               \hbox to \Wdd@{\hfil}
               \hRule@\hRule@\hRule@\hRule@
               \vfil\hhrule\vfil}%
             \vrule\hskip-.4pt}
    }%
  \gridcount=0%
  \setbox\vGridbox@=\hbox{%
      \vbox{\offinterlineskip\parindent=0pt\hsize=0pt
         \vskip-.4pt\hrule%
         \hbox to \Wdd@{%
                 \vtop to \Htt@{\vfil}%
                 \vRule@\vRule@\vRule@\vRule@
                 \hfil\vvrule\hfil}%
         \hrule\vskip-.4pt}}%
  \wd\hGridbox@=0pt\ht\hGridbox@=0pt
  \wd\vGridbox@=0pt\ht\vGridbox@=0pt
  \hbox{\box\hGridbox@\box\vGridbox@}%
  }

 \def\LabelsGlobal{\def\LabGl@b{\global}}
 \def\LabelsLocal{\def\LabGl@b{}}
 \LabelsGlobal 

 \def\SetLabels#1\endSetLabels{%
   \LabGl@b\Labeltoks@={#1()\\}%
   }

 \LabGl@b\Labeltoks@={()\\}

 \def\ShowGrid{\LabGl@b\let\PlaceGrid@\PlaceGrid@@}
 \def\HideGrid{\LabGl@b\let\PlaceGrid@\relax}
 \def\Grids{\ShowGrid\LabGl@b\let\GridSwitch@\ShowGrid}
 \def\noGrids{\HideGrid\LabGl@b\let\GridSwitch@\HideGrid}

 \noGrids

 \def\bAdjust@@{%
     \setbox\auxLabelbox@=\hbox{\raise \dp\auxLabelbox@
            \box\auxLabelbox@}}
 \def\bAdjust@{\let\vAdjust@\bAdjust@@}

 \def\eAdjust@@{\dimen0=-.5\ht\auxLabelbox@
     \advance\dimen0 by .5\dp\auxLabelbox@
     \setbox\auxLabelbox@=
            \hbox{\raise\dimen0\box\auxLabelbox@}}
 \def\eAdjust@{\let\vAdjust@\eAdjust@@}

 \def\tAdjust@@{%
     \setbox\auxLabelbox@=\hbox{\raise-\ht\auxLabelbox@
            \box\auxLabelbox@}}
 \def\tAdjust@{\let\vAdjust@\tAdjust@@}

 \let\vAdjust@\relax

 \def\lAdjust@{\let\hAdjust@\rlap}
 \def\rAdjust@{\let\hAdjust@\llap}

 \let\hAdjust@\relax\let\vAdjust@\relax

 \def\FetchLabel@#1(#2)#3\\{%
     \IN@0#2@@\ifIN@
        \setbox0=\hbox{\ignorespaces#1#3\unskip}%
        \ifdim\wd0>0pt
           \ms@g{}%
           \ms@g{ !!! Bad label(s)? !!!}%
           \message{ #1(#2)#3}%
        \fi
        \def\LabelMole@##1\endFetchLabel@{%
            \IN@0()\\@##1@%
            \ifIN@\def\Temp@{\FetchLabel@##1\endFetchLabel@}%
            \else\def\Temp@{}%
            \fi
            \Temp@
           }%
     \else
       \ignorespaces#1\unskip
       \setbox\auxLabelbox@=%
         \hbox to 0pt{\hss\ignorespaces\hAdjust@
          {\ignorespaces#3\unskip}\hss}%
       \vAdjust@
       \let\hAdjust@\relax\let\vAdjust@\relax
       \AugmentLabelBox@@{#2}%
       \ht\Labelbox@=0pt\dp\Labelbox@=0pt
       \let\LabelMole@\FetchLabel@%
     \fi\LabelMole@}

 \newtoks\XYSep@ 
 \def\SetXYSeparator#1{%
     \IN@0#1@@\ifIN@\XYSep@{*}%
     \else
     \XYSep@{#1}%
     \fi
     }

 \SetXYSeparator*

 \def\AugmentLabelBox@@#1{%
     \IN@0\the\XYSep@ @#1@\ifIN@
       \SPLIT@0\the\XYSep@ @#1@%
       \setbox\Labelbox@=\hbox to 0pt{%
         \unhbox\Labelbox@
         \Shifted@@{\the\Initialtoks@\Wddd@}%
         {\the\Terminaltoks@\Httt@}%
         {\box\auxLabelbox@}}%
     \else
         \ms@g{}%
         \ms@g{ !!! Bad insertion point. !!!}%
         \message{ (#1\ this point was rejected.)}%
     \fi
    }

 \def\FetchOption@#1[#2]#3\endFetchOption@{%
    \def\temp{#1}
    \ifx\temp\empty
       \Edgect@=#2\relax
       \let\EdgGl@b\relax
       \SetEdges@
       \Cleaner@#3%
    \fi}

 \def\Cleaner@#1[@]{\Labeltoks@{#1}}
     
 \def\PlaceLabels@@{\mathsurround=0pt
     \def\Cr@{\\}%
     \let\L\lAdjust@\let\R\rAdjust@
     \let\B\bAdjust@\let\E\eAdjust@\let\T\tAdjust@
     \expandafter\FetchOption@\the\Labeltoks@[@]\endFetchOption@
     \Wddd@=\Wdd@ \Edgescale@\Wddd@ 
     \Httt@=\Htt@ \Edgescale@\Httt@
     \expandafter\FetchLabel@\the\Labeltoks@\endFetchLabel@
     \box\Labelbox@
     }%

 \let \PlaceLabels@\PlaceLabels@@

 \def\AffixLabels#1{\setbox\Coordinatebox@=\hbox{#1}%
      \Wdd@=\wd\Coordinatebox@ \Htt@=\ht\Coordinatebox@
      \advance\Htt@ \dp\Coordinatebox@
      \hbox{\copy\Coordinatebox@\kern-\Wdd@ 
           \Shifted@@{0pt}{-\dp\Coordinatebox@}%
           {\PlaceLabels@\PlaceGrid@}%
           \kern\Wdd@}%
      \GridSwitch@ 
      \LabGl@b\Labeltoks@{()\\}%
      }
 
   \let\wlog\LabelFigwlog@ld   
   \catcode`\@=\LabelFigCatAt  


 
                                By

              Raymond S\'eroul <A18645@FRCCSC21.BITNET>
                                and 
              Laurent Siebenmann <lcs@topo.math.u-psud.fr>
    
              VERSIONS: July 1991, Oct 1991, Jan 1992, July 1992

INTRODUCTION

      This labelling package is intended for TeX users who
rely on non-TeX sources for for their graphics inserts.  It
provides means for adding TeX labels to such inserts with a
minimum of fuss. 

       For most labels, TeX users have in the past found it
reasonably convenient to rely on non-TeX sources. Typical
occasions when an inescapable need for TeX labels seemed to
arise are

 (a) when the graphics program lacks certain exotic or complex
mathematical symbols

 (b) when the very highest typographical quality is wanted for the
labels

 (c) when labels included with the graphics fail to print, 
 and you cannot figure out why (cf. boxedeps.doc).  The labels
 provided by labelfig.tex are 100

       Since this package first appeared, many users, who in the
past scarcely dreamed of using TeX labels, have come to use
nothing but.  So it is now appropriate to add

Intoxication Warning:  TeX labels may be addictive and expensive. 

     If you have a fast preview you may disagree, and even find
that this package provides an agreeable paste-up environment; see
extra applications at end.

     Note to publishers: It is possible and convenient to ultimately
export the TeX labels produced by labelfig.tex to become an integral
part of the EPS file. This is often desired by a publisher who typically
uses an "upmarket" graphics or page layout program, with which the
staff is skilled in perfecting figures.  See Appendix I for
a recipe.

     The authors are grateful to Patrick Ion of Math Reviews for
helpful comments and encouragement.

BASIC INSTRUCTIONS

    After reading in the macro file using

preview or proof your figure with a coordinate grid printed on
top, by typing the following:

    \ShowGrid  
    \AffixLabels{<the graphics insertion>}

Here <the graphics insertion> is what you would type to insert
the graphics object alone without the grid.  This must provide
for the space around it. For example <the graphics insertion>
might well be \BoxedEPSF{MyFigure scaled 700} using the
boxedeps.tex macro package (from same source); this provides a
TeX box containing the encapsulated PostScript insert specified by
the file MyFigure. \AffixLabels{...} provides the grid (supposing
\ShowGrid is present) and later, once you have specified labels
using the grid, it will "tack on" the labels.

     The grid is a sort of (usually elongated) checkerboard of
ten rows and ten columns and its (internal) partitions are by
default numbered  .1, ... ,.9  both horizontally (X-coordinate
running left to right) and vertically (Y-coordinate running bottom
to top).  Thus the points enclosed by the grid correspond to the
points of the unit square in the cartesian "X-Y" plane, the lower
left corner corresponding to the origin (0,0).  By extrapolation,
the full page corresponds to a larger rectangle in the plane.

     These coordinates serve to position labels as follows.
Before the \AffixLabels{...} command type label specifications:

  \SetLabels
   (<X-coordinate>*<Y-coordinate>) <first label> \\
   .
   .
   .
   (<X-coordinate>*<Y-coordinate>)  <last label> \\
  \endSetLabels

Each row specifies one label and is terminated by \\.  In each
row, the position indicator comes first; it is written as a
standard cartesian point except that the X- and Y- coordinates
are separated by * rather than a comma because TeX allows a
comma as decimal point. There are no dimension units to specify
as the unit is the grid itself.

     By default, this cartesian point specifies where the middle
of the baseline of the label will be located.  However if you precede
the point by \L [or \R] the left [or right] edge of the baseline will
be located there. Similarly you may also precede the point by \T, \E,
or \B to vertically align the top equator or bottom of the label box
at the specified point.  This gives nine standard positions of
the label with respect to the insertion point --- corresponding to
the eight principle points of the compas and the center

                     \L\T     \T      \R\T

                     \L\E     \E      \R\E

                     \L\B     \B      \R\B

But this neglects the default "baseline" level of TeX,
giving potentially three more positions

                     \L    <no tag>   \R

For text, the baseline level is often the preferred. Its relation to
the others is variable. It will often coincide with the bottom level,
as happens for "X".  But it is often distinct, as for "g", in which
case you have in all 12 distinct positions rather than 9.

     It is convenient to think of this specification of label
position as attaching the label by a thumb-tack to the coordinate
grid. There are up to twelve positions of the thumb-tack on the
label, while the position of the thumb-tack on the coordinate grid is
arbitrary.  Normally, one choses the position of the thumb-tack on
the label to be the one that is the closest to the item being
labeled.  There are good reasons for this "rule of thumb":

   (a)  It facilitates correct positioning at first try.

   (b)  If the scale of the figure must be altered after labels
have been affixed, the labels have a good chance of remaining well
positioned.

   (c)  The visible grid need not extend beyond the "bounding box"
for the figure, because the best preferred position is always
(at least almost) within the bounding box .

The second reason is particularly important. Indeed it often
happens that scale has to be altered after labelling begins, in
order to either provide space for the labels, or to adjust
proportions between the labels and the figure.  (The size of labels
is unaffected by scaling.)

     Here is an artificial but self-contained test which uses
TeX rules to make a graphics object.

TEST

    Do not skip this!



 \def\FrameIt#1{\hbox{\vrule$\vcenter {\hrule\kern3pt%
             \hbox {\kern3pt #1\kern3pt}%
               \kern3pt\hrule}$\relax\vrule}}

 \def\Caption#1#2{\FrameIt{%
       \vtop {\hsize=#1\relax \parindent=0pt
         \leftskip=0pt \rightskip=0pt plus15pt
         \parfillskip=0pt
         \lineskip=1pt\baselineskip=0pt
         #2}}}

 \def\FirstQuadrant{\hbox to 100pt{\vrule\vbox to 100pt{%
        \hbox to 100pt{\hfil}\vfil\hrule}\hss}}


  \SetLabels
    \R(.5*.2) $\zeta\,\cdot$\\
    (.9*-.10) $\xi$\\
    \R(-.03*.9) $\eta$\\
    \T(.5*.9) \Caption{70pt}{%
          \it The norm of
          $g(\xi+i\eta)$ is indicated on
          contours of this invisible surface.}\\
  \endSetLabels

  \AffixLabels{\FirstQuadrant} \end

  Note that the coordinates to use for labels are indicated on the
edges of the grid (when visible) corresponding to the conventional
x- and y- axes of the Cartesian plane. By default the grid is
1-by-1. However, by the command \Edges{100}, you can change this
to 100-by-100 and many users find this alternative most
convenient. Place the command \Edges{...} in your style file (or
header) since its effect is is global. Other possible edge values
are 10 and 1000.

  If you use the command \Edges{...} at all, do so with care.  For
if you accidentally delete an \Edges{...} command your labels will
abruptly be badly misplaced and may logically but mysteriously
generate "dimension too big" errors under TeX and "off page" errors
under your driver.  

  You can dictate the edgescale for an individual figure by giving
the scale in brackets immediately after \SetLabels.  Thus, to
import into an article using say \Edge{100} a figure labelled using
another edgescale, say the original 1-by-1 default, you can use
\SetLabels[1]...\endSetLabels.


GETTING IT DOWN PAT

     Complicated labeling deserves the same respect as
complicated mathematics.  Do not expect it to come out perfect the
first time!  What is needed in either case is a mechanism to
repeatedly typeset troublesome pieces.

     One mechanism is always available.  One does complicated
labelling in a separate "test" file involving just the figure being
labelled;  a texpert will know how to \dump TeX's current state as
a temporary format that restarts rapidly at each retry.  Usually,
one then pastes the completed labelled figure back into the main
TeX file, but, of course, one can also \input it as an auxiliary
file.

     If you do not have a TeXpert at handy, here is a first
approximation to an efficient setup. By deletions reduce a copy
of your article to just a few lines before and after the figure.
Now label the figure, and finally, copy and paste the labelled
figure to the original article. Then copy the next figure to label
into this testbed and repeat. The TeXpert can improve the  speed
at which TeX starts up, by compiling a format specifically for
your article; just one caution: best NOT include in the format
ephemeral details of setup like \Set<mydriver>ArtSpecials (from
boxedeps.tex because this reads  figure dimensions which you may
change during your work session.

     An improved mechanism to repeatedly typeset troublesome
pieces is now available on the Macintosh; it is called LinoTeX;
see the same ftp sources.  It could be set up on many types
of computer.

     Before using labelfig.tex to attach labels to a graphics
object inserted using boxedeps.tex or BoxedArt.tex, make it a
firm rule to carefully adjust the bounding box using the trimming
commands of these packages, and also at least tentatively scale
and position the object. Beware of changing the grid inadvertently
after the labels have been positioned.  For example, correcting
the bounding box of a PostScript graphics object can foul up the
labels by changing the coordinate grid to which the labels are
attached. This is particularly true for the trimming  commands of
boxedeps.tex and BoxedArt.tex. However, as noted already, change
of scale is much less disruptive, and modest adjustments should be
well tolerated.

     Sometimes the labels protrude so far from the bounding box
of a figure that the figure has to be repositioned.  Best do this
by ad hoc spacing, say using \hglue and \vglue; altering the
bounding box would create a vicious circle.

     Remember that you are responsible for preventing labels
from overlapping. You are responsible for all label typography
including size and style. A label is really just about anything
that can be put in a TeX box. Note that spaces at the beginning
and end of labels will normally be suppressed; if you really want
them you must protect them with TeX braces.

     This package temporarily sets the \mathsurround parameter
of TeX to zero  while the labels are being affixed. This is done
because nonzero \mathsurround space would influence the position
of left and right aligned labels; then, when a texpert or printer
modifies mathsurround, diagram labeling might be disastrously
altered. There is a small price to pay involving labels that are
formatted as caption boxes including mathematics: you  may want or
need to specify an explicit mathsurround space within the caption
box; it will not influence anything outside.

     Those hostile to the use of * as separator between
the X and Y coordinates of label insertion points, are free to
impose another using \SetXYSeparator{<the new separator>}.  
Americans may prefer "," to "*" since they never use a 
comma as a decimal point; on the other hand, * may be more visible.

APPENDIX (I)  MERGING labelfig.tex LABELS INTO AN EPSF GRAPHICS OBJECT.

     As promised in the introduction, here is a recipe useful for
publishers. It works at least on Macintosh and at least for vectorized
graphics and Adobe type1 fonts.  (There is surely a similar recipe for
PCs under MSWindows.)

 (a)  Use boxedeps.tex utility to integrate the figure given by the eps
file, "x.eps" say, with a visible frame around it.  See
\ShowDisplacementBoxes command in boxedeps.tex.  To get precise results
automatically it is important to use the \Trim... commands of
boxedeps.tex making the "DisplacementBox" neatly fit the figure.

 (b)  Use the TeX printer driver and LaserWriter (versions >= 8.1.1) to
export to an EPSF the DVI page containing the integrated, labelled
figure. You now have an EPS file  "xx.eps"  that contains too much, and at
the wrong scale, and at wrong position.

 (c)  Convert the EPSF to an Adode Illustrator format EPSF using
the shareware utility called epsConvert by Sam Weiss
1993-- (currently $25).

 (d)  In Illustrator (or a compatible program), group the labels and the
"DisplacementBox"; copy them to the clipboard and paste them into "x.ps".
This step requires that all the label fonts be "visible to the Macintosh.

 (e)  Translate and scale the pasted group consisting of the labels plus
the "DisplacementBox" so as to make the "DisplacementBox" the bounding
box of (labelless) figure represented by "x.eps".  At this point the
labels will be correctly placed on the figure "x.eps".

 (f)  Ungroup and delete the "DisplacementBox".  The result is the
desired single EPS file, "x+.eps" say, It contains the original figure
plus its labels.  

     Using grouping and ungrouping appropriately in "x+.eps", a
publisher's staff can very efficiently improve label positions etc.

APPENDIX II)  SOME EXOTIC APPLICATIONS

     The grid of labelfig.tex is analogous to a light-table in
classical page makeup with wax or latex glue.  In principle, you
can use it to compose any page from its indivisible parts.  This
even has some of the artisanal charm of classical paste-up
provided you have a fast screen preview to make the process
"interactive".

     In practice labelfig.tex is a tool for nonstandard jobs.
Here are a few going beyond the labelling already discussed.

(I)  GRAPHICS INTEGRATION.

     This is accomplished by treating the imported graphics
objects as labels.  The underlying graphics object is then
typically an empty  \vbox to <dimension>{\vfill} in a TeX
\midinsert...\endinsert construction.  A label line
might be of the form

   (.1*.1) \special{... MyFigure ...}\\

The exact form of the special command varies from driver to
driver.  However, in the case of encapsulated PostScript graphics
(EPSF norm), by relying on boxedeps.tex, one can have the
following standard syntax (independant of driver  (see
boxedeps.doc for details.
  
  (.1*.1) \BoxedEPSF{MyFigure scaled <scale in mils>}\\

This may be slow since it requires TeX to read the PostScript
file to read bounding box using many complex macros.  So you
may want to try

  (.1*.1) \EPSFSpecial{MyFigure}{<scale in mils>}\\

which is fast and driver independant, but it squashes the
bounding box, normally to its lower left corner.

     Similarly for graphics of the Macintosh PICT norm ---
using BoxedArt.tex (same sources) in place of boxedeps.tex.

     This approach to integration is to be recommended when
one is assembling a composite graphics object.

 (II)  COMMUTATIVE DIAGRAM ENHANCEMENT

     Commutative diagrams or arrays of mathematical objects
connected by arrows of various sorts are common in mathematics.
The mathematical objects require the use of TeX.  Recently TeX
acquired a good collection of arrows of all slopes --- that of
LamSTeX --- plus pwerful macros to build the diagrams.

     However, even the LamSTeX collection is often
inadequate; it lacks for example double shafted arrows, dotted
arrows and curved arrows. Fortunately it is possible to produce
such arrows on an individual basis using sophisticated graphics
programs such as Illustrator and AldusFreehand (both serving
the EPSF norm) or using Metafont (with its public domain norm).
Since the creation of each new arrow is a work of love, you
probably want to limit the number of arrows by using LamSTeX
for most arrows. The 40K commutative diagram module of LamSTeX
has been adapted to work with AmSTeX and a copy may be posted
with LabelFig and related files. Unfortunately no one has yet
offered a version that works with Plain TeX or LaTeX.

       Suffice it here to say that when the exotic arrow has
been somehow imported into TeX, labelfig.tex treats it as a
label that one affixes to the commutative diagram.  Two other
steps will be treated in separate notes, namely the matter of
extracting the dimension specifications for the arrow and the
construction of the arrow --- for these steps are far from
unique and often depend intimately on your computer environment. 
Notes for the Macintosh-Textures-Illustrator combination are
found in the file ExoticArrows.doc.

 (III) NESTING 

Ingenuity pays off in exploiting labelfig.tex. One can
mix graphics and typography quite freely.  labelfig.tex is good
for freeform or overlapping arrangements, while boxedeps.tex (or
BoxedArt.tex) is best for regimented non-overlapping
arrangements --- and the two can be combined.

     The default behavior of labelfig.tex is not ideal 
for nesting objects, because to prevent trouble for beginners
the register for labels is globally cleared when \AffixLabels
concludes.  But there are switches available

      \LabelsGlobal      \LabelsLocal

which change this.  To understand this, extend the above test 
by something like:


 \LabelsLocal

 \SetLabels
    (.5*.5) AAA\\
 \endSetLabels

 {
 \SetLabels
    (.5*.5) ZZZ\\
 \endSetLabels
   \AffixLabels{\FirstQuadrant}
 }

   \AffixLabels{\FirstQuadrant}


     There are however potential pitfalls.  Neither
labelfig.tex nor boxedeps.tex has been tested under extreme
conditions. Problems may occur if their procedures are
indiscriminately nested. For boxedeps.tex (not labelfig.tex)
there is a precise cause for worry, namely many of its
variables are "global", which means that TeX braces will not
provide the protection one might expect.

COMMAND SUMMARY FOR labelfig.tex

  Here [...] means optional (one or zero)
       [...]* means any number of such constructs

  \SetLabels
    [[<P>](<X><Sep><Y>) <label> \\]*
  \endSetLabels
  \ShowGrid  
  \AffixLabels{<the figure>}

   --- <P> is tack position, one of eleven or empty
              order irrelevant

                   \L\T      \T      \R\T

                   \L\E      \E      \R\E

                     \L               \R

                   \L\B      \B      \R\B

   --- (<X><Sep><Y>) insertion point;
  <Sep> is separator, = * by default;
  \SetXYSeparator{<Sep>} changes it.
   <X> and <Y> are real numbers

  --- <label> a label to attach 

  --- <the figure> the figure to label 

  \GlobalLabels (default)     
  \LocalLabels  setting for nested constructs.

 \Grids makes ALL grids appear; \HideGrid then makes just next disappear.
 \noGrids returns to default.  The commands are always global.

 \GridLineWidth{<dimension>} adjusts width of grid lines. Default is very
small, to give "hairline" effect. If your grid lines are missing try
setting \GridLineWidth{1pt}.

 \Edges#1 globally changes the edge size of all grids to the numerical 
value #1, which must be 1, 10, 100, or 1000.  The default is 1.

VERSION HISTORY.
 --- Jan 1993: \Edges#1 and [??] option after \SetLabels
 --- July 1992: \Grids, \noGrids, \HideGrid;
       Gridlines become hairlines; \GridLineWidth{<dimension>}.
 --- Oct 1991, Jan 1992: \SetXYSeparator{<Sep>},  \LabelsGlobal,
       \LabelsLocal.
 --- July 1991: first release

Address for bugs and other feedback:

        Raymond S\'eroul
        IREM and Lab. de Typographie Informatise
        Univ. Rene Descartes
        Strasbourg

    Tel 33-88-41-63-45
    Email:  A18645@FRCCSC21.BITNET

        Laurent Siebenmann
        Mathematique, Bat. 425,
        Univ de Paris-Sud,
        91405-Orsay,
        France

    Tel 33-1-6941-7949; 
    Email: lcs@topo.math.u-psud.fr  

\chardef\newinsCatAt\the\catcode `\@
\catcode `\@=11
%
%
%
\newskip\insertskipamount\newskip\inserthardskipamount
\insertskipamount 12pt plus2pt     
\inserthardskipamount 4pt          
\def\insertskip{\vskip\insertskipamount}
%
%
\newskip\LastSkip
\def\SaveLastSkip{\LastSkip\lastskip}
\def\RestoreLastSkip{\nobreak\vskip-\LastSkip\vskip\LastSkip}
%
%
\newcount\SplitTest
\def\SetSplitTest{\SplitTest\insertpenalties
  \insert\topins{\floatingpenalty1}%
  \advance\SplitTest-\insertpenalties}
%
%
\def\midinsert{\par
 \SaveLastSkip\penalty-150\SetSplitTest\RestoreLastSkip
 \ifnum\SplitTest=-1
  \@midfalse\p@gefalse\else\@midtrue\fi\@ins}
\def\@ins{\par\begingroup\setbox\z@\vbox\bgroup%
  \vglue\inserthardskipamount}
\def\endinsert{\egroup 
  \if@mid \dimen@\ht\z@ \advance\dimen@\dp\z@
    \advance\dimen@\insertskipamount
    \advance\dimen@\pagetotal\advance\dimen@-\pageshrink
    \ifdim\dimen@>\pagegoal\@midfalse\p@gefalse\fi\fi
  \if@mid%
    \ifdim\lastskip<\insertskipamount\removelastskip\insertskip\fi
    \nointerlineskip\box\z@\penalty-200\insertskip
  \else%
    \SaveLastSkip
    \insert\topins{\penalty100 
    \splittopskip\z@skip
    \splitmaxdepth\maxdimen \floatingpenalty\z@
    \ifp@ge \dimen@\dp\z@
    \vbox to\vsize{\unvbox\z@\kern-\dimen@}
    \else \box\z@\nobreak\insertskip\fi}
    \RestoreLastSkip
   \fi\endgroup}
%
\catcode `\@=\newinsCatAt

\checkdefinedreferencetrue
\theoremcountingtrue
\sectionnumberstrue
\figuresectionnumberstrue
\forwardreferencetrue
\tocgenerationtrue
\citationgenerationtrue
\hyperstrue
\initialeqmacro


\long\def\comment#1{\relax}

\newcount\ftntno
\ftntno=0
\def\Ftnote#1{\global\advance\ftntno by 1\ftnote{{$^{\the\ftntno}$}}{#1}}
\def\subsection#1{\bigskip\noindent {\bf #1.}\par\nobreak}

\def\itemrm#1{\item{{\rm #1}}}
\def\beginitems{\begingroup\parindent=25pt}
\def\enditems{\vskip0pt\endgroup\noindent}



\def\U{{\Bbb U}}
\let\SSS\S
\def\SS{{\Bbb S}}
\def\C{{\Bbb C}}
\def\H{{\Bbb H}}
\def\Np{{\N_+}}

\def\Sp#1{{\SS^{#1}}}
\def\SpD{{\SS^2_D}}

\def\verts{{\ss V}}

\def\edges{{\ss E}}

\def\LE{{\ss LE}}
\def\RW{{\ss RW}}

\def\B{{\ss B}}
\def\Tr{{\ss trunk}}     
\def\WT{{{\ss WT}}}
\def\FT{{\ss FT}}

\def\T{{\ss T}}


\def\dHa{d_{\Ha}} 
\def\dd_#1{d_{\Ha(#1)}}
\def\dlr{\hat d} 
\def\dlru{\dlr_{\U}} 
\def\dHc#1{d_{\Ha(#1)}} 
\def\toP{{\,\mathop{\to}\limits^{\P}\,}}
\def\dsp{d_{sp}} 
\def\ppn#1{\left|#1\right|_{2\pi}}


\def\mea{\mu}
\def\mmu{\P}
\def\law{\mea}
\def\dlaw{\law_\delta}
\def\mun{\P_\delta}


\def\Peano{\theta}
\def\path{\omega}
\def\lep{\xi}
\def\lsl{\sigma}
\def\sle{\sigma}
\def\Pa{\omega}
\def\Paa(#1,#2){\Pa_{#1,#2}}
\def\Paq{\Pa^q}
\def\ray{\QQ}
\def\rayz{\ray^0}


\def\Gd{G_\delta(\H)}
\def\hGd{\hat G_\delta(\H)}
\def\G{G}

\def\d#1{#1^\dagger}    
\def\dZ{\delta\Z^2}
\def\Vi{{\verts_I}}
\def\Vb{{\verts_{\bd}}}
\def\Vba{{\verts_{\bd}^0}}
\def\zerp{{o_1}}
\def\onep{{o_2}}

\def\zer{{o_1^*}}
\def\one{{o_2^*}}
\def\aa{\zerp}
\def\aap{\aa}
\def\bb{\onep}


\def\lec{\kappa}
\def\lecv{\kappa}
\def\lecl{{\kappa_o}}
\def\crit{4}
\def\xplec{2}


\def\eps{\epsilon}
\def\Cal#1{{\cal #1}}

\def\leq{\leqslant}

\def\geq{\geqslant}
\def\Ha{{\Cal H}}
\def\br{\overline}

\def\KK{\frak K}
\def\dist{{\rm dist}}

\def\diam{\operatorname{diam}}
\def\hatzet{\hat\zeta}
\def\winding{W}

\def\pdt{{\partial_t}}
\def\ev#1{{\Cal #1}}
\def\bd{\partial}
\def\st{:\,}
\def\operatorname#1{\mathop{{\rm #1}}}
\def\Rsp{{\An_{\rm sp}}}
\def\QQ{Q}
\def\An{A}
\def\im{\operatorname{Im}}

\def\hat{\widehat}
\def\ST{{\frak T}} 
\def\WST{{\frak W\frak T}} 
\def\FST{{\frak F\frak T}} 
\def\hT{{\hat \T}}

\def\lowner{L\"owner}

\ifproofmode \relax \else\head{} {\versiondate}\fi
\vglue20pt

\title{Scaling limits of loop-erased random walks}
\title{and uniform spanning trees}

\author{Oded Schramm}

\abstract{
The uniform spanning tree (UST) and the loop-erased
random walk (LERW) are strongly related probabilistic processes.
We consider the limits of these models on a fine grid in the plane,
as the mesh goes to zero. 
Although the existence of scaling limits is still unproven,
subsequential scaling limits can be defined in various ways,
and do exist.  We establish some basic a.s.\ properties of these
subsequential scaling limits in the plane. It is proved
that any LERW subsequential scaling limit is a simple
path, and that the trunk of any UST subsequential scaling limit 
is a topological tree, which is
dense in the plane.

The scaling limits of these processes are conjectured to be conformally
invariant in dimension 2.  We make a precise 
statement of the conformal invariance conjecture for the LERW,
and show that this conjecture implies
an explicit construction of the scaling limit, as follows.
Consider the \lowner\ differential equation
$$
{\partial f\over\partial t} 
= z \, {{\zeta(t)+z}\over{\zeta(t)-z}} \, {\partial f\over\partial z}
\,,
$$
with boundary values $f(z,0)=z$, in the range $z\in\U=\{w\in\C\st |w|<1\}$,
$t\leq 0$.
We choose $\zeta(t):= \B(-2t)$, where $\B(t)$ is Brownian motion
on $\partial \U$ starting at a random-uniform point in $\partial \U$.
Assuming the conformal invariance of the LERW scaling limit
in the plane, we prove that the scaling limit
of LERW from $0$ to $\partial\U$ has the same
law as that of the path $f(\zeta(t),t)$
(where $f(z,t)$ is extended continuously to $\partial\U\times (-\infty,0]$).
We believe that a variation of this process gives the scaling
limit of the boundary of macroscopic critical percolation clusters.
}

\bottom{Primary
60D05
. %
Secondary
60J15
.} 
{\lowner{} differential equation, topological trees,
loop-erased Brownian motion, conformal invariance.}
{Research supported by the Sam and Ayala Zacks Professorial Chair.}

\vfill\eject

\articletoc

\midinsert
\centerline{\epsfysize=2.4in\epsfbox{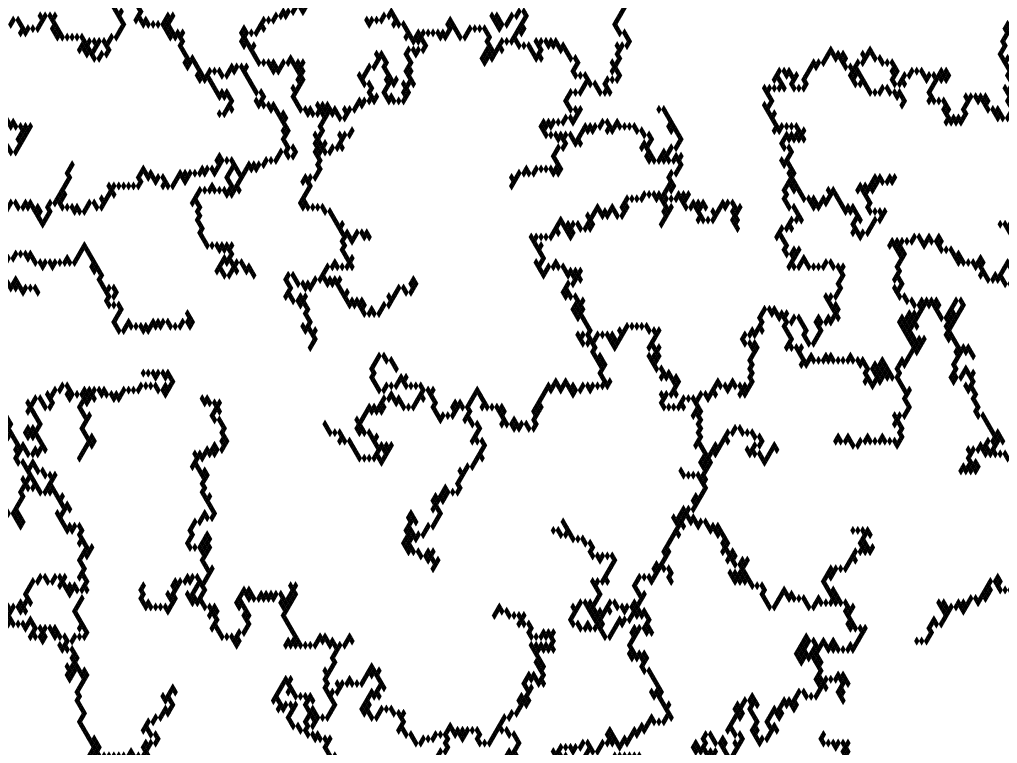}}
\caption{
An approximation of the UST scaling limit trunk.}
\endinsert

\bsection {Introduction}{s.intro}

\subsection{General remarks about scaling limits}
It is often the case that grid-based probabilistic models should
be considered as a mere substitute, or simplification, of a continuous
process.  There are definite advantages for working in the discrete
setting, where unpleasant technicalities can frequently be avoided,
simulations are possible, and the setup is easier to comprehend. 
On the other hand, one is often required to pay some price for the 
simplification.
When we adopt the grid-based world, we sacrifice rotational
of conformal symmetries which the continuous model may enjoy,
and often have to accept some arbitrariness in the formulation
of the model.  There are also numerous examples where the continuous process
is easier to analyze than the discrete process, and in such
situations the continuous may be a useful simplification
of the discrete.

Understanding the connections between grid-based models and
continuous processes is a project of fundamental importance,
and so far has only limited success.  
As mathematicians, we should not content ourselves with the vague 
notion that the
discrete and continuous models behave ``essentially the same'',
but strive to make the relations concrete and precise.


One reasonable way to define a continuous process, is by taking
a {\bf scaling limit} of a grid process.  This means making
sense of the limit of a sequence of grid processes on finer
and finer grids.  Recently, Aizenman \ref b.Aizenman:web/ 
has proposed a definition for the scaling limit of percolation,
and we shall propose a somewhat different definition \ref b.Schramm:inprep/.

Although, in general, the understanding of the connections
between grid-based models and continuous models is lacking,
there have been some successes.  The classical and archetypical
example is the relation between simple random walk (SRW) and
Brownian motion, which is well studied and quite well understood.
See, for example, the discussion of Donsker's Theorem in
\ref b.Durrett:probability/.
We also mention that recently, T\'oth and Werner \ref b.TW/ 
have described the scaling limit of a certain self-repelling walk on $\Z$.

The present paper deals with the scaling limits of two very
closely related processes, the loop-erased random walk (LERW)
and the uniform spanning tree (UST). 
While these processes are interesting
also in dimensions $3$ and higher, we restrict attention to two
dimensions.
In the plane, the scaling limits are conjectured to possess
conformal invariance (precise statements appear below),
and this can serve as one justification for the special interest
in dimension $2$. The recent preprint by Aizenman, Burchard,
Newman and Wilson \ref b.ABNW/ discusses scaling limits of
random tree processes in two dimensions, including the UST.
The present work answers some of the questions left open
in \ref b.ABNW/.

The most fundamental task in studying a scaling limit process is to set up
a conceptual foundation for the scaling limit.  This means answering
the following two questions: what kind of object is the scaling limit,
and what does it mean to be the scaling limit?  For the first question,
there's often more than one \lq\lq right\rq\rq{} answer.  For example,
the scaling limit of (two-sided) simple random walk in $\Z$ is usually 
defined as a probability measure on $C(\R)$, the space of continuous,
real-valued functions on $\R$, but it could also be defined
as a measure on the space of closed subsets of $\R\times\R$, with
an appropriately chosen metric.  There's also a \lq\lq wrong\rq\rq{}
answer here.  It is not a good idea to consider the scaling
limit of SRW as a probability measure on $\R^\R$, though this might
seem at first as more natural.

After the conceptual framework is fixed, the next natural question
is the existence of the scaling limit.  Unfortunately, we cannot
report on any progress here.  There is every reason to believe that
UST and LERW have scaling limits, but a proof is still lacking.
However, in the setup we propose below, the existence of subsequential
scaling limits is almost a triviality:
for every sequence of positive $\delta_j$ tending to $0$, there is
a subsequence $\delta_{j_n}$ such that the UST and the LERW
on the grids $\delta_{j_n}\Z^2$ do converge to a limit as
$n\to\infty$.  Such a limit is called a subsequential
scaling limit of the model, and is a probability measure on some space.

All the above discussion concerns foundational issues, which are important. 
But it is not  less important to prove
properties of the (subsequential) scaling limit%
\Ftnote{Sometimes, this can be done without even defining the scaling limit.
For example, the Russo-Seymour-Welsh Theorem
\ref b.Russo/, \ref b.SW/ in percolation theory implies properties of any
reasonably defined percolation scaling limit.  Similarly, 
Benjamini's preprint \ref b.Benjamini:ustsl/ does not explicitly discuss
the scaling limit of UST, but has implication to the UST scaling limit.}.
In this paper, we prove several almost sure properties of the UST and
LERW subsequential scaling limit.  
We now describe these models and explain the results.

\subsection{The LERW model, and its scaling limit}
Consider some set of vertices $K\neq\emptyset$, in a recurrent graph $G$,
and a vertex $v_0$.
The LERW from $v_0$ to $K$ in $G$ is obtained by running simple random
walk (SRW) from $v_0$, erasing loops as they are created, and stopping
when $K$ is hit.
Here's a more precise description.
Let $\RW$ be simple random walk starting at $\RW(0)=v_0$ and stopped at the
first time $\tau$ such that $\RW(\tau)\in K$.
Its loop-erasure, $\LE=\LE_{\RW}$, is defined inductively as follows:
$\LE(0):=v_0$, and $\LE(j+1)=\RW(t+1)$ if $t$ is the last time
less than $\tau$ such that $\RW(t)=\LE(j)$.  The walk $\LE$
stops when it gets to $\RW(\tau)\in K$.
Note that this LERW is a random simple path from $v_0$ to $K$.

On a transient graph, it may happen that $\RW$ does not hit $K$.
However, one can discuss the loop-erasure of the walk continued
indefinitely, since it a.s.\ visits any vertex only finitely many
times.

LERW on $\Z^d$ was studied extensively by Greg Lawler (see
the survey paper \ref b.Lawler:survey/ and the references therein),  
who considered LERW as a simpler substitute for the self-avoiding random
walk (see the survey \ref b.Slade:survey/),
which is harder to analyze.  However, we believe
that the LERW model is just as interesting mathematically, because
of its strong ties with SRW and UST.

For compactness's sake, in the following we consider the plane $\C=\R^2$
as a subset of the two-sphere, $\Sp2=\C\cup\{\infty\}$,
which is the one point compactification
of the plane, and work with the spherical metric $\dsp$ on $\Sp2$.
Let $D$ be a domain (nonempty open connected set) in the plane $\C=\R^2$.
We consider a graph $\G=\G(D,\delta)$,
which is an approximation of the domain $D$ in the square grid $\dZ$ of
mesh $\delta$. 
The {\bf interior vertices}, $\Vi(\G)$, of $\G$ are the vertices of
$\dZ$ which are in $D$, and the {\bf boundary vertices}, $\Vb(\G)$,
are the intersections of edges of $\dZ$ with $\bd D$.
(The precise definition of $\G$ appears in \ref s.back/.)
Suppose that each component of $\bd D$ has
positive diameter.
Let $a\in D$, and let $\LE=\LE_{a,D,\delta}$ be 
LERW from a vertex $a'\in\dZ\cap D$ closest to $a$ to $\Vb(\G)$ in $\G$. 

To make sense of the concept of the {\bf scaling limit} of LERW in
$D$, we think of $\LE$ as a random set in $\br D$.
Recall that the Hausdorff distance $\dHa(X,Y)$ between two closed
nonempty sets $X$ and $Y$ in a compact metric space $Z$ is the
least $t\geq 0$ such
that each point $x\in X$ is within distance $t$ from $Y$ and each point
$y\in Y$ is within distance $t$ from $X$; that is,
$\dHa(X,Y)=\inf\left\{t\geq 0\st X\subset\bigcup_{x\in X} B(x,t),\
Y\subset\bigcup_{y\in Y} B(y,t)\right\}$.
On the collection $\Ha(D)$ of closed subsets of $D$, we use the
metric $\dHc D(X,Y):=\dHa(X\cup\bd D,Y\cup\bd D)$, and
$\Ha(D)$ is compact with this metric.
Then $\LE\cap D$ is a random element in $\Ha(D)$, and its distribution
$\mea_\delta=\mea_{\delta,D}$ is a probability measure on $\Ha(D)$.
Because the space of Borel probability measures on a compact space is
compact in the weak topology%
\Ftnote{We review the notion of weak convergence in \ref s.back/.}%
, there is a sequence $\delta_j\to 0$ 
such that the weak limit $\mea_0:=\lim_{j\to\infty}\mea_{\delta_j}$
exists.  Such a measure $\mea_0$ will be called
a {\bf subsequential scaling limit} measure of LERW from $a$ to $\bd D$.
If $\mea_0=\lim_{\delta\to0}\mea_\delta$, then we say
that $\mea_0$ is the {\bf scaling limit} measure of LERW from $a$ to $\bd D$.

Similarly, we may consider the scaling limit of LERW
between two distinct points $a,b\in\Sp2$, as follows.  For $\delta>0$,
we take $\LE$ to be the loop-erasure of SRW on $\dZ$ starting from a vertex
of $\dZ$ within distance $2\delta$ of $a$ and stopped when
it first hits a vertex within distance $2\delta$ of $b$.
Since $\LE$ is a.s.\ compact,
its distribution is an element of the Hausdorff space $\Ha(\Sp2)$,
and there exists a Borel probability measure on $\Ha(\Sp2)$,
which is a subsequential scaling limit measure of the law of $\LE$.

\procl t.noloop
Let $D$ be a domain in $\Sp2$ such that each
connected component of $\bd D$ has positive diameter, and let 
$a\in D$.  Then every subsequential scaling limit measure
of LERW from $a$ to $\bd D$ is supported on simple paths.

Similarly, if $a,b$ are distinct points in $\Sp2$,
then every subsequential scaling limit of the LERW from $a$ to $b$
in $\dZ$ is supported on simple paths.
\endprocl

Saying that the measure is supported on simple paths means that
there's a collection of simple paths whose complement has zero measure.

\subsection{The conformal invariance conjecture for LERW}
Consider two domains $D,D'\subset \Sp2$.  Every homeomorphism $f:D\to D'$ 
induces a homeomorphism $\Ha(D)\mapsto\Ha(D')$.  Consequently,
if $\mea$ is a probability measure on $\Ha(D)$, there is
an induced probability measure $f_*\mea$ on $\Ha(D')$.

\procl g.confinv
Let $D\subsetneqq\C$ be a simply connected domain in $\C$, and let
$a\in D$. 
Then the scaling limit of LERW from $a$ to $\bd D$ exists.
Moreover, suppose that
$f:D\to D'$ is a conformal homeomorphism onto a domain $D'\subset\C$.
Then $f_*\mea_{a,D}=\mea_{f(a),D'}$, where $\mea_{a,D}$ 
is the scaling limit measure of LERW from $a$ to $\bd D$,
and $\mea_{f(a),D'}$ is the scaling limit measure of LERW from $f(a)$
to $\bd D'$.
\endprocl

Although conformal invariance conjectures have been \lq\lq floating in the
air\rq\rq{} in the physics literature for quite some time now,
we believe that this precise statement has not yet appeared explicitly.
Support for this conjecture comes from simulations which we have
performed, and from the work of Rick Kenyon~\ref b.Kenyon:conf/,
\ref b.Kenyon:5o4/, \ref b.Kenyon:longrange/.

We prove that \ref g.confinv/ implies an explicit description of
the LERW scaling limit
in terms of solutions of \lowner{}'s differential equation
with a Brownian motion parameter.  We now give a brief explanation
of this. 

\midinsert
\SetLabels
\R(0.49*0.5)$0$\\
\T\R(0.55*0.4)$\beta$\\
\L\B(0.75*0.36)$q$\\
\T\R(0.79*0.24)$\beta_q$\\
\endSetLabels
\centerline{\AffixLabels{\epsfysize=2.4in\epsfbox{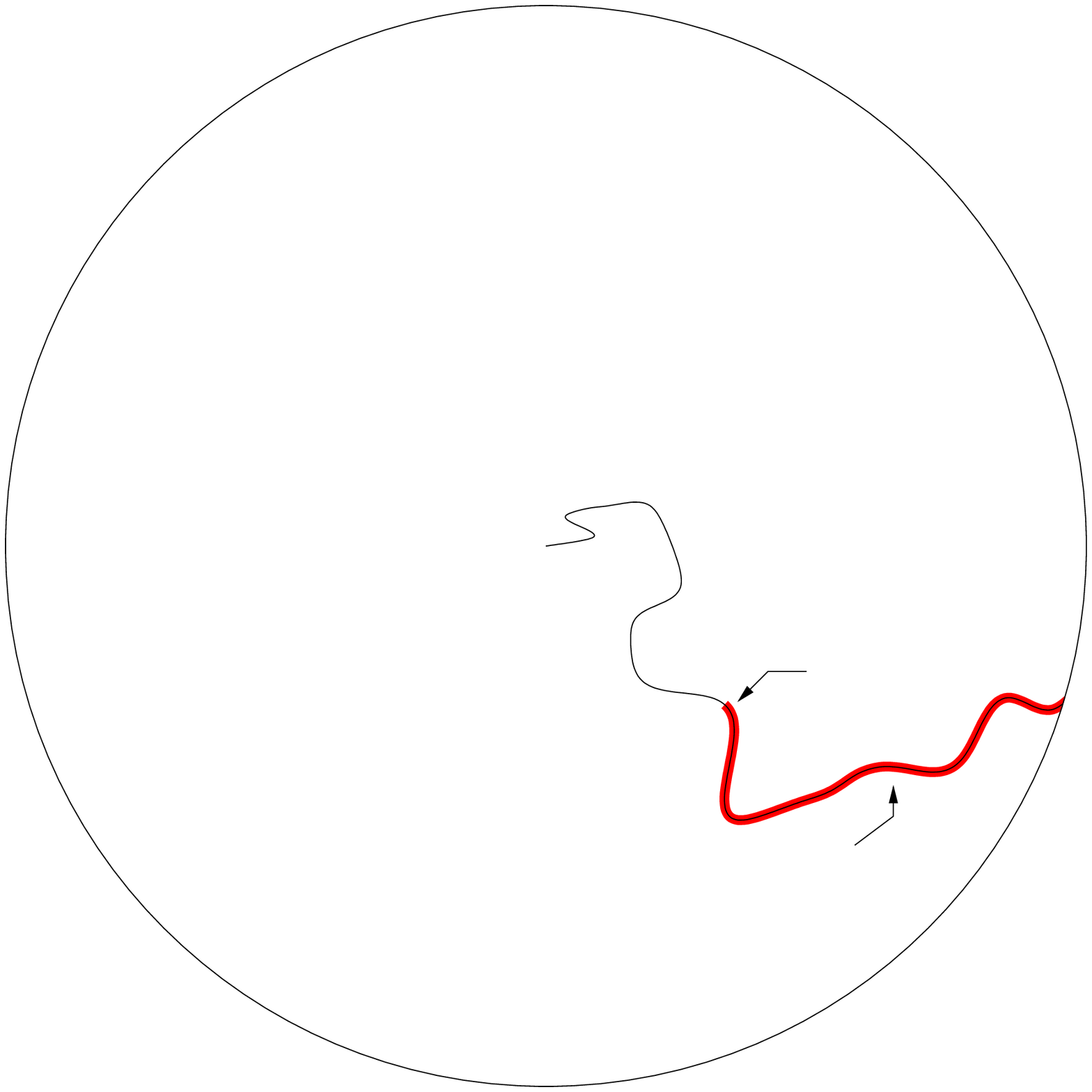}}}
\caption{\figlabel{gamma}\enspace}
\endinsert

Let $\U:=\{z\in\C\st |z|<1\}$, the unit disk.
If $\gamma$ is a compact simple path in $\br U-\{0\}$,
such that $\gamma\cap\bd \U$ is an endpoint of $\gamma$,
then there is a unique
conformal homeomorphism $f_\gamma:\U\to\U-\gamma$ such that $f_\gamma(0)=0$
and $f_\gamma'(0)>0$ (that is, $f_\gamma'(0)$ is real and positive).
Moreover, if $\hat \gamma$ is another such path, and
$\hat\gamma\supset\gamma$, then $f_\gamma'(0)>f_{\hat\gamma}'(0)$.
Now suppose that $\beta$ is a compact simple path in $\br\U$ such
that $\bd\U\cap\beta$ is an endpoint of $\beta$ and $0$ is the other
endpoint of $\beta$, as in \ref f.gamma/.  For each point $q\in\beta-\{0\}$,
let $\beta_q$ be the arc of $\beta$ extending from $q$ to $\bd\U$,
and let $h(q):=\log f'_{\beta_q}(0)$.
(If $q$ is the endpoint of $\beta$ on $\bd\U$, then
$f_{\beta_q}(z)=z$.)  It turns out that
$h$ is a homeomorphism from $\beta-\{0\}$ onto $(-\infty,0]$.
We let $q(t)$ denote the inverse map $q:(-\infty,0]\to\beta$,
and set $f(z,t)=f_t(z):=f_{\beta_{q(t)}}(z)$.
In this setting, \lowner{}'s Slit Mapping Theorem \ref b.Lowner/
(see also \ref b.Pommerenke:Lowner/)
states that $f_t(z)$ is the solution of \lowner{}'s equation
$$
{\bd\over\bd t}f_t(z) = 
z f_t'(z){\zeta(t)+z\over\zeta(t)-z}
\,,\qquad\forall z\in\U\,,\forall t\in(-\infty,0]
\,,
\label e.loew
$$
where 
$\zeta:(-\infty,0]\to\bd\U$ is some continuous function.
In fact, $\zeta(t)$ is defined by the equation
$$
f\bigl(\zeta(t),t\bigr)=q(t)\,.
$$
(The left hand side makes sense, since there is
a unique continuous extension of $f_t$ to $\bd\U$.)
Note that $f$ also satisfies
$$
f(z,0)=z\,,\qquad\forall z\in\U
\,.
\label e.bdval
$$ 

\procl t.lesl \procname{The differential equation for the LERW scaling limit}
Assume \ref g.confinv/. 
Let $\B(t)$, $t\geq 0$ be Brownian motion on $\bd\U$ starting
from a uniform-random point on $\bd\U$.
Let $f(z,t)$ be the solution of \ref e.loew/ and \ref e.bdval/, with
$\zeta(t):=\B(-\xplec t)$.
Set
$$
\sle(t)=f(\zeta(t),t)\,,\qquad t\leq 0
\,.
\label e.sledef
$$
Then $\{0\}\cup\sle\bigl((-\infty,0)\bigr)$ has the same distribution
as the scaling limit of LERW from $0$ to $\bd\U$.
\endprocl

The Brownian motion $\B(t)$ in the theorem can be defined as
$\B(t):=\exp(i\hat \B(t))$, where $\hat\B(t)$ is ordinary
Brownian motion on $\R$, starting at a uniform-random point
in $[0,2\pi)$.

\procl r.ode 
Although \ref e.loew/ may look like a PDE, it can in fact
be presented as an ODE. 
Set $\Phi(z,t,s) :=  f_s^{-1}(f_t(z))$ when $t\leq s\leq 0$.
Then 
\begineqalno
\Phi(z,t,0)&=f_t(z)\,,\quad
\cr
\Phi(z,t,t)&=z
\,.
\label e.phis
\endeqalno
It is immediate to see that
$f_t$ satisfies \ref e.loew/ iff $\Phi=\Phi(z,t,s)$ satisfies
$$
{\bd\Phi \over\bd s}=
-\Phi {\zeta(s)+\Phi\over\zeta(s)-\Phi}
\,.
\label e.phiode
$$ 
Therefore, $f_t(z)$ can obtained by solving the ODE \ref e.phiode/
with $t$ fixed and $s\in[t,0]$.

Note that $(\zeta+\Phi)/(\zeta-\Phi)$ has positive real part
when $\Phi\in\U$ and $\zeta\in\partial\U$.
Therefore, \ref e.phiode/ implies that $|\Phi(z,t,s)|$ is monotone
decreasing as a function of $s$.  {}From this it can be deduced
that there is a unique solution to the
system \ref e.phis/ and \ref e.phiode/ in the interval $s\in[t,0]$.
\endprocl

Obviously, \ref t.lesl/ together with \ref g.confinv/ describe
the LERW scaling limit in any simply connected domain $D\subsetneqq\C$,
since such domains are conformally equivalent to $\U$.

It would be interesting to extract properties of the LERW scaling
limit from \ref t.lesl/.  

At the heart of the proof of \ref t.lesl/ lies the following simple
combinatorial fact about LERW.  Conditioned on a subarc $\beta'$ of
the LERW $\beta$ from $0$ to $\bd D$, which extends from 
some point $q\in\beta$ to $\bd D$, the distribution of $\beta-\beta'$
is the same as that of LERW from $0$ to $\bd (D-\beta')$, conditioned
to hit $q$. (See \ref l.dprod/.)
When we take the scaling limit of this property, and apply the
conformal map from $\bd (D-\beta')$ to $\U$, this translates
into the Markov property and stationarity of the associated \lowner{}
parameter $\zeta$.

\subsection{The uniform spanning tree and its scaling limit}

Shortly, a definition of the uniform spanning tree (UST) on $\Z^2$ 
will be given.  The UST is a statistical-physics model. It lies
in the boundary of the two-parameter family of random-cluster measures,
which includes Bernoulli percolation and the Ising model \ref b.Hag:ustlim/.
The UST is very interesting mathematically, partly because
it is closely related to the theory of resistor networks,
potential theory, random walks, LERW, and in dimension $2$,
also domino tilings.
The paper \ref b.BLPS:usf/ gives a comprehensive study of 
uniform spanning trees (and forests), following earlier pioneering work
\ref b.Aldous:ust/, \ref b.Broder:ust/, \ref b.Pemantle:dichot/, 
\ref b.BP:ust/, \ref b.Hag:ustlim/.
A survey of current UST theory can be found in \ref b.Lyons:usfsurvey/.

Let $G$ be a connected graph.  A {\bf forest}
is a subgraph of $G$ that has no cycles.
A {\bf tree} is a connected forest.  A subgraph of $G$
is {\bf spanning} if it contains $\verts(G)$, the set of vertices of $G$.
We will be concerned with spanning trees.
Since a spanning tree is determined by its edges, we often don't make
a distinction between the spanning tree $\T$ and its set of edges $\edges(T)$. 

If $G$ is finite,
a uniform spanning tree (UST) in $G$ is a random spanning
tree $\T\subset G$, selected according to the uniform measure.
(That is, $\P[\T=T_1]=\P[\T=T_2]$, whenever $T_1$ and $T_2$ are
spanning trees of $G$.)

It turns out that UST's are very closely related to LERW's.
If $a,b\in\verts(G)$, then the (unique) path in the UST joining
$a$ and $b$ has the same law as the LERW from $a$ to $b$ in
$G$.  (This, in particular, implies that the LERW from $a$ to
$b$ has the same law as the LERW from $b$ to $a$.)
Wilson's algorithm \ref b.Wilson:alg/, which will be described in \ref s.back/,
is a very useful method to build the UST by running LERW's.

R.~Lyons proposed (see \ref b.Pemantle:dichot/)
to extend the notion of UST to infinite graphs.
Let $G$ be an infinite connected graph.  Consider a nested
sequence of connected finite sugraphs $G_1\subset G_2\subset\cdots\subset G$
such that $G=\bigcup_j G_j$.  For each $j$, the uniform spanning
tree measure $\mea_j$ on $G_j$ may be considered as a measure
on $2^{\edges(G)}$, the $\sigma$-field of subsets of the edges of $G$,
generated by the sets of the form $\bigl\{F\subset \edges(G)\st e\in F\bigr\}$,
$e\in\edges(G)$.  
Using monotonicity properties, it can be shown that the weak limit
$\mea:=\lim_{j\to\infty}\mea_j$ exists, and does not depend on
the sequence $\{G_j\}$.  It is called
the free uniform spanning forest measure (FSF) on $G$.
(The reason for the word `free', is that there's another natural
kind of limit, the wired uniform spanning forest (WSF).
On $\Z^d$, these two measures agree.)
R.~Pemantle \ref b.Pemantle:dichot/
proved that if $d\leq 4$, the FSF measure on $\Z^d$ is supported
on spanning trees (that is, if $\T$ is random
and its law is the FSF measure, then $\T$ is a.s.\ a spanning
tree), while if $d\geq 5$, the measure is supported on
disconnected spanning forests.

Let us now restrict attention to the case $G=\Z^2$.
Since it is supported on
spanning trees, we call the FSF measure on $\Z^2$ 
the {\bf uniform spanning tree} (UST) on $\Z^2$.
I.~Benjamini \ref b.Benjamini:ustsl/ and R.~Kenyon~\ref b.Kenyon:longrange/
studied asymptotic properties of the UST on a rescaled
grid $\delta\Z^2$, with $\delta$ small, but did not attempt to
define the scaling limit.
Aizenman, Burchard, Newman and Wilson \ref b.ABNW/ defined the
scaling limit of UST (and other tree processes) in $\Z^2$,
and studied some of their properties.

We present a different definition for the scaling limit
of the UST.  Let $\delta>0$.  Again, we think of
$\dZ$ as a subset of the sphere $\Sp2=\R^2\cup\{\infty\}$.
Let $\hT_\delta$ be the UST on $\dZ$, union with the point
at infinity.  Then $\hT_\delta$ can be thought of as a 
random compact subset of $\Sp2$.  However, it is fruitless
to consider the weak limit as $\delta\to0$ of the law of $\hT_\delta$ as
a measure on the Hausdorff space $\Ha(\Sp2)$, since
the limit measure is an atomic measure supported on the
single point in $\Ha(\Sp2)$, which is all of $\Sp2$.

Given two points $a,b\in\hT_\delta$, let
$\path_{a,b}$ be the unique path in $\hT_\delta$ with endpoints
$a$ and $b$; allowing for the possibility $\path_{a,b}=\{a\}$, when
$a=b$.  (It was proved by R.~Pemantle that a.s.\ the UST $\T$ in $\Z^2$ has
a single end; that is, there is a unique infinite ray in $\T$ starting at $0$.
This implies that indeed $\path_{a,b}$ exists and is unique, not only for
$\T$, but also for $\T\cup\{\infty\}$.)
Let $\ST_\delta=\ST_\delta(\hT_\delta)$ be the collection of all triplets,
$(a,b,\path_{a,b})$, where $a,b\in\hT_\delta$.
Then $\ST_\delta$ is a closed subset of $\Sp2\times\Sp2\times\Ha(\Sp2)$,
and the law $\dlaw$ of $\ST_\delta$ is a probability measure
on the compact space $\Ha\bigl(\Sp2\times\Sp2\times\Ha(\Sp2)\bigr)$.
By compactness, there is a subsequential weak limit
$\law$ of $\dlaw$ as $\delta\to 0$, which is a probability measure
on $\Ha\bigl(\Sp2\times\Sp2\times\Ha(\Sp2)\bigr)$.
We call $\law$ a {\bf subsequential UST scaling limit in $\Z^2$}.  
If $\law:=\lim_{\delta\to 0}\dlaw$, as a weak limit, then  $\law$ is the
{\bf UST scaling limit}.

We prove

\procl t.ustlim
Let $\law$ be a subsequential UST scaling limit in $\Z^2$,
and let $\ST\in \Ha\bigl(\Sp2\times\Sp2\times\Ha(\Sp2)\bigr)$ be
a random variable with law $\law$.  Then the following holds
a.s.
\beginitems
\itemrm{(i)} For every $(a,b)\in \Sp2\times\Sp2$, there
is some $\path\in\Ha(\Sp2)$ such that $(a,b,\path)\in\ST$.
For almost every $(a,b)\in\Sp2\times\Sp2$, this $\path$ is unique.
\itemrm{(ii)} For every $(a,b,\path)\in\ST$, if $a\neq b$, then
$\path$ is a simple path; that is, homeomorphic
to $[0,1]$.  If $a=b$, then $\path$ is a single point or
homeomorphic to a circle.  For almost every $a\in\Sp2$,
the only $\path$ such that $(a,a,\path)\in\ST$
is $\{a\}$.
\itemrm{(iii)} The {\bf trunk},
$$
\Tr:=\bigcup_{(a,b,\path)\in\ST}\bigl(\path-\{a,b\}\bigr)
\,,
$$
is a topological tree (in the sense of \ref d.toptree/), which
is dense in $\Sp2$.
\itemrm{(iv)} For each $x\in\Tr$, there are at most
three connected components of $\Tr-\{x\}$.
\enditems
\endprocl

This theorem basically answers all the topological questions about
the UST scaling limit on $\Z^2$.  It is sharp, in the sense
that all the ``almost every'' clauses cannot be replaced by
``every''.  Benjamini~\ref b.Benjamini:ustsl/ proved a result which
is closely related to item (iv) of the theorem, and
\ref b.ABNW/ proved (in a different language)
that (iv) holds with ``three'' replaced by some unspecified
constant.

The {\bf dual} of a spanning tree $T\subset \Z^2$ is the
spanning subgraph of the
dual graph $(1/2,1/2)+\Z^2$ containing all edges that do not
intersect edges in $T$.  It turns out that duality is  
measure preserving from the UST on $\Z^2$ to the UST on the dual grid.
The key to the proof of \ref t.ustlim/ is the statement that the
trunk is disjoint from the trunk of the dual UST scaling limit.

Let us stress that \ref t.ustlim/ is not contingent on \ref g.confinv/.
The only contingent theorems proved in this paper are
\ref t.lesl/, and \ref t.confust/, which says that
\ref g.confinv/ implies conformal invariance 
for the scaling limit of the UST on subdomains of $\C$.

Recent work of R.~Kenyon \ref b.Kenyon:conf/, \ref b.Kenyon:5o4/
proves some conformal invariance results
for domino tilings of domains in the plane.  
There is an explicit correspondence between the UST in $\Z^2$ and
domino tilings of a finer grid.
Based on this correspondence, some
properties of the UST can be proved using Kenyon's machinery.
For example, Kenyon has shown \ref b.Kenyon:5o4/ that the expected number
of edges
in a LERW joining two boundary vertices in $(\dZ)\cap [0,1]^2$ 
(whose distance from each other is bounded from below)
grows like $\delta^{-5/4}$, as $\delta\to 0$.
He can also show \ref b.Kenyon:longrange/ that the weak limit as $\delta\to0$ 
of the distribution
of the UST meeting point of three boundary vertices of $D\cap(\dZ)$
is equivariant with respect to conformal maps.
That can be viewed as a partial conformal invariance result
for the UST scaling limit.
Another example for the applications of Kenyon's work to
the UST appears in \ref s.windl/.
It seems plausible that perhaps soon there would be a proof
of \ref g.confinv/.

\subsection{SLE with other parameters, critical percolation, and the
UST Peano curve}
Let $\lec\geq 0$, and take $\zeta(t):= \B(-\lec t)$, where $\B$
is as above, Brownian motion on $\bd\U$, started from a uniform
random point.  Then there is a solution $f(z,t)$ of 
\ref e.loew/ and \ref e.bdval/, and for each $t\leq 0$
$f_t = f(\cdot,z)$ is a conformal map from $\U$ into some subdomain
$D_t\subset\U$.  We call the process
$\sle^\lec_t:=\U-D_t$, $t\leq 0$, the
{\bf stochastic \lowner\ evolution} (SLE)
with parameter $\lec$.
It is not always the case that $\sle^\lec_t$ is a simple path.
Let $\KK$ be the set of all $\lec\geq 0$ such that for all $t<0$ the
set $\sle^\lec_t$ is a.s.\ a simple path.
We show in \ref s.crit/ that $\sup\KK\leq\crit$,
and conjecture that $\KK=[0,\crit]$.

In the past, there has been some work on the question of
which \lowner{} parameters $\zeta$ produce slitted disk mappings
(\ref b.Kufarev/, \ref b.Pommerenke:Lowner/),
but only limited progress has been made.
Partly motivated by the present work,
Marshall and Rohde \ref b.MR:lowner/ have looked into this problem
again, and have shown that when $\zeta$ satisfies a 
H\"older condition with exponent $1/2$ and H\"older(1/2) norm
less than some constant, the maps $f_t$ are onto slitted disks,
and this may fail when $\zeta$ has finite but large
H\"older(1/2) norm.

Given some $\lec>0$, even if $\lec\notin\KK$, the process
$\sle^\lec_t$, $t\leq 0$, is quite interesting.
It is a celebrated conjecture that critical Bernoulli percolation
on lattices in $\R^2$ exhibits conformal invariance in the
scaling limit \ref b.Langlands:Bull/.
Assuming such a conjecture, we plan to prove in a subsequent work that
a process similar to SLE describes the scaling limit of the outer boundary
of the union of all critical percolation clusters in a domain 
$D$ which intersect a fixed arc on the boundary of $D$.
We also plan to prove that this implies Cardy's \ref b.Cardy/
conjectured formula for the limiting crossing probabilities
of critical percolation, and higher order generalizations
of this formula.

Let us now briefly explain this.
In \ref f.pcurve/, each of the hexagons is colored black
with probability $1/2$, independently, except that the hexagons
intersecting the positive real ray are all white, and the hexagons
intersecting the negative real ray are all black.  
Then there is a boundary path $\beta$, passing through $0$ 
and separating the black and the white regions adjacent to $0$.
Note that the percolation in the figure is equivalent
to Bernoulli$(1/2)$ percolation on the triangular grid,
which is critical.
(See \ref b.Grimmett:book/ for background on percolation.)
The intersection of $\beta$ with the upper half plane,
$\H:=\{z\in\C\st \im z>0\}$, which is indicated in the
picture, is a random path in $\H$ connecting the
boundary points $0$ and $\infty$.

\midinsert
\centerline{\epsfysize=2.4in\epsfbox{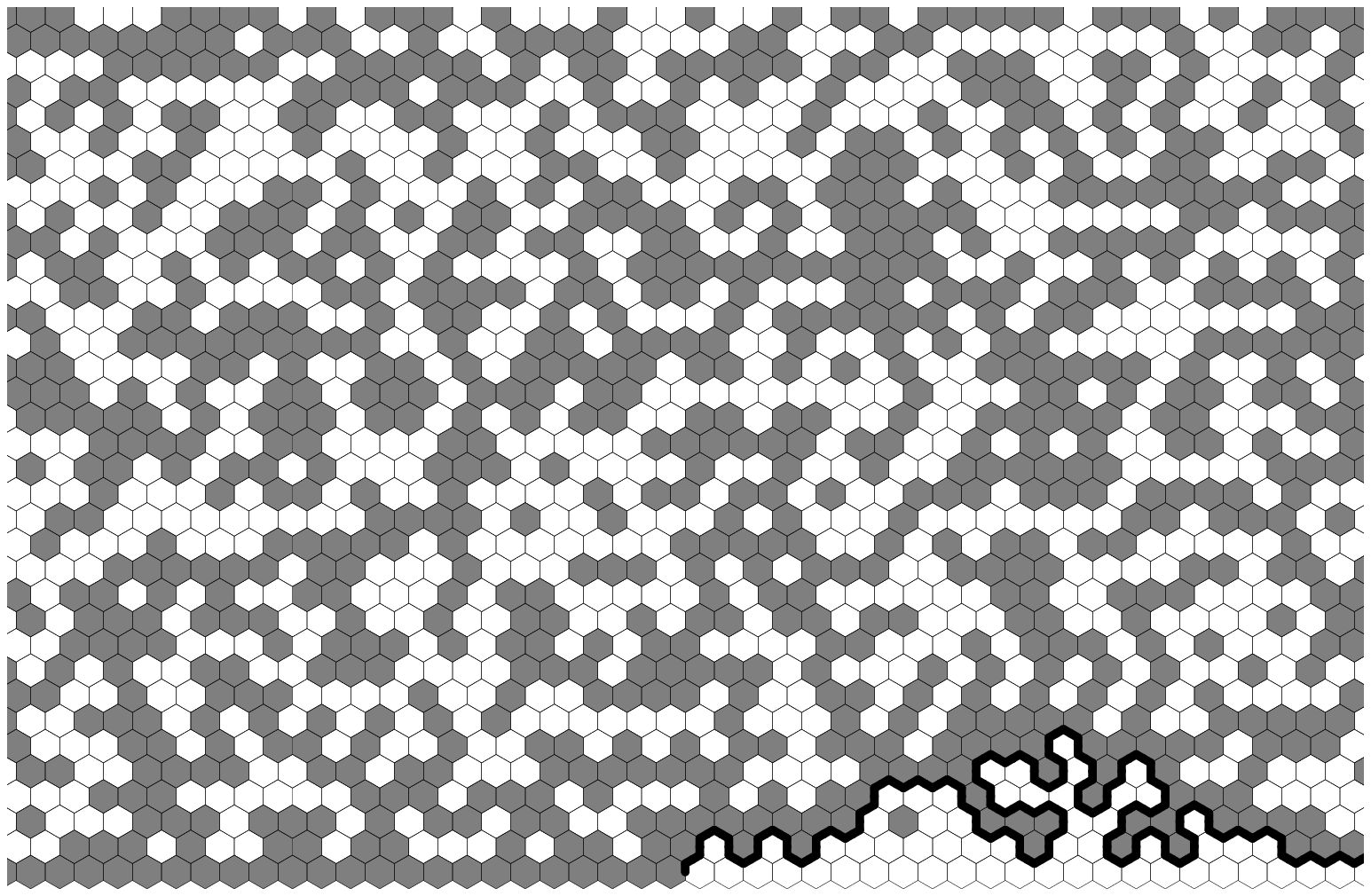}}
\caption{\figlabel{pcurve}\enspace
The boundary curve for critical percolation with mixed boundary
conditions.}
\endinsert

A subsequential scaling limit of $\beta\cap\H$ exists,
by compactness, and naturally, we believe that the weak limit exists.
Let $\gamma$ be the scaling limit curve.
The physics wisdom (unproven, perhaps not even precisely formulated,
but well supported) is that the scaling limit of the
``external boundary'' of macroscopic critical percolation clusters
in two dimensions has dimension $7/4$ and is {\bf not} a simple path
\ref b.ADA/, and we believe that this is true for $\gamma$.

In a subsequent paper, we plan to prove (by adapting the proof
of \ref t.lesl/), under the assumption
of a conformal invariance conjecture for the scaling limit
of critical percolation, that $\gamma$ can be described using a 
\lowner-like differential equation in the upper half
plane with Brownian motion parameter, as follows.
Consider the differential equation
$$
{\bd\over\bd t}f_t(z) = 
{-2f_t'(z)\over\zeta(t)-z}
\,,\qquad\forall z\in\H\,,\forall t\in(-\infty,0]
\,,
\label e.loewh
$$ 
where $\zeta(t)=\B(-\lec t)$, $\B$ is Brownian
motion on $\R$ starting at $\B(0)=0$, and $f_0(z)=z$.
Then $f_t$ is a conformal mapping from $\H$ onto a subdomain
of $\H$, which is normalized by the so-called {\bf hydrodynamic
normalization} 
$$
\lim_{z\to\infty}f_t(z)-z=0\,.
\label e.hydro
$$
The claim is that for $\lec=6$, the image of the path
$t\mapsto f_t\bigl(\zeta(t)\bigr)$ has the same distribution as $\gamma$.
{}From this, one can derive Cardy's \ref b.Cardy/ conjectured formula for the
limiting crossing probabilities of critical percolation, as well as
some higher order generalizations.  This will be done in subsequent
work, but basically depends on the ideas appearing in \ref s.crit/ below.

A similar representation applies to the scaling limit of the Peano curve which
winds around the UST (this curve was discussed in \ref b.Duplantier:Peano/ and
mentioned in \ref b.BLPS:usf/), but with $\lec=8$. 
Given a domain $D\subset\Sp2$, whose boundary is a simple closed path,
and given two distinct points $a,b\in\bd D$, there is
a naturally defined (subsequential) scaling limit of the Peano
curve of the UST in $D$, with appropriate boundary conditions,
and the scaling limit is an (unparameterized) curve from $a$ to $b$,
whose image covers $\br D$.
One can show that \ref g.confinv/ implies a conformal invariance
property for the scaling limit.  Based on this,
it should be possible to adapt the proof of \ref t.lesl/ to 
show that \ref g.confinv/ implies a representation of the
form \ref e.loewh/ for this Peano scaling limit when $D=\H$, $a=0$,
$b=\infty$, and $\lec=8$.
We give a brief overview of this in \ref s.peano/, and hope to
give a more thorough treatment in a subsequent paper.

The differential equation \ref e.loewh/ is very similar to 
\lowner{}'s equation, and the only essential difference is that
a normalization at an interior point for the maps $f_t$
is replaced by the hydrodynamic normalization \ref e.hydro/ at a 
boundary point ($\infty$).  The interior point normalization 
is natural for the LERW scaling limit, because the LERW is
a path from an interior point to the boundary of the domain.
The Peano curve and the boundary of percolation clusters, as discussed
above, are paths joining two boundary points, and hence the
hydrodynamic normalization is more appropriate for them.

Although \lowner{}'s Slit Mapping Theorem mentioned above applies
to domains of the form $\U-\alpha$, where $\alpha$ is a simple path
in $\U-\{0\}$ with one endpoint in $\bd \U$, 
Pommerenke \ref b.Pommerenke:Lowner/ has a generalization, which
is valid for some paths $\alpha$ which are not simple paths.
It is this genaralization (or rather, its version in $\H$ with
the hydrodynamic normalization) which will substitute
\lowner's Slit Mapping Theorem for the treatment of the percolation
boundary or Peano curve scaling limits.

The emerging picture is that different values of $\lec$ in
the differential equations \ref e.loew/ or \ref e.loewh/
produce paths which are scaling limits of naturally defined
processes, and that these paths can be space-filling,
or simple paths, or neither, depending on the parameter $\lec$.

\subsection{Acknowledgement}
I wish to express gratitude to Itai Benjamini,
Rick Kenyon 
and David Wilson for inspiring discussions and helpful information.
\ref l.tightdiam/ has been obtained jointly with Itai Benjamini.
Mladen Bestvina, Brian Bodwitch,
Steve Evans, Yakar Kannai, Greg Lawler, Russ Lyons, Yuval Peres,
Steffen Rohde, Jeff Steif, Benjamin Weiss and Wendelin Werner
have provided very helpful advice.

\bsection {Some background and terminology}{s.back}
\resetConsts

This section will introduce some notations which will be used,
and discuss some of the necessary background. 
We begin with a review of uniform spanning trees and forests.
The reader may consult \ref b.BLPS:usf/ for a comprehensive
treatment of that subject.

\subsection{The domination principle}
Suppose that $H$ and $H'$ are two random subsets
of some set.  We say that $H'$ {\bf stochastically
dominates} $H$ if there is a probability measure $\mu$ on
pairs $(A,B)$ such that $A$ has the same law as $H$,
$B$ has the same law as $H'$, and
$\mu\bigl\{(A,B)\st B\supset A\bigr\}=1$.
Such a $\mu$ is called a {\bf monotone coupling} of $H$ and $H'$.

Let $G$ be a finite connected graph, and let $G_0$ be
a connected nonempty subgraph.  Let $M$ be the set of
vertices of $G_0$ that are incident 
with some edge in $\edges(G)-\edges(G_0)$.
(We let $\edges(G)$ and $\verts(G)$ denote the edges and
vertices of $G$, respectively.)
Let $G_0^W$ be the
graph obtained from $G$ by identifying all the vertices
in $M$ to a single vertex, called
the {\bf wired vertex}.
Then $G_0^W$ is called the {\bf wired graph}
associated to the pair $(G_0,G)$.
Let $\T$ be the UST on $G$, let $\T_0^F$ be the
UST on $G_0$, and let $\T_0^W$ be the UST on $G_0^W$.
Then $\T_0^F$ is called the {\bf free spanning tree}
of the pair $(G_0,G)$, and $\T_0^W$ is the
{\bf wired spanning tree} of the pair $(G_0,G)$.
Sometimes, we call $\T_0^F$ [respectively, $T_0^W$] the
UST on $G_0$ with free [respectively, wired] boundary conditions.
The {\bf domination principle} states that
$\T_0^F\cap G_0$ stochastically dominates $\T\cap G_0$,
and that $\T\cap G_0$, stochastically dominates $\T_0^W\cap G_0$.

Now let $G$ be an infinite connected graph, and let
$G_1\subset G_2\subset \cdots$ be an infinite
sequence of finite connected subgraphs satisfying $\bigcup_j G_j=G$.
Let $\mu_j^W$ be the law of the wired spanning tree of $(G_j,G)$.
Based on the domination principle, it is easy to verify that
the weak limit $\mu^W$ of $\mu_j^W$ exists, and is a probability measure
on spanning forests of $G$.  It is called
the {\bf wired spanning forest of $G$} (WSF).

The domination principle, when appropriatly interpreted,
carries over to infinite and to disconnected graphs as well.
If $G$ is a disconnected graph, we take the FSF [respectively, WSF] on
$G$ to be the spanning forest of $G$ whose intersection with every
component of $G$ is the FSF [respectively, WSF] of that component,
and with the restriction to the different components being independent.
The more general formulation of the domination principle states
that when $G_0$ is a subgraph of $G$, then 
$\T_0^F\cap G_0$ stochastically dominates $\T\cap G_0$ 
and $\T\cap G_0$ stochastically dominates $\T_0^W\cap G_0$,
where $\T,\T_0^F$ and $\T_0^W$ are the FSF on $G,G_0$ and $G_0^W$,
respectively, and the same statement holds when FSF is replaced by WSF.

On all recurrent connected graphs, the WSF is equal to the FSF, 
and both are trees.  Therefore, on recurrent graphs we shall refer to
this measure as the uniform spanning tree (UST).

\subsection{Grid approximations of domains}
Let $D\subset\C$ be a domain, that is, an open, connected set.
Given $\delta>0$, we define a graph $G=\G(D,\delta)$, which
is a discrete approximation of the domain $D$ in the grid $\dZ$,
as follows.  The {\bf interior vertices} $\Vi(G)=\Vi(D,\delta)$ of $G$
are the vertices of $\dZ$ which are in $D$.  The
{\bf boundary vertices} $\Vb(G)=\Vb(D,\delta)$ are the points of intersection
of the edges of the grid $\dZ$ with $\bd D$, the boundary of $D$.
The vertices of $G$ are $\verts(G)=\Vb(G)\cup\Vi(G)$.
If $a,b\in\verts(G)$ are distinct, then $[a,b]$ is an edge of $G$
iff there is an edge $e\in\edges(\dZ)$
such that the open segment $\bigl\{ta+(1-t)b\st t\in(0,1)\bigr\}$
is contained in $D\cap e$.

We will often be considering random walks on $\G(D,\delta)$ starting
at $0$, when $0\in D$.  It will be useful to denote by $\Vba$ the set
of vertices $v\in\Vb(D,\delta)$ such that there is a path
from $0$ in $\G(D,\delta)$ whose intersection with $\Vb$ is $v$.  

The wired graph, $\G^W(D,\delta)$, associated with $D$ is $\G(D,\delta)$ with
all the vertices $\Vb(D,\delta)$ collapsed to a single vertex,
which we simply denote $\bd D$.

We shall often not distinguish between a graph and its
planar embedding, if it has an obvious planar embedding. 
For example, the UST on $G$
will also be interpreted as a random set in the plane.

\subsection{Wilson's algorithm}
Let $G$ be a finite graph.  Wilson's
algorithm \ref b.Wilson:alg/ for generating a UST in $G$ proceeds
as follows.
Let $v_0\in\verts(G)$ be an arbitrary vertex (which we call the root),
and set $T_0:=\{v_0\}$.
Inductively, assume that a tree $T_j\subset G$ has been constructed.
If $\verts(T_j)\neq\verts(G)$,
choose a vertex $v_{j+1}\in\verts(G)-\verts(T_j)$, 
let $W_{j+1}$ be LERW from $v$ to $T_j$ in $G$, and set
$T_{j+1}:=T_j\cup W_{j+1}$.
Otherwise $\verts(T_j)=\verts(G)$, and the algorithm stops and outputs
$T_j$.  
It is somewhat surprising, but true, that no matter how the
choices of the vertices $v_j$ are made, the output of the algorithm
is a tree chosen according to the uniform measure.

If $G$ is infinite, connected and recurrent, Wilson's algorithm also ``works''.
When the subtree $T_j$ generated by the algorithm includes all the
vertices in a certain finite set $K\subset\verts(G)$,
the subtree of $T_j$ spanned by $K$ (that is, the minimal connected
subgraph of $T_j$ that contains $K$) has the same law as the
subtree of the UST of $G$ spanned by $K$.
(There is also a version of Wilson's algorithm which is useful for
generating the WSF of a transient graph, but we shall not need this.)

\subsection{Harmonic measure estimates}
Because of Wilson's algorithm, many questions about the UST can
be reduced to questions about simple random walks (SRW's).
It is therefore hardly surprising that we often need to
obtain a {\bf harmonic measure estimate}; that is,
an estimate on the probability that SRW starting from a vertex $v$
will hit a certain set of vertices $K_1$ before hitting
another set $K_0$.  As a function of $v$, this probability is
harmonic\Ftnote{A function $h$ 
is {\bf harmonic} at $v$, if $h(v)$ is the
average of the value of $h$ on the neighbors of $v$.}
away from $K_0\cup K_1$.

Almost all the harmonic measure estimates which we will use
are entirely elementary, and follow from the following easy fact.
Consider  an annulus $A=A(p,r,R)$, with center $p$,
inner radius $r$ and outer radius $R>r$.  Suppose that $\delta$
is sufficiently small so that there is a path
in $\dZ\cap A$ which separates the boundary components of $A$.
Let $q\in\dZ\cap A$ be some vertex such that the distance
from $q$ to the boundary $\bd A$ is at least $r/c$, where $c>0$ is some
constant.  Let $X$ be the image of SRW starting from $q$,
which is stopped when it first leaves $A$.  Then the
probability that $X$ contains a path separating the boundary
components of $A$ is bounded below by some positive function of $c$.
(This can be proved directly using only the Markov property 
and the invariance of SRW under the automorphisms of $\dZ$.)
One consequence of this fact and the Markov property, which
we will often use, is as follows.

\procl l.harmm
Suppose that $K$ is a connected subgraph
of $\dZ$ of diameter at least $R$, and $v\in\dZ$. 
Then the probability that SRW starting from $v$ will exit the
ball $B(v,R)$ before hitting $K$ is at most
$\nco a\bigl(\dist(v,K)/R\bigr)^{\nco b}$, where
$\co a,\co b>0$ are absolute constants.
\endprocl

This lemma holds for the spherical as well as Euclidean metric.

\subsection{Laplacian random walk}
Although this will not be needed in the paper, we have to
mention another interpretation of LERW.
Let $G$ be a finite connected graph, let $K\subset\verts(G)$ be a 
set of vertices, and let $a\in\verts(G)-K$.
The LERW from $a$ to $K$ can also
be inductively constructed, as follows.
Suppose that the first $n$ vertices $a=\LE(1),\LE(2),\dots,\LE(n)$ have been
determined and $\LE(n)\notin K$. 
Let $h_n:\verts(G)\to[0,1]$ be the function 
which is $1$ on $K$, $0$ on $\{\LE(1),\dots,\LE(n)\}$
and harmonic
on $\verts(G)-\bigl(K\cup \{\LE(1),\dots,\LE(n)\}\bigr)$.
Then $\LE(n+1)$ is chosen among the neighbors
$w$ of $\LE(n)$, with probability proportional to $h_n(w)$.

This formulation of the LERW may serve as a heuristic for
\ref g.confinv/, since discrete harmonic functions are good
approximation for continuous harmonic functions, and continuous
harmonic functions in 2D have conformal invariance properties.
However, this heuristic is quite weak, since near a non-smooth boundary
of a domain, the approximation is not good.

\subsection{Weak convergence of measures}
We now recall several facts and definitions regarding weak convergence.
The reader may consult \crfl{EK:markov}{Chap.~3} for proofs and
further references.
Let $(X,d)$ be a compact metric space, and let $\{\mu_j\}$ be
a sequence of Borel probability measures on $X$.
$\mu$ be a Borel probability measure on $X$.
Saying that the sequence $\mu_j$ converges
weakly to $\mu$ means that
$\lim_j \int f\, d\mu_j=\int f \,d\mu$ for all continuous $f:X\to\R$.
The Prohorov metric on the space of Borel probability
measures on $X$ is defined by
$$
\dlr(\mu,\mu'):= \inf\bigl\{\eps>0\st
\mu(K)\leq \mu'(N_\eps (K)) + \eps
\hbox{ for all closed }K\subset X\bigr\}
\,,
$$
where $N_\eps(K):=\bigcup_{x\in K} B(x,\eps)$
is the $\eps$-neighborhood of $K$.
The space of Borel probability measures on $X$ is compact
with respect to the Prohorov metric, and weak convergence
is equivalent to convergence in the Prohorov metric.

Let $\Cal M(\mu,\mu')$ be the collection of all Borel
measures $\nu$ on $X\times X$ such that $\nu(A\times X)=\mu(A)$
and $\nu(X\times A)=\mu'(A)$ for all measurable $A\subset X$.
Such a $\nu$ is called a {\bf coupling} of $\mu$ and $\mu'$.
The Prohorov metric satisfies
$$
\dlr(\mu,\mu') =
\inf_{\nu\in\Cal M(\mu,\mu')}
\inf \Bigl\{\eps>0\st \nu\{(x,y)\st d(x,y)\geq\eps\}\leq\eps\Bigr\}
\,.
\label e.pcoup
$$
In other words, $\dlr(\mu,\mu')<\eps$
means that one can find a probability space $(\Omega,\P)$ and two $X$ valued
random variables $x,y:\Omega\to X$,
such that $\P[d(x,y)\geq\eps]\leq\eps$ and such that $x$ has
law $\mu$ and $y$ has law $\mu'$.  This is obtained by
taking an appropriate $\P=\nu\in\Cal M(\mu,\mu')$,
$\Omega:=X\times X$, and letting $x$ and $y$ be
the projections on the first and second  factors, respectively.

\subsection{Conformal maps}
We review some elementary facts about conformal (aka univalent) mappings,
as may be found in \ref b.Duren:book/, for example.
Let $D\subset\C$ be some domain.
A continuous map $f:D\to \C$,
which is injective and complex-differentiable is {\bf conformal}. 
If $f$ is conformal, then $f^{-1}:f(D)\to\C$ is also conformal.

Let $D\subsetneqq\C$ be simply connected.
Then Riemann's Mapping Theorem states that 
there is a conformal homeomorphism $f=f_D$ from $\U$ onto $D$. 
Suppose also that $0\in D$, then $f$ can be chosen
to satisfy the normalizations $f(0)=0$ and $f'(0)>0$,
which render $f$ unique.
In this case, the number $f'(0)$ is called the
{\bf conformal radius of $D$} (with respect to $0$).
The Schwarz Lemma implies that 
$$
f'(0)\geq \inf\bigl\{|z|\st z\notin f(\U)\bigr\}
\,,
\label e.schw
$$
while on the other hand, the Koebe $1/4$ Theorem gives
$$
f'(0)\leq 4 \inf\bigl\{|z|\st z\notin f(\U)\bigr\}
\,.
\label e.k1q
$$
Hence, up to a factor of $4$, the conformal radius can be determined
from the in-radius.

If $0<a<b<\infty$, then the set of all conformal maps $f:\U\to\C$
satisfying $f(0)=0$ and
$a\leq |f'(0)|\leq b$ is compact, in the topology of uniform
convergence on compact subsets of $\U$.
If $f_j:\U\to\C$ are conformal, $f_j(0)=0$, and $f_j\to f$ locally uniformly,
then the image $f(\U)$ can be described in terms of the images
$D_j:=f_j(\U)$.  Let $D$ be the maximal open connected
set containing $0$ and contained in
$\bigcup_{n=1}^\infty\bigcap_{j=n}^\infty D_j$.
If $D\neq\emptyset$, then $f(\U)=D$;
otherwise, $f(\U)=\{0\}$.  This is called Carath\'eodory's Kernel
Convergence Theorem.

If $f:\U\to D$ is a conformal homeomorphism onto $D$, and
$\bd D$ is a simple closed curve, then $f$ extends continuously
to $\bd \U$.  The same is true if $D=\U-\beta$, where $\beta$
is a simple path.

\bsection {No Loops}{s.noloops}
\resetConsts

In this section, we prove that any LERW subsequential
scaling limit is supported on
the set of simple paths.  That is, we prove \ref t.noloop/.

Let $\zerp$ and $\onep$ be distinct points $\R^2$.
For each $\delta>0$, let $\zer$ and $\one$ be
vertices of $\dZ$ closest to $\zerp$ and $\onep$, respectively.
Note that in $\dZ$ the combinatorial distance between
two vertices $v,v'\in\dZ$ is
$\delta^{-1}\|v-v'\|_1$.
However, all metric notions we use will refer to the Euclidean 
or spherical distance.  In this section, we will mainly use the
Euclidean metric.

Let $\RW$ be a random walk on $\dZ$ starting at $\one$
and stopped when $\zer$ is reached for the first time.
Let $\Pa=\Pa_\delta$ denote the loop-erasure of $\RW$, and let $\mun$
denote the law of $\Pa_\delta$.
The following lemma will show that the diameter of
$\Pa_\delta$ is ``tight''.

\procl l.tightdiam 
$$
\mun\bigl[{\diam \Pa> s\, \dist(\zer,\one)}\bigr]
\leq \nco c s^{-\nco d}
\,,
$$
where $\co c,\co d>0$ are absolute constants.
\endprocl

The proof is based on Wilson's algorithm and an elementary 
harmonic measure estimate.
A more precise estimate can be obtained by using the discrete Beurling
Projection Theorem (see \ref b.Kesten:hitting/ or \ref b.Lawler:makarov/). 

\proof
Set $r:=\dist(\zer,\one)$, $R:=s\, \dist(\zer,\one)/4$,
and let $z\in \dZ$ be some vertex such that $\dist(\zer,z)\geq10R$.
For $a,b\in\dZ$,
let $\Paa(a,b)$ denote the path in the UST of $\dZ$
joining $a$ and $b$.
Let $w$ be the meeting point of $\zer,\one,z$ in the UST;
that is, the vertex in
$\Paa(\zer,\one)\cap\Paa(\zer,z)\cap\Paa(\one,z)$.
By Wilson's algorithm,
the distribution of $\Pa$ is identical with the
distribution of $\Paa(\zer,w)\cup\Paa(\one,w)$.

We now estimate $\PBig{\diam\bigl(\Paa(a,w)\bigr)> s\, \dist(\zer,\one)}$.
Using Wilson's algorithm, we may generate $\Paa(\zer,w)$ by letting
$\Paa(\one,z)$ be LERW from $\one$ to $z$,
and letting $\Paa(\zer,w)$ be LERW from $\zer$ to $\Paa(\one,z)$.
Condition on $\Paa(\one,z)$.
By \ref l.harmm/,
the probability that SRW starting at $\zer$ will
exist $B(\zer,R)$ before hitting $\Paa(\one,z)$ is
at most $\nco a (r/R)^{\nco b}$, for some constants
$\co a,\co b>0$.  Therefore,
$\P[\diam \Paa(\zer,w)>2 R]\leq \co a (r/R)^{\co b}$.
Moreover, the same estimate holds for $\diam \Paa(\one,w)$.
Since $\Paa(\zer,\one)=\Paa(\zer,w)\cup\Paa(\one,w)$,
we get
$\P[\diam \Paa(\zer,\one)>4 R]\leq 2\co a (r/R)^{\co b}$.
This completes the proof of the lemma.
\Qed

\procl r.closetopt
The proof of the lemma can be easily adapted to show that
if $a,b\in\dZ$ and $K\subset\dZ$, then the probability that
LERW from $a$ to $K$ will intersect $B(b,r)$ 
is at most $\nco c  \bigl(r/\dist(a,b)\bigr)^{\nco d}$,
provided $r\geq \delta$, where
$\co c,\co d>0$ are absolute constants.
\endprocl

\procl d.qlp
Let $z_0\in\R^2$, $r,\eps>0$.  An $(z_0,r,\eps)$-quasi-loop
in a path $\Pa$ is a pair $a,b\in\Pa$ with
$a,b\in B(z_0,r)$, $\dist(a,b)\leq\eps$,
such that the subarc of $\Pa$ with endpoints $a,b$ is
not contained in $B(z_0,2r)$.
Let $\ev A(z_0,r,\eps)$ denote the set of simple
paths in $\R^2$ that have a $(z_0,r,\eps)$-quasi-loop.
\endprocl

\procl l.noloop
Let $c$ be the distance from $\zer$ to $\one$,
let $r\in(0,c/4)$, $\eps>0$ and $z_0\in\R^2$.
Then $\lim_{\eps\to 0} \mun[{\ev A(z_0,r,\eps)}] = 0$, uniformly
in $\delta$.
\endprocl

\proof
Let $B_1=B(z_0,r)$ and $B_2=B(z_0,2r)$.
The distance from $B_2$ to at least one of the points $\zer,\one$
is at least $c/2$.
By symmetry, we assume with no loss of generality that
$\dist\bigl(\one,B_2\bigr)\geq c/2$.
Let $\eps_1\in(0,c/4)$, and let $q$ be a vertex in $\dZ$ such that
$\dist\bigl(q,\one\bigr)\in[\eps_1-\delta,\eps_1]$.

Let $\Paq$ be a LERW from $q$ to $\one$ in $\dZ$.
Let $\RW$ be an independent simple random walk from $\zer$.
Let $\RW'$ be the part of the walk $\RW$ until $\Paq$ is first hit.
Then, by Wilson's algorithm, $\Pa_\delta$ has the same distribution
as the arc connecting $\zer$ to $\one$ in $\LE(\RW')\cup\Paq$.

Let $\ev A'$ be the event that $\LE(\RW')$ has a
$(z_0,r,\eps)$-quasi-loop, and let $\ev C$ be the
event that $\Paq$ intersects $B_2$. 
Then
$$
\mun\bigl[{\ev A(z_0,r,\eps)}\bigr] \leq \P[\ev A'] + \P[\ev C]
\,.
\label e.cas
$$
Since $\dist\bigl(q,\one\bigr)\leq \eps_1$ and
$\dist\bigl(B_2,\one\bigr)\geq c/2$,
from \ref l.tightdiam/ we get an estimate of the form 
$$
\P[\ev C] \leq \nco {aa} (\eps_1/c)^{\nco {ab}}
\,.
\label e.pc
$$ 

We now find an upper bound for $\P[\ev A']$.
Let $s_1$ be the first time $s\geq 0$ such that $\RW(s)\in B_1$.
Let $t_1$ be the first time $t\geq s_1$ such that $\RW(t_1)\notin B_2$.
Inductively, define $s_j$ to be the first time $s\geq t_{j-1}$
such that $\RW(s)\in B_1$ and $t_j$ to be the first time $t\geq s_j$
such that $\RW(t)\notin B_2$.  Let $\tau$ be the first time $t\geq 0$
such that $\RW(t)\in\Paq$.  Finally, for each $s\geq 0$ let
$\RW^s$ be the restriction of $\RW$ to the interval $t\in[0,s]$.

For each $j=1,2,\dots$, we consider several events depending on
$\RW^{t_j}$ and $\Paq$.  Let $\ev Y_j$ be the
event that $\LE(\RW^{t_j})$ has a $(z_0,r,\eps)$-quasi-loop.
Let $\ev T_j$ be the event that $\tau> t_j$.
It is easy to see that 
$$
\ev A' \subset \bigcup_{j=1}^\infty (\ev Y_j \cap \ev T_j)
\,.
$$
Since $\ev T_j\supset\ev T_{j+1}$ for each $j$, this implies 
$$
\ev A' \subset \ev T_{m+1} \cup \bigcup_{j=1}^m \ev Y_j
\,.
\label e.decomp
$$
for every $m$.

We first estimate
$\P[\ev T_{j+1} \mid \RW^{t_j},\ \Paq]$.
Conditioned on any $\Paq$,
the probability that a SRW starting at any vertex outside of
$B_2$ will hit $\Paq$ before hitting $B_1$ is at least
$$ 
\nco 1 (\eps_1)^{\nco 2}
\,,
$$
where $\co 1>0$ depends only on $c$ and $r$, and
$\co 2>0$ is an absolute constant.
This is based on the fact that
$\Paq$ is connected, contains $\one$, and has
diameter at least $\eps_1-\delta$.
Applying this to the walk $\RW$ from time $t_j$ on, we therefore get
$$
\P[\ev T_{j+1} \mid \RW^{t_j},\ \Paq] \leq 1 - \co 1 (\eps_1)^{\co 2}
\,.
$$
By induction, we therefore find that
$$
\P[\ev T_m\mid\Paq] 
\leq \Bigl( 1 - \co 1 (\eps_1)^{\co 2}\Bigr)^{m-1}
\,.
\label e.hitp
$$

We now estimate $\P[\ev Y_{j+1}\mid \neg\ev Y_j, \RW^{t_j}]$.
Let $Q_j$ be the set of components of $\LE(\RW^{s_{j+1}})\cap B_2$
that do not contain $\RW(s_{j+1})$.
Observe that for $\ev Y_{j+1}$ to occur, there must
be a $K\in Q_j$ such that the random walk $\RW$ comes at some
time $t\in[s_{j+1},t_{j+1}]$ within distance $\eps$
of $K\cap B_1$ but $\RW(t)\notin K$ for all $t\in[s_{j+1},t_{j+1}]$.
But if $\RW(t)$ is close to $K$, $t\in[s_{j+1},t_{j+1}]$,
then \ref l.harmm/ can be applied, to estimate the probability
that $\RW$ will not hit $K$ before time $t_{j+1}$.
That is,
conditioned on $\RW^{s_{j+1}}$, for each given $K\in Q_j$,
the probability that
$\RW([s_{j+1},t_{j+1}])$ gets to within distance $\eps$ of
$K$ but does not hit $K$ is at most
$\nco 3 (\eps/r)^{\nco 4}$, where $\co 3,\co 4>0$ are absolute constants. 
Consequently, we get,
$$
\P[\ev Y_{j+1}\mid \neg\ev Y_j, \RW^{t_j}]
\leq \co 3|Q_j|(\eps/r)^{\co 4}
\,.
$$
Observe that $|Q_j|$, the cardinality of $Q_j$, is at most $j$.
Therefore,
$$
\P[\ev Y_{j+1}\mid \neg\ev Y_j]
\leq  \co 3 j(\eps/r)^{\co 4}
\,.
$$
This gives
\begineqalno
\P\Bigl[\bigcup_{j=1}^m\ev Y_j\Bigr]
\leq \sum_{j=1}^{m-1} \P[\ev Y_{j+1} \cap\neg\ev Y_{j}]
&
\leq \sum_{j=1}^{m-1} \P[\ev Y_{j+1} \mid \neg\ev Y_{j}]
\cr &
\leq \sum_{j=1}^{m-1} j \co 3 (\eps/r)^{\co 4}
\leq \nco 5 m^2 (\eps/r)^{\co 4}
\label e.yest
\endeqalno 
Combining this with \ref e.cas/, \ref e.pc/,
\ref e.decomp/ and \ref e.hitp/, we find that
$$
\mun\bigl[{\ev A(z_0,r,\eps)}\bigr]
\leq
 \co {aa} (\eps_1/c)^{\co {ab}} +
\co 5 m^2 (\eps/r)^{\co 4} +
\Bigl( 1 - \co 1 \eps_1^{\co 2}\Bigr)^{m-1}
\,.
$$
The lemma follows by taking
$m := \lfloor\eps^{-\co 4/3}\rfloor$
and $\eps_1 :=-1/\log\eps$, say.
\Qed

\procl t.ql
Let $(X,d)$ be a compact metric space, let $\aa,\bb\in X$, let $f:(0,\infty) \to (0,\infty)$ be monotone
increasing and continuous, and let $\Gamma=\Gamma(f)$ be the set of all compact simple paths
$\gamma\subset X$ with endpoints $\aa$ and $\bb$ which
satify the following property.  Whenever $x,y$ are points in $\gamma$ and $D(x,y)$ is the
diameter of the arc of $\gamma$ joining $x$ and $y$ we have
$$
d(x,y) \geq f(D(x,y))\,.
\label e.qlc
$$
Then $\Gamma$ is compact in the Hausdorff metric.
\endprocl

For this we will need the following Janiszewski's~\ref b.Janiszewski/
topological characterization of $[0,1]$ (see \crfl{Newman:planesets}{IV.5}):

\procl l.topchar \procname{Topological Characterization of Arcs}
Let $K$ be a compact, connected metric space, and let $\aa,\bb\in K$. 
Suppose that for every $x\in K-\{\aa,\bb\}$ the set $K-\{x\}$ is disconnected.
Then $K$ is homeomorphic to $[0,1]$.
\Qed
\endprocl

\proofof t.ql
Let $\Ha=\Ha(X)$ denote the space of compact nonempty subsets of $X$ with the
Hausdorff metric $\dHa$.
Let $\gamma$ be in the closure of $\Gamma$ in $\Ha$.
Then $\gamma$ is connected, compact, and $\gamma\supset\{\aa,\bb\}$.
We now use \ref l.topchar/ to show that $\gamma$ is a simple path.
Indeed, suppose that $x\in\gamma-\{\aa,\bb\}$.

\comment{
We first show that $\gamma-\{x\}$ has at most two components.  Let $y\in\gamma$.
Let $B$ be an open ball about $x$ with radius ${1\over 2}f\bigl(d(x,y)/2\bigr)$.
Let $\{\gamma_n\}$ be a sequence in $\Gamma$ such that $d_{\Ha}(\gamma_n,\gamma)<1/n$.
Then there are sequences $x_n,y_n\in\gamma_n$ with $x_n\to x$ and $y_n\to y$
as $n\to\infty$.  Let $N_\aa$ be the set of $n$'s such that  $y_n$ is on the arc
of $\gamma_n$ joining $\aa$ to $x_n$, and let $N_\bb$ be the set of $n$'s such that
$y_n$ is on the arc of $\gamma_n$ joining $x_n$ to $\bb$.
Then $N_\aa\cup N_\bb$ contains the positive integers, so at least one of these
sets is infinite.  Suppose, for example, that $N_\aa$ is infinite.
Let $\gamma_n'$ be the closed arc of $\gamma_n$ joining $y_n$ and $\aa$.

Suppose that $z_n\in\gamma_n\cap B$ and $x_n\in B$. 
Then the arc in $\gamma_n$ joining $z$ and $x_n$ has diameter at most
$d(x,y)/2$, since $\gamma_n\in \Gamma(f)$.
Consequently, for large $n$ it cannot contain $y_n$.
Therefore, for large $n$ the arc $\gamma_n'$ does not intersect $B$. 
Let $\gamma'$ be a Hausdorff limit of a subsequence of $\gamma_n'$.
Then $\gamma'$ is connected, disjoint from $B$, and contains $\aa$ and $y$.
So $y$ and $\aa$ are in the same connected component of $\gamma-\{x\}$.
Similarly, if $N_\bb$ is infinite, then $y$ and $\bb$ are in the same connected
component of $\gamma-\{x\}$.
We conclude that each connected component of $\gamma-\{x\}$ intersects $\{\aa,\bb\}$,
and in particular there are at most two components.
}

We show that $\aa$ and $\bb$ are in distinct components of $\gamma-\{x\}$.
Let $\{\gamma_n\}$ be a sequence in $\Gamma$ such
that $\dHa(\gamma_n,\gamma)<1/n$, and let $\{x_n\}$ be
a sequence with $x_n\in\gamma_n$ and $d(x_n,x)<1/n$.
For each $n$ let $\gamma^\aa_n$ be the closed arc of
$\gamma_n$ with endpoints $\aa$ and $x_n$, and let $\gamma^\bb_n$ be the
closed arc of $\gamma_n$ with endpoints $x_n$ and $\bb$.
By passing to a subsequence, if necessary, assume with no
loss of generality that the Hausdorff limits
$\gamma^\aa=\lim \gamma^\aa_n$ and $\gamma^\bb=\lim\gamma^\bb_n$ exist.
If $p_n\in\gamma^\aa_n$, $q_n\in\gamma^\bb_n$, then
$d(p_n,q_n)\geq f\bigl(d(p_n,x_n)\bigr)$, since $\gamma_n\in\Gamma(f)$.
By taking limits we find that if $p\in \gamma^\aa$ and
$q\in\gamma^\bb$, then $d(p,q)\geq f\bigl(d(p,x)\bigr)$.
Consequently, $\gamma^\aa-\{x\}$ and $\gamma^\bb-\{x\}$ are disjoint.
Because $\gamma^\aa\cup\gamma^\bb=\gamma$ and $\gamma^\aa,\gamma^\bb$
are compact, the set $\gamma-\{x\}$ is not connected.
Hence, by \ref l.topchar/, $\gamma$, is a simple path.

It remains to prove that $\gamma\in\Gamma(f)$.
For any simple path $\beta\subset X$ with endpoints $\aa,\bb$, let
$R(\beta)$ be the set of all $(x,y)\in\beta\times\beta$ such
that $x$ belongs to the subarc of $\beta$ with endpoints $\aa$
and $y$.  With arguments as above,  it is not hard to show that
$\lim R(\gamma_n)=R(\gamma)$, in the Hausdorff metric on $X\times X$.
Since $f$ is continuous, it then easily follows that
$\gamma\in\Gamma(f)$.  The details are left to the reader.
(Actually, one can see that this statement is not essential for the proof of
\ref t.noloop/.  There, we only need the fact that
the Hausdorff closure of $\Gamma(f)$ is contained
in the set of simple paths.)
\Qed

Although this will not be needed here, we note the following
variation on \ref t.ql/.

\procl t.qlg
Let $(X,d)$ be a compact metric space, let $f:(0,\infty) \to (0,\infty)$ be
monotone increasing and continuous, and let $\Gamma'=\Gamma'(f)$ be the set of
all compact subsets of $\gamma\subset X$ that are simple paths satifying the
following property.  Whenever $x,y$ are points in $\gamma$ and $D(x,y)$ is the
diameter of the arc of $\gamma$ joining $x$ and $y$ we have
$$
d(x,y) \geq f(D(x,y))\,.
\label e.qlc
$$
Set $\Gamma_0:=\bigl\{\{x\}\st x\in X\bigr\}$.
Then $\Gamma'\cup\Gamma_0$ is compact in the Hausdorff metric.
\endprocl

The proof does not require much more than the proof of \ref t.ql/.
We omit the details.

\proofof t.noloop
We start with the proof of the second statement,
and first assume that $a,b\neq\infty$.
Let $\zer$ and $\one$ be vertices of $\dZ$ closest to $a$ and
$b$, respectively.  Let $\Pa_\delta$ be LERW from $\zer$ to $\one$
in $\dZ$, and let $\Pa'_\delta$ be a path from $a$ to $b$,
obtained by taking the line segment joining $a$ to a closest
point $a'$ on $\Pa_\delta$, taking the line segment
joining $b$ to a closest point $b'$ on $\Pa_\delta$,
and taking the path in $\Pa_\delta$ joining $a'$ and $b'$.
Then the Hausdorff distance from $\Pa_\delta$ to $\Pa'_\delta$
is less than $2\delta$.
Let $\mun'$ be the law of $\Pa'_\delta$.
Let $\mu$ be some subsequential weak limit of $\mun$ as
$\delta\to 0$.  Then it is also a subsequential scaling
limit of $\mun'$.  Let $m$ be large.
By \ref l.tightdiam/, there is an $R_m$ such that with probability
at least $1-2^{-m-1}$ we have $\Pa_\delta\subset B(\aap,R_m)$.
For each $j\in\N$,
let $z^j_1,\dots,z^j_{k_j}$ be a finite set of points in $\R^2$ such
that the open balls of radius $2^{-j-2}$ about these points cover
$\overline B(\aap,R_m)$.
For each $j\in\N$, let $\eps^m_j\in(0,1)$ be sufficiently small so that
$$
\mun\bigl[ \ev A(z_i^j,2^{-j},\eps^m_j)\bigr]< 2^{-m-2-j}/k_j
\,,
\label e.rare
$$
for all $i=1,\dots, k_j$, and for all $\delta>0$.
Such $\eps_j^m$ exist, by \ref l.noloop/.
Finally, let $f_m:(0,\infty)\to(0,\infty)$
be a continuous monotone increasing function satisfying
$f_m(2^{2-j})\leq 2^{-3}\min(\eps_j^m,2^{-j})$ for each $j=1,2,\dots$
and $\sup_{s>0} f_m(s)\leq 2^{-3}\min(\eps_1^m,1/2)$.

Let $\ev X_m$ be the space of all compact nonempty
subsets of $\overline B(\aap,R_m)$, and set
$$
\ev Q_m := \ev X_m - \bigcup_{j=1}^\infty\bigcup_{i=1}^{k_j}
\ev A(z_i^j,2^{-j-1},\eps^m_j)
\,.
$$
Note that
$$
\mun'(\neg \ev Q_m) \leq 2^{-m}
\,,
\label e.musup
$$
for all $\delta$ and $m$.
Also note that if we set $X:=\overline B(\aap,R_m)$
 and $f:=f_m$ in \ref t.ql/, then 
$$
\ev Q_m\subset\Gamma_m:=\Gamma(f_m)
\,.
\label e.contain
$$
Indeed, suppose that $\gamma\in \ev X_m$ is a path in $X$ joining $a$ and $b$
and $\gamma\notin \Gamma_m$.
Then there are $x,y\in\gamma$ and $w$ contained in the arc
of $\gamma$ joining $x$ and $y$ such that $\dist(x,y)< 2f(\dist(x,w))$.
Since $\sup f\leq 2^{-4}$, we have $\dist(x,y)< 2^{-3}$.
Let $j$ be such that $2^{-j+1}\leq \dist(x,w)\leq 2^{-j+2}$. 
Then
$$
\dist(x,y)<2f(\dist(x,w))\leq 2f(2^{-j+2})\leq 2^{-j-2},\qquad
\dist(x,y)<2^{-3}\eps_j^m
\,.
$$
Let $i\in\{1,\dots,k_j\}$ be such that $\dist(x,z^j_i)\leq 2^{-j-2}$.
Then $\dist(z^j_i,w)\geq 2^{-j}$ and $x,y\in B(z^j_i,2^{-j-1})$.
Consequently, since $d(x,y) < \eps_j^m$, we have $\gamma\in \ev  A(z_i^j,2^{-j-1},\eps_j^m)$.
This proves \ref e.contain/.

By \ref e.musup/ and \ref e.contain/, we get $\mun'(\neg\Gamma_m)\leq 2^{-m}$.
\ref t.ql/ tells us that $\Gamma_m$ is compact.
Therefore, we also have $\mu(\neg\Gamma_m)\leq 2^{-m}$, and so 
$$
\mu\left(\neg\bigcup_m \Gamma_m\right)=0
\,.
$$
This completes the proof for the case $a,b\neq\infty$, because
each element of $\bigcup_m\Gamma_m$ is a simple path.

The proof when $a$ or $b$ in $\infty$ is similar.
One only needs to note that \ref l.noloop/ is
valid when $z_0=\infty$ and the distances are measured in the spherical
metric.  Indeed, the basic harmonic measure estimate \ref l.harmm/
is also valid in the context of the spherical metric.

The proof of the first statement of \ref t.noloop/ is also similar.
The details are left to the reader.
\Qed

Note that the first statement of \ref t.noloop/ implies that
for almost every subsequential scaling limit path
$\gamma$ from $a$ to $\bd D$, the closure of $\gamma\cap D$
intersects $\bd D$ in a single point.  This fact is easy
to deduce directly, since it is also true for the image of
SRW starting near $a$ and stopped when $\bd D$ is hit.

\bsection {First steps in the proof of \ref t.lesl/}{s.markov}
\resetConsts

Throughout this section we assume \ref g.confinv/.
Let $\lsl$ be random, with the law of the scaling limit of
LERW from $0$ to $\bd \U$.  
{}From \ref t.noloop/ we know that $\lsl$ is a.s.\ a simple path.

Recall that if $D\subset D'\subsetneqq\C$ are simply connected
domains with $0\in D$, then the conformal radius of $D'$
is at least as large as the conformal radius of $D$.  This
follows from the Schwarz Lemma applied to the map
$f_{D'}^{-1}\circ f_D:\U\to\U$.

For each $t\in(-\infty,0]$,
let $\lsl_t$ be the subarc of $\lsl$ with one endpoint
in $\bd \U$ such that
the conformal radius of $\U-\lsl_t$ is $\exp t$.
It is clear that $\lsl_t$ varies continuously in $t$.
Let $f_t:\U\to\U-\lsl_t$ be the conformal map
satisfying $f_t(0)=0$ and $f_t'(0)=\exp t$.
By \lowner{}'s slit mapping theorem \ref b.Lowner/,
there is a unique continuous $\zeta=\zeta_\lsl:(-\infty,0]\to\bd\U$ 
such that the differential equation \ref e.loew/ holds.
Let $\hatzet=\hatzet_\lsl:(-\infty,0]\to\R$ be the
continuous function satisfying $\zeta(t)=\exp\bigl(i\hatzet(t)\bigr)$
and $\hatzet(0)\in[0,2\pi)$.  Our goal is to prove

\procl p.redu 
The law of $\hatzet$ is stationary, and $\hatzet$ has independent increments.
\endprocl

This means that for each $s<0$ 
the law of the map $t\mapsto\hatzet(s+t)$ restricted to
$(-\infty,0]$ is the same as the law of $\hatzet$,
and that for every $n\in\N$ and $t_0\leq t_1\leq\cdots\leq t_n\leq 0$,
the increments
$\hatzet(t_1)-\hatzet(t_0),\hatzet(t_2)-\hatzet(t_1),\dots,
\hatzet(t_n)-\hatzet(t_{n-1})$ are independent.
The proof of this proposition, as well as the next, will be completed in
later sections.

Note that \ref g.confinv/ implies that the distribution of $\lsl$ is invariant
under rotations
of $\U$ about $0$.  Let $\lsl^1$ be random with the
law of $\lsl$ conditioned to hit $\bd U$ at $1$. 
If $\lambda$ denotes the (random) point
in $\lsl\cap\partial \U$, then $\lsl^1$ has the same law as
$\lambda^{-1}\lsl$.  
It turns out that \ref p.redu/ will follow quite easily from 

\procl p.prod
Assume \ref g.confinv/.  Fix some $t<0$.
Take $\lsl^1$ and $\lsl$ to be independent.
As above, let $\lsl_t$ be
the compact arc of $\lsl$ that has one endpoint on $\partial \U$ and
such that the conformal radius of
$\U-\lsl_t$ is $\exp(t)$.  Let $q(t)$ be the endpoint of $\lsl_t$ that
is in $\U$.  Let $\phi$ be the conformal map from $\U$ onto
$\U-\lsl_t$ satisfying $\phi(0)=0$ and $\phi(1)=q(t)$.
Then $\lsl_t\cup\phi(\lsl^1)$ has the same law as $\lsl$.
\endprocl

Now comes an easy lemma about LERW, and \ref p.prod/ will be
obtained from this lemma by passing to the scaling limit.  The passage
to the scaling limit is quite delicate.  Recall the definition of
the graph $\G(D,\delta)$ approximating a domain $D$, from
\ref s.back/.

\procl l.dprod
Let $\delta>0$ and $t<0$ be fixed, and let $D\subsetneqq\C$
be a simply connected domain with  $0\in D$.
Let $\beta$ be LERW from $0$ to $\bd D$ in $\G(D,\delta)$.
Let $\beta_t$ be the compact arc in $\beta$ such that $\beta_t\cap\partial D$
is an endpoint of $\beta_t$ and such that the conformal radius of
$D-\beta_t$ is $\exp(t)$.  Let $q_1(t)$ be the
endpoint of $\beta_t$ that is not on $\partial\U$.
Set $D_t=D-\beta_t$.
Then the law of $\beta-\beta_t$ conditioned on $\beta_t$
is equal to the law of LERW from $0$ to $\bd D_t$, conditioned
to hit $q_1(t)$.
\endprocl

\proof
There are several different ways to prove this lemma.
We prove it using the relation between LERW and the UST.
Suppose that $\alpha$ is a path such that $\beta_t=\alpha$
has positive probability. 
We assume for now that the endpoint $q$ of
$\alpha$ which is in $D$ lies in the relative interior of an
edge $e$ of $\dZ$  (this must be true except for at most a countable
possible choices of $t$), and set $\tilde\alpha:=\alpha-e$.
Let $\tilde q$ be the endpoint of $\tilde \alpha$ in $D$.
Let $\T$ be the UST on $\G^W(D,\delta)$, the wired graph of $D$.
Then $\beta$ may be taken as the path in $\T$ from $0$ to $\bd D$. 
We may generate $\T$ using
Wilson's algorithm with root $\bd D$, and starting with vertices
$v_1=\tilde q$ and $v_2=0$.  
Conditioning on $\beta_t$ being equal to $\alpha$
is the same as conditioning on $\alpha\subset \beta$,
which is the same as conditioning on the LERW $\LE_1$ from $v_1$ to
$\bd D$ to be $\tilde\alpha$ and that the LERW from $v_2$
to $\LE_1$ hits $\tilde q$ through the edge $e$.
This completes the proof in the case where $q$ is not a vertex
of $\dZ$.  The case where $q\in\verts(\dZ)$ is treated similarly.
\Qed

\bsection {Getting uniform convergence}{s.gu}
\resetConsts

The principle goal of this section is to state and prove
\ref c.unifsl/ below.
The main point there is that \ref g.confinv/ implies
that the weak convergence
of loop erased random walk in a domain $D\subset\U$ is uniform in $D$. 

\procl l.apxx
Let $K$ be a compact connected set in $\overline \U$ that contains
$\partial \U$ but with $0\notin K$.  Let $F$ be a compact subset
of $\U-K$, and let $\eps>0$. 
Then there is a $\delta_1>0$ with the following property. 
Let $K'$ be a compact connected set in $\overline \U$ with
$K'\supset\partial\U$ and $d_{\Ha(\U)}(K,K')<\delta_1/5$.
Set $D=\U-K'$. Let $\delta\in(0,\delta_1/5)$,
and let $Q\subset \Vba(D,\delta)$ be nonempty.
Let $\RW^Q$ be the random walk on $\G(D,\delta)$ starting at $0$ that
stops when it hits $\Vb(D,\delta)$, conditioned to hit $Q$.
Then the probability that $\RW^Q$ will reach $F$
after visiting some vertex within distance 
$\delta_1$ of $K$ is less than $\eps$.  
\endprocl

\proof
We need to recall some basic facts relating the conditioned
random walk $\RW^Q$ to the unconditioned random walk $\RW$ (that
stops when hitting $\Vb(D,\delta)$).
First recall that $\RW^Q$ is a Markov chain.
(This is easy to prove directly.  See also
the discussion of Doob's $h$-transform in \crfl{Durrett:BMM}{\SSS 3.1}.)
Let $v_0$ be some vertex, $t_0\in\N$,
and $W$ a set of vertices.  Let $\tau=\tau_W$ be the
least $t\geq t_0$ such that $\RW^Q(t)\in W$, if
such exists, and otherwise set $\tau=\infty$.
For $w\in W$ let
$$
a_W(v_0,w):=\P[\tau<\infty, \RW^Q(\tau)=w\mid\RW^Q(t_0)=v_0]
\,,
$$
and let $\tilde a_W(v_0,w)$ be the corresponding quantity for $\RW$.
Then 
$$
a_W(v_0,w) = {\tilde a_W(v_0,w) h(w)\over h(v_0)}
\,,
\label e.hest
$$
where $h(v)$ is the probability that $\RW$ hits $Q$
when it starts at $v$.  This formula is easy to verify.

Let $W_1$ be the set of vertices $v$ of $\G(D,\delta)$ such that
$h(v)< \eps h(0)/5$.
By \ref e.hest/,
$$
\sum_{w\in W_1} a_{W_1}(0,w)\leq (\eps/5)\sum_{w\in W_1} \tilde a_{W_1}(0,w)
\leq \eps/5
\,.
$$
Consequently, the probability that $\RW^Q$ visits $W_1$ is at most $\eps/5$.

Let $\rho$ be the distance from $F$ to $K$, and assume that $100\delta_1<\rho$.
Then an easy discrete Harnack inequality shows that 
$h(v)/h(0)<\nco a$ for all vertices $v$ in $F$,
where $\co a$ is some constant which does not
depend on $K'$ or $\delta$, but may depend on $F$.
There is a first vertex,
say $\hat v$, visited by $\RW^Q$ such that the distance from $\hat v$ to
$K$ is at most $\delta_1$.
Let $W_2$ be the set of vertices of $\G(D,\delta)$ in $F$.
If $\hat v\notin W_1$, then $h(\hat v)\geq \eps h(0)/5$ and hence
$$ 
\sum_{w\in W_2} a_{W_2}(\hat v,w) 
=\sum_{w\in W_2}\tilde a_{W_2}(\hat v,w) {h(w)\over h(0)}{h(0)\over h(\hat v)}
\leq 5\co a\eps^{-1} \sum_{w\in W_2} \tilde a_{W_2}(\hat v, w)
\,.
\label e.eee
$$ 
But since $K$ is connected, \ref l.harmm/ shows that
the probability that a simple random walk
starting at $\hat v$ will get to distance $\rho/2$ from $\hat v$ without
hitting $V_\partial(D,\delta)$ is bounded by 
$\nco b (\delta_1/\rho)^{\nco c}$, where $\co b,\co c>0$ are absolute
constants. 
By \ref e.eee/,
$$
\sum_{w\in W_2} a_{W_2}(\hat v,w) 
\leq 5 \co a \eps^{-1} 
\co b (\delta_1/\rho)^{\co c}
\,.
$$
Consequently, if $\delta_1$ is chosen sufficiently small,
the probability that $\RW^Q$ starting at $\hat v$ will hit $W_2$
is less than $\eps/5$, provided $\hat v\notin W_1$.
But $\P[\hat v\in W_1]\leq \eps/5$, since the probability that
$W_1$ is visited is at most $\eps/5$.
The lemma follows.
\Qed

In the following, we let $\RW^{D,\delta}$ denote SRW on
$\G(D,\delta)$ that stops when it hits $\Vb(D,\delta)$,
and let $\RW^{D,\delta,Q}$ denote $\RW^{D,\delta}$ conditioned
to hit $Q$, if $Q\subset\Vba(D,\delta)$. 
Suppose that $\nu$ is a probability measure on $\Vba(D,\delta)$ and $p$ is
random with law $\nu$, then
$\RW^{D,\delta,\nu}$ will denote $\RW^{D,\delta}$ conditioned
to hit $p$ given $p$.  In other words, the law of
$\RW^{D,\delta,\nu}$ is the convex combination of
the laws of the walks $\RW^{D,\delta,\{p\}}$, with coefficients $\nu(\{p\})$.

\procl l.doub
Assume \ref g.confinv/.  Let $D\subset\U$ be a Jordan domain with $0\in D$.
Let $\phi: D \to\U$ be the conformal homeomorphism from $D$ to $\U$ satisfying
$\phi(0)=0$, $\phi'(0)>0$.  Then as $\delta\to 0$ the law of the
pair $\bigl(\phi(\LE(\RW^{D,\delta})),\partial\U\cap\phi(\RW^{D,\delta})\bigr)$
tends weakly to the law of the pair $(\lsl,\partial\U\cap\lsl)$.
\endprocl

This lemma is easily proved using arguments as in the proof of
\ref l.apxx/, and is therefore left to the reader.

If $X$ and $Y$ are random closed subsets of $\U$,
we let $\dlru(X,Y)$ denote the Prohorov distance between
the law of $X\cup\bd\U$ and the law of $Y\cup\bd\U$
(see \ref s.back/, towards the end), where the metric
$\dd_{\br\U}$ is used on $\Ha(\br\U)$.
If $F$ is a subset of $\br\U$, we set
$\dlr_F(X,Y):=\dlru\left(X\cup\br{\U-F},Y\cup\br{\U-F}\right)$.
This is a measure of how much $X$ and $Y$ differ inside $F$.

\procl l.apxxx
Let $\eps>0$.
Then there is a $\delta_0>0$ and a finite collection of smooth Jordan domains
$D_1,D_2,\dots,D_n\subset\U$ with $0\in D_j$ for all $j$, and
with the following property.
Let $D\subset \U$ be a simply connected domain with $B(0,\eps)\subset D$,
let $\delta\in (0,\delta_0)$ and let $Q\subset \Vba(D,\delta)$ be nonempty.
Let $F_\eps$ be the component of $0$ in the set of points in $D$
that have distance at least $\eps$ to $\partial D$.
Then there is a $D'\in \{D_1,\dots,D_n\}$
and a probability measure $\nu$
on $\Vba(D',\delta)$
such that
$$
\dlr_{F_\eps}(\LE(\RW^{D,\delta,Q}),\LE(\RW^{D',\delta,\nu}))
<\eps
\,.
\label e.apxxx
$$
Moreover, we may require that the Hausdorff distance from
$\partial D$ to $\partial D'$ is at most $\eps$.
\endprocl

\proof  
Fix some $D$ and $Q$ as above.
By \ref l.apxx/, there is a $\delta_1>0$ such that with probability
$\geq 1-\eps/5$ the walk $\RW^{D,\delta,Q}$ does not reach
$F_\eps$ after exiting $F_{\delta_1}$, provided $\delta\in (0,\delta_1/5)$.
Also, there is a $\delta_1'<\delta_1$ such that with probability
$\geq 1-\eps/5$ this walk does not reach $F_{\delta_1}$ after
 exiting $F_{\delta_1'}$.
Let $D'\subset D$ be a smooth Jordan domain with $D'\supset F_{\delta_1'}$,
and such that the Hausdorff distance from $\partial D'$ to
$\partial D$ is less than $\eps$.
For every $\delta<\delta_1'/5$, let
$\nu=\nu_\delta$ be the hitting measure of $\RW^{D,\delta,Q}$ on
$\Vba(D',\delta)$.
Observe that we may think of $\RW^{D',\delta,\nu}$ as equal to
$\RW^{D,\delta,Q}$ stopped when $V_\partial(D',\delta)$ is hit.

Let $\ev A_1$ be the event that
$\RW^{D,\delta,Q}$ does not visit $F_\eps$ after exiting $F_{\delta_1}$,
and let $\ev A_2$ be the event that
$\RW^{D,\delta,Q}$ does not visit $F_{\delta_1}$ after exiting $F_{\delta_1'}$.
Note that on the event $\ev A_1\cap\ev A_2$, after exiting
$F_{\delta_1'}$ the random walk does not visit any vertex $v$ which
was already visited prior to the last visit to $F_\eps$.  Consequently,
the intersection of $F_\eps$ with the loop erasure of the walk
does not change after the first exit of $F_{\delta_1'}$.
Since we may couple $\RW^{D',\delta,\nu}$ to equal $\RW^{D,\delta,Q}$
stopped on
$V_\partial(D',\delta)$, this means that we may obtain a coupling giving
$F_\eps\cap \LE(\RW^{D,\delta,Q})=F_\eps\cap\LE(\RW^{D',\delta,\nu})$ on
$\ev A_1\cap A_2$.
Since $\P[\ev A_1\cap \ev A_2]\geq 1-\eps/2$, this proves the lemma for a 
single $D$.  However, the same solution would stand for every $D''$ with
$\partial D''$ sufficiently close to $\partial D$ in the Hausdorff metric.
Hence, the compactness of the
Hausdorff space of compact, connected subsets of $\overline\U-B(0,\eps)$
completes the proof.
\Qed

\procl l.unifhit
Let $D\subset\U$ be a smooth Jordan domain with
$0\in D$, let $\delta>0$,
let $q\in \Vba(D,\delta)$, and let $\RW^q:= \RW^{D,\delta,\{q\}}$.
Then the law of $\LE(\RW^q)$ is uniformly continuous in $q$.
That is, for every $\eps>0$ there is a $\delta_1>0$ such that
$$
\dlr_D(\LE(\RW^q),\LE(\RW^{q'}))<\eps
$$
provided $\delta\in(0,\delta_1)$, and $|q-q'|<\delta_1$.
\endprocl

\proof
Let $\delta_0>0$ be very small.  It is easy to see that when
$|q-q'|$ is small we may couple $\RW^q$ and $\RW^{q'}$ so that
with probability at least $1-\eps/5$ they are equal until they both
come within distance $\delta_0$ of $\partial D$.  Hence,
the lemma follows by using an argument similar to the one used in
the proof of \ref l.apxxx/.
\Qed

\procl c.unifsl
Assume \ref g.confinv/, and let $\eps>0$.  
Then there is a $\delta_1>0$ with the following property.
Let $D\subset\U$ be a simply connected domain with $B(0,\eps)\subset D$,
and let $F_\eps$ be the connected component of $0$ in the set of
all points $z$ with $d(z,\partial D)\geq\eps$.
Let $\delta\in (0,\delta_1)$, and
let $Q\subset \Vba(D,\delta)$ be nonempty.
Let $\phi$ be the conformal homeomorphism from $\U$ to $D$ that satisfies
$\phi(0)=0$ and $\phi'(0)>0$.  Then there is a random $\lambda\in\partial\U$
independent from $\lsl^1$ such that 
$$
\dlr_{F_\eps}\left(\phi(\lambda\lsl^1),\LE(\RW^{D,\delta,Q})\right)<\eps
\,.
$$
\endprocl

The main point here is that $\delta_1$ does not depend on $D$
or on $Q$.

\proof
Let $\eps_1>0$ be much smaller than $\eps$.
Suppose that $D'$ is a domain in the list
appearing in \ref l.apxxx/ that satisfies the requirements
there with $\eps_1$ in place of $\eps$, and let $\nu$ be as in that lemma.
Fix some small $\delta_1$.
For each $q\in \Vba(D',\delta)$ let
$\nu_q$ be restriction of the hitting measure
of $\RW^{D',\delta}$ to $B(q,\delta_1)\cap V_\partial(D,\delta)$,
normalized to be a probability measure.
\ref l.unifhit/ implies that we may replace $\nu$ by a probability measure
$\nu'$, which is a convex combination of such $\nu_q$, while having
$$
\dlr_{F_{\eps_1}}(\LE(\RW^{D',\delta,\nu}),\LE(\RW^{D',\delta,\nu'}))<\eps_1
\,,
\label e.cua
$$
provided $\delta_1$ is sufficiently small.
Moreover, by \ref g.confinv/, provided $\delta_1$ is sufficiently small
and $\delta\in(0,\delta_1)$, we have 
$$
\dlr_{F_{\eps_1}}(\LE(\RW^{D',\delta,\nu_q}),\psi(\lambda_q\lsl^1))\leq \eps_1
\,,
\label e.cub
$$
where $\lambda_q\in\partial\U$ is random and independent from $\lsl^1$,
and $\psi$ is the conformal homeomorphism from $\U$ to $D'$ satisfying
$\psi(0)=0$ and $\psi'(0)>0$.
Since the list $D_1,\dots,D_n$ in \ref l.apxxx/ is
finite, we may take $\delta_1$ to be independent of $D$.
Consequently, there is a random $\lambda\in\partial\U$ independent
from $\lsl^1$ with
$$
\dlr_{F_{\eps_1}}(\LE(\RW^{D',\delta,\nu'}),\psi(\lambda\lsl^1))\leq \eps_1
\,.
\label e.cuc
$$
Provided we have chosen $\eps_1$ sufficiently small,
we have that $|\psi(\phi^{-1}(z))-z|<\eps_1$ for
$z\in F_{\eps/2}$. 
\ref c.unifsl/ now follows from \ref e.apxxx/ with $\eps_1$
in place of $\eps$ and from \ref e.cua/, \ref e.cub/ and \ref e.cuc/.
\Qed


\bsection {Recognizing the \lowner\ parameter as Brownian motion}{s.cont}
\resetConsts

\proofof p.prod
Let $\beta:=\LE(\RW^{\U,\delta})$ and
let $\beta_t$, $D_t$ and $q_1(t)$ be defined as in \ref l.dprod/.
Set $\gamma:=\LE( \RW^{D_t,\delta,q_1(t)})$, where
$\RW^{D_t,\delta,q_1(t)}$ is taken to be independent from
$\beta$ conditioned on $\beta_t$.
Using \ref g.confinv/, 
$\dlru(\beta,\lsl)\to 0$ as $\delta\to 0$.
By \ref e.pcoup/, this means that we may couple
$\beta$ and $\lsl$ (that is, make them defined on the
same probability space, where they are not necessarily independent)
such that $\beta\toP \lsl$ in $\Ha(\U)$, where $\toP$ denotes
convergence in probability as $\delta\to0$.
Since $\lsl$ is a.s.\ a simple path, this also implies that
$\beta_t\toP\lsl_t$.


Let $\hat\phi$ be the conformal map 
from $\U$ onto $\U-\lsl_t$
that satisfies $\hat\phi(0)=0$ and $\hat\phi'(0)>0$, and
let $\hat\psi$ be the similarly normalized conformal map from
$\U$ onto $\U-\beta_t$.
Because $\beta_t\toP\lsl_t$, it follows that $\hat\psi\toP\hat\phi$,
in the topology of uniform convergence on compact subsets of $\U$.

Set $\hat\psi_\lambda(z):=\hat\psi(\lambda z)$
and $\hat\phi_\lambda(z):=\hat\phi(\lambda z)$ for $\lambda\in\partial\U$.
Given every $\eps>0$ and a closed set $A\subset\U$ let
$F_\eps(A)$ be the connected component of $0$ in the
set of points with distance at least $\eps$ from $A$ (or the empty
set, if $d(0,A)<\eps$),
and let $W_\eps(A):=\br{\U-F_\eps(A)}$.

By \ref c.unifsl/, for every $\eps>0$
there is a random $\lambda\in\partial\U$ independent
from $\lsl^1$ (but not from $\beta_t$) such that
$\dlr_{F_\eps(\beta_t)}(\hat\psi_\lambda(\lsl^1),\gamma)\to 0$.
 (The law of $\lambda$ may depend on $\delta$ and $\eps$.)
Observe that $\P[F_{2\eps}(\lsl_t)\subset F_{\eps}(\beta_t)]\to 1$
as $\delta\to 0$, because $\beta_t\toP\lsl_t$.
Therefore, we may conclude that
$\dlr_{F_{2\eps}(\lsl_t)}(\hat\psi_\lambda(\lsl^1),\gamma)\to 0$.
Since this is true for every $\eps>0$, it follows that we
may choose $\lambda=\lambda_\delta$ so that
$\dlr_{\U-\lsl_t}(\hat\psi_\lambda(\lsl^1),\gamma)\to 0$,
as $\delta\to 0$.
Because $\hat\psi_\lambda\toP\hat\phi_\lambda$, we therefore also have
$\dlr_{\U-\lsl_t}(\hat\phi_\lambda(\lsl^1),\gamma)\to 0$;
that is,
$\dlr_{\U}(\lsl_t\cup\hat\phi_\lambda(\lsl^1),\lsl_t\cup\gamma)\to 0$.
Since $\beta_t\cup\gamma$ has the same law
as $\beta$ (by \ref l.dprod/), and since $\beta_t\toP\lsl_t$, this gives,
$$
\dlr_\U(\lsl_t\cup\hat\phi_\lambda(\lsl^1),\beta)=
\dlr_\U(\lsl_t\cup\hat\phi_\lambda(\lsl^1),\beta_t\cup\gamma)\to 0
\,.
\label e.prodc
$$

Let $\lambda^*$ be random in $\partial\U$ with a law that
is some weak (subsequential) limit of the law of $\lambda$ as $\delta\to0$.
 It follows from \ref e.prodc/ that
$\lsl_t\cup\hat\phi_{\lambda^*}(\lsl^1)$ has the same law
as $\lsl$.  In particular, it is a simple path.
The only possibility is therefore that
$\hat\phi_{\lambda^*}=\phi$ a.s., which completes the proof
of \ref p.prod/.
\Qed

\proofof p.redu 
Recall that 
$\zeta=\zeta_\lsl:(-\infty,0]\to\bd\U$ is
the \lowner\ parameter associated to the LERW scaling limit
$\lsl\subset\U$.
Let $\tilde\zeta$ be the \lowner\ parameter associated
with the path $\lsl^1$, and let
$\tilde f_t$ be the associated solution of the \lowner\ system.
Note that $\tilde\zeta(0)=1$, since $\lsl^1\cap\bd\U=\{1\}$.
Fix some $t_0<0$.  Using 
\ref p.prod/ and its notations (with $t_0$ replacing $t$),
we know that the path $\check\lsl:=\lsl_{t_0}\cup\phi(\lsl^1)$
has the same law as $\lsl$.
Let $\check f_t$ be the solution of the \lowner\ system
associated with the path $\check\lsl$,
and let $\check\zeta$ be the associated \lowner\ parameter.
Then $\check\zeta$ has the same law as $\zeta$, by \ref p.prod/.
Let $\phi$ be as in \ref p.prod/, and set $\lambda := |\phi'(0)|/\phi'(0)$
When $t<t_0$, we have 
$$
\check f_t(z) = \phi\circ \tilde f_{t-t_0}(\lambda z)\,,
$$
because the right hand side is a suitably normalized conformal
map from $\U$ onto $\U-\check\lsl_{t}$.
We differentiate with respect to $t$, and use \ref e.loew/, to get
\begineqalno
{\bd \over\bd t} \check f_t(z)
&
= \phi'(\tilde f_{t-t_0}(\lambda z)) 
{\bd \over\bd t} \tilde f_{t-t_0}(\lambda z)
\cr &
= \phi'(\tilde f_{t-t_0}(\lambda z)) 
\lambda z\tilde f_{t-t_0}'(\lambda z){\tilde\zeta(t-t_0)+
\lambda z\over\tilde\zeta(t-t_0)-\lambda z}
= 
z\check f_t'(z){\lambda^{-1}\tilde\zeta(t-t_0)+
z\over\lambda^{-1}\tilde\zeta(t-t_0)-z}
\,.
\endeqalno
Consequently, it follows that
$\check \zeta(t) = \lambda^{-1}\tilde\zeta(t-t_0)$
for $t<t_0$.  It is clear that $\check\zeta = \zeta$ for $t\in[t_0,0]$.
Continuity of $\check \zeta$ gives $\lambda^{-1} =\zeta(t_0)/\tilde\zeta(0)=
\zeta(t_0)$.
Since $\zeta$ and $\tilde\zeta$ are independent, and $\tilde\zeta$
has the same law as $\zeta$ conditioned on $\zeta(0)=1$,
\ref p.redu/ follows.
\Qed

We shall need the following

\procl t.ito
Let $a(t)$, $t\geq 0$, be a real valued process (that is, a random
function $a:[0,\infty)\to\R$).
Suppose that $a$ is continuous a.s.\ and for every $n\in\N$ and
every $(n+1)$-tuple $0=t_0\leq t_1\leq t_2\leq\cdots\leq t_n$, the increments
$a(t_{j})-a(t_{j-1})$, $j=1,\dots,n$, are independent.
Then for every fixed $s_0\in(0,\infty)$, the random variable $a(s_0)$ is
Gaussian.
\endprocl

This theorem follows from the general theory of 
L\'evy processes. An entirely elementary proof can be found in
Section~4.2 of \ref b.Ito:lectures/.

\procl c.brows
There is a constant $c>0$ such that the process
$\zeta(t)$ has the same law as $\B(-c t)$,
where $\B(t)$ is Brownian motion on $\bd \U$ started at a uniform random
point.
\endprocl

\proof
That $\zeta(t)$ has the same law as $\B(-ct)$ for some
$c\geq 0$ follows immediately from \ref p.redu/ and \ref t.ito/.
The fact that $c>0$ is clear, since the LERW scaling limit
is not equal a.s.\ to a line segment.
\Qed

\bsection{The winding number of SLE}{s.wind}
\resetConsts

Let $\lecv\geq 0$, let $\B(t)$ be Brownian motion on $\bd \U$ started
at a uniform random point on $\bd\U$, and set 
$$
\zeta=\zeta_\lecv:=\B(-\lecv t)
\,.
\label e.zetadef
$$

\procl d.lep
Let $\KK$ denote the set of all
$\lecv\geq 0$ such that the \lowner{} evolution
$f_t$ defined by \ref e.zetadef/,
\ref e.loew/ and \ref e.bdval/ is a.s.\ for every $t<0$
a Riemann map to a slitted disk.  For $\lecv\in\KK$,
let $\lep_\lecv$ denote the (random) path defined by
$\lep_\lecv(t)=f_t(\zeta(t))$.  That is,
$\lep_\lecv$ is the path in $\br\U$ such that
$f_t$ is the nomalized Riemann map to 
$\U-\lep_\lecv\bigl([t,0]\bigr)$.

The random process $\lep_\lecv\bigl([t,0]\bigr)$, $t\leq 0$, will be called
{\bf stochastic \lowner\ evolution} (SLE) with constant $\lecv$.
\endprocl

As before, we let $\hat\B:[0,\infty)\to\R$ be the continuous map
satisfying $\B=\exp i \hat \B$ and $\hat\B(0)\in[0,2\pi)$.

\procl t.wind
Let $\lecv\in\KK$.  Let $T\leq 0$, and let
$\theta_\lecv(T)$ be the winding number of the path
$\lep_\lecv\bigl([T,0]\bigr)$ around $0$; that is
$\theta_\lecv(T)=\arg(\lep_\lecv(0))-\arg(\lep_\lecv(T))$,
with $\arg$ chosen continuous along $\lep_\lecv$.
Then for all $s>0$,
$$
\P\Bigl[\bigl|T - \log\left|\lep_\lecv(T)\right|\bigr|>s\Bigr] 
\leq \nco f\exp (-\nco e s)
\,,
\label e.1e
$$
and
$$
\P\Bigl[\bigl|\theta_\lecv(T)-\hat\B(0)+\hat\B(-\lecv T)\bigr|>s\Bigr] 
\leq \co f\exp (-\co e s)
\,,
\label e.2e
$$
where $\co f,\co e>0$ are constants, which depend only on $\lecv$.
\endprocl

Loosely speaking, the theorem says  that
$t+i \hat \B(-\lecv t)$ is a good approximation of the
path $\log \lep_\lecv(t)$.
A consequence of the theorem is that $\theta_\lecv(t)/\sqrt{\lecv t}$
converges to a gaussian of unit variance as $T\to-\infty$.

\proof
Let $f_t$ be defined by
\ref e.zetadef/, \ref e.loew/ and \ref e.bdval/.
Set $\lep:=\lep_\lecv$.
Let $w(t,z):= f_t^{-1}\bigl(f_T(z)\bigr)$,
and let $y=y(t,z):=\arg w(t,z)$, where
$\arg w$ is chosen to be continuous in $t$.

By \ref r.ode/, $w$ satisfies the
differential equation
$$
\pdt w =
-w {\B(-\lecv t)+w\over
\B(-\lecv t)-w}
\,,
\label e.zode
$$
where $\pdt$ denotes differentiation with respect to $t$.
Set $x=x(t,z):=\log|w(t,z)|$.  Then $w=\exp(x+iy)$, and
\ref e.zode/ can be rewritten, 
$$
\pdt x+i\pdt y =
{
\sinh x + i \sin\bigl(\hat\B(-\lecv t)-y\bigr) 
\over
\cosh x - \cos \bigl(\hat\B(-\lecv t)-y\bigr) }
\,.
\label e.ab
$$ 

Let $z_1$ be a random point on $\partial \U$, chosen
uniformly, and independent from the Brownian motion $B$.
Then $w(0,z_1)=f_T(z_1)$ is some point on the boundary
of $D_T:=f_T(\U)$.
Note that $\partial D_T$ is a connected set
that contains $\partial \U$ and intersects the
circle $\partial B\bigl(0,\exp(T)\bigr)$, by \ref e.schw/.
Set $A_s:=\bigl\{z\in \partial D_T\st |z| > \exp(T+s)\bigr\}$.
It follows from the continuous version of \ref l.harmm/ for Brownian motion
that the harmonic measure of
$A_s$ in $D_T$ at $0$ is bounded by $O(1)\exp(-\nco 1 s)$,
for some constant $\co 1>0$ and every $s\in\R$;
that is, at zero, the bounded harmonic function on $D_T$ that
has boundary values $1$ on $A_s$
and has boundary values $0$ on $\partial D_T-\overline A_s$
is bounded from above by $O(1)\exp(-\co 1 s)$.
Since harmonic measure is invariant under conformal maps,
we conclude that the measure of  $f_T^{-1}(A_s)$ is at most
$O(1)\exp(-\co 1 s)$.  This means that
$$
\P\bigl[\log |w(0,z_1)|-T>s\bigr] =
\P\bigl[\log |f_T(z_1)|-T>s\bigr] \leq O(1)\exp(-\co 1 s)
\,.
\label e.ww
$$

Now set $z_0=\B(-\lecv T)$.
Then $\lep(T)=w(0,z_0)$,
and so we need to relate $|w(0,z_0)|$
and $|w(0,z_1)|$.
Let $\tau$ be the least $t\in[T,0]$ such that
$w(t,z_1)=\B(-\lecv t)$,
if such a $t$ exists, and set $\tau=0$ if not.
Note that $|w(t,z_1)|=1$ while $t<\tau$, and $|w(t,z_1)|<1$
for $t\in(\tau,0]$. 
Also observe that conditioned on $\tau<0$, the law of
the process $\bigl(w(t,z_1)\st t\in [\tau,0]\bigr)$
is the same as the law of the process
$\bigl(w(t+T-\tau,z_0)\st t\in [\tau,0]\bigr)$.
Consequently, the random variable
$w(\tau-T,z_1)$ (where $\B$ is taken as two-sided Brownian motion
and \ref e.zode/ is extended to the range $t>0$), conditioned on $\tau<0$,
has the same distribution as the random variable
$w(0,z_0)$.
By \ref e.ab/, $\pdt x\leq 0$, and therefore 
$|w(\tau-T,z_1)|\leq |w(0,z_1)|$ on the event $\tau<0$.
Thus, for every $s\in\R$ we have
$\P\bigl[|w(0,z_1)|>s\mid\tau<0\bigr]\geq \P\bigl[|w(0,z_0)|>s\bigr]$.
Because $|w(0,z_1)|=1$ when $\tau=0$, we may drop the
conditioning on $\tau<0$.
Now \ref e.ww/ gives
$$
\P\bigl[\log |\lep(T)|-T>s\bigr] =
\P\bigl[\log |w(0,z_0)|-T>s\bigr] 
\leq O(1)\exp(-\co 1 s)
\,.
$$
On the other hand, the Koebe $1/4$ Theorem \ref e.k1q/ gives
$$
\exp T=|f'_T(0)|\leq 4 \inf\bigl\{|z|\st z\notin f_T(\U)\bigr\}
\leq 4|\lep(T)|\,,
$$
and so $\log\left|\lep(T)\right|+\log 4\geq T$ always.
This completes the proof of \ref e.1e/.

Now let $\tau_1$ be the least $t\in[T,0]$ such that
$x(t)=w(t,z_0)\leq -1$, and set $\tau_1:=0$ if such
a $t$ does not exist.  Since $x(t)$ is monotone
decreasing, we may write $y(t)$ as a function of $x$:
$y=g(x)$.
By \ref e.ab/,
$$
g'\bigl(x(t)\bigr) =
{ \sin\bigl(\hat \B(-\lecv t)-y(t)\bigr) \over \sinh x(t) }
$$
and hence
$$ \left|g'\bigl(x(t)\bigr)\right|\leq \bigl|\sinh x(t)\bigr|^{-1}
\,.
$$
And so we get 
$$
|y(0)-y(\tau_1)|
= \int_{x(0)}^{x(\tau_1)} |g'(x)|\,dx
\leq \int_{-\infty}^{-1} \bigl|\sinh x\bigr|^{-1}\,dx
<\infty
\,.
\label e.bdy
$$ 

Let $\phi(s,t):= f_t ^{-1} \bigl(\lep(s)\bigr)$
for $T\leq s\leq t\leq 0$.  Then $\phi$ is continuous and
its image does not contain $0$.  Hence, it may be considered
as a homotopy in $\C-\{0\}$ from the path
$\phi(s,0)=\lep(s)$, $s\in[T,0]$, to the
concatenation of the inverse of the path
$\phi(T,t)= f_t ^{-1}\bigl(\lep(T)\bigr)=w(t,z_0)$,
$t\in[T,0]$, with the path
$\phi(t,t)=\B(-\lecv t)$, $t\in[T,0]$.
Therefore, its winding number is the sum of the corresponding winding numbers.
This means that 
$$
\theta_\lecv(T) = \hat\B(0)-\hat\B(-\lecv T)+y(T)-y(0)
\,.
$$
By \ref e.bdy/, it therefore suffuces to prove the appropriate bound
on the tail of $|y(\tau_1)-y(T)|$.

Let $\ppn y:= \min\bigl\{|y-2 \pi n|\st n\in\Z\bigr\}$.
Set $t_0=T$, inductively, let $t_j$ be the first
$t\in[t_{j-1},0]$ such that $\pi/2 = \ppn {\hat\B(-\lecv t) - \hat\B(-\lecv t_{j-1})}$,
and set $t_j=0$ if no such $t$ exists.
Equation \ref e.ab/ shows that for every $s\in(T,0)$
$$
\pdt \ppn {y(t)-\hat\B(-\lecv s)}\leq 0,\qquad\hbox{at }t=s
\,.
$$
Consequently, for every $j\in\N$, if there is an $s\in[t_j,t_{j+1}]$
such that $\ppn{y(s)-\hat\B(-\lecv t_j)}<\pi/2$, then this is satisfied
also for all $s'\in[s,t_{j+1})$, because $y(t)$ cannot get out
of the set
$\Bigl\{p\in\R\st \bigl|{p-\hat\B(-\lecv t_j)}\bigr|_{2\pi}<\pi/2\Bigr\}$
while $\hat\B(-\lecv t)$ is in it.  This implies
that $|y(t_{j+1})-y(t_j)|<2\pi$.
Hence $|y(\tau_1)-y(T)|\leq 2\pi \min\{j\in\N \st t_j>\tau_1\}$.
Therefore, for every $a>0$ and $n\in\N$,
$$
\Pbig{|y(\tau_1)-y(T)|\geq 2\pi n}
\leq
\Pbig{\tau_1-T\geq a} + \Pbig{t_n\leq T+a}
\,.
\label e.casa
$$
The first summand on the right hand side is
bounded by $O(1)\exp(-\nco s\, a)$, for some constant $\co s>0$,
by \ref e.1e/.
To estimate the second summand, observe that
conditioned on $t_n\leq T+a$, we have probability at least
$2^{-n-1}$ for the event
$$
B\bigl(-\lecv\, (T+a)\bigr)-\hat\B(-\lecv T)\geq n\pi/2
\,,
\label e.evf
$$
because when $t_n\leq T+a$ and
$B\bigl(-\lecv\, (T+a)\bigr)\geq B\bigl(-\lecv t_n)
\geq\cdots\geq B\bigl(-\lecv t_0)$, we have \ref e.evf/.
However, \ref e.evf/ has probability
$$
(2\pi a\lecv)^{-1/2}
\int_{n\pi/2}^\infty \exp\bigl(-s^2/(2\lecv a)\bigr)\,ds
\leq O(1) \exp\bigl(-n^2/\nco n a\bigr)
\,,
$$
and hence
$$
\Pbig{t_n\leq T+a}
\leq O(1)  2^n\exp\bigl(-n^2/\co n a\bigr)
\,.
$$
We choose $a$ to be $n$ times a very small constant.
Then our above estimates, together with~\ref e.casa/,
give
$$
\Pbig{|y(\tau_1)-y(T)|\geq 2\pi n}
\leq O(1)\exp(-\nco m n)
\,,
$$
with the constants depending only on $\lecv$.
This completes the proof of the theorem.
\Qed

\bsection{The twisting constant of LERW}{s.windl}
\resetConsts 

Consider some scaling limit measure $\mmu$ of LERW from $0$
to $\bd\U$, and let $\gamma$
be random with law $\mmu$.
Assuming \ref g.confinv/, we have established that
SLE with some constant $\lecl$ has law $\mmu$.
In this section we show that $\lecl=\xplec$, and
thereby complete the proof of \ref t.lesl/.
 
Let $\eps\in(0,1)$, let $\gamma_\eps$ be the connected component of
$\gamma-B(0,\eps)$ which has a point in $\partial U$,
and let $\winding(\gamma_\eps,0)$ be the winding number of $\gamma_\eps$
around $0$, in radians.  That is, $\winding(\gamma_\eps,0)$ is
the imaginary part of $\int_{\gamma_\eps}z^{-1} dz$.
By symmetry, it is clear that $\E[\winding(\gamma_\eps,0)]=0$.
We shall show that
$$
\Ebig{ \winding(\gamma_\eps,0)^2}= 2\log (1/\eps) + O(1)\sqrt{\log (1/\eps)}\,.
\label e.wnv
$$
Based on this and the results of \ref s.wind/, it will follow that
$\lecl = \xplec$. 

The proof of \ref e.wnv/ will use Kenyon's work \ref b.Kenyon:conf/.
The overall idea of the proof is very simple, and based on the
relations between UST and domino tilings.
We now briefly review the relations between the UST on $\Z^2$ and
domino tilings, and the height function for domino tilings.
For a more thorough discussion, the reader should consult \ref b.Kenyon:conf/.

A {\bf domino tiling} of the grid $\Z^2$ is a tiling
of $\R^2$ by tiles of the forms $[k,k+1]\times[j,j+2]$
and $[k,k+2]\times [j,j+1]$, where $k,j\in\Z$.  
A domino tiling of $\Z^2$ may also be thought of as a
{\bf perfect matching} of the dual grid $\bigl((1/2)+\Z\bigr)^2$.  
(A perfect matching of a graph $G$
is a set of edges $M\subset\edges(G)$ such that every vertex
is incident with precisly one edge in $M$.)

Let us start with finite graphs.  
Let $D$ be a simply connected domain in $\R^2$ whose
boundary is a simple closed curve in the grid $\Z^2$,
and let $G:=\Z^2\cap\br D$.
Let $\rho_0$ be some vertex in $\bd D\cap\Z^2$, which we call
the {\bf root}.
Let $\check G$ be the graph $\bigl((1/2)\Z^2)\cap \br D$ with
$\rho_0$ and its incident edges removed.
Then there is a bijection, discovered by Temperley,
between the set of perfect matchings on $\check G$ and
spanning trees of $G$. 

Temperley's bijection (see \ref f.temperley/)
works as follows.  For every edge $[v,u]$ in the matching $M$
such that $v\in\Z^2$, we put in the tree the edge $e_u$ whose
center is $u$.  This gives the set of edges in the tree $T$.
If $[v,u]\in M$ is as above, we may orient the edge $e_u$
away from $v$, and then the tree $T$ will be oriented towards the root
$\rho_0$.

\midinsert
\SetLabels
\R\E(.24*-.01)$\rho_0$\\
\endSetLabels
\tabskip=1em plus2em minus.5em
\halign to \hsize{\hfil # \hfil&\hfil # \hfil&\hfil # \hfil\cr
\AffixLabels{\epsfxsize=1.5in\epsfbox{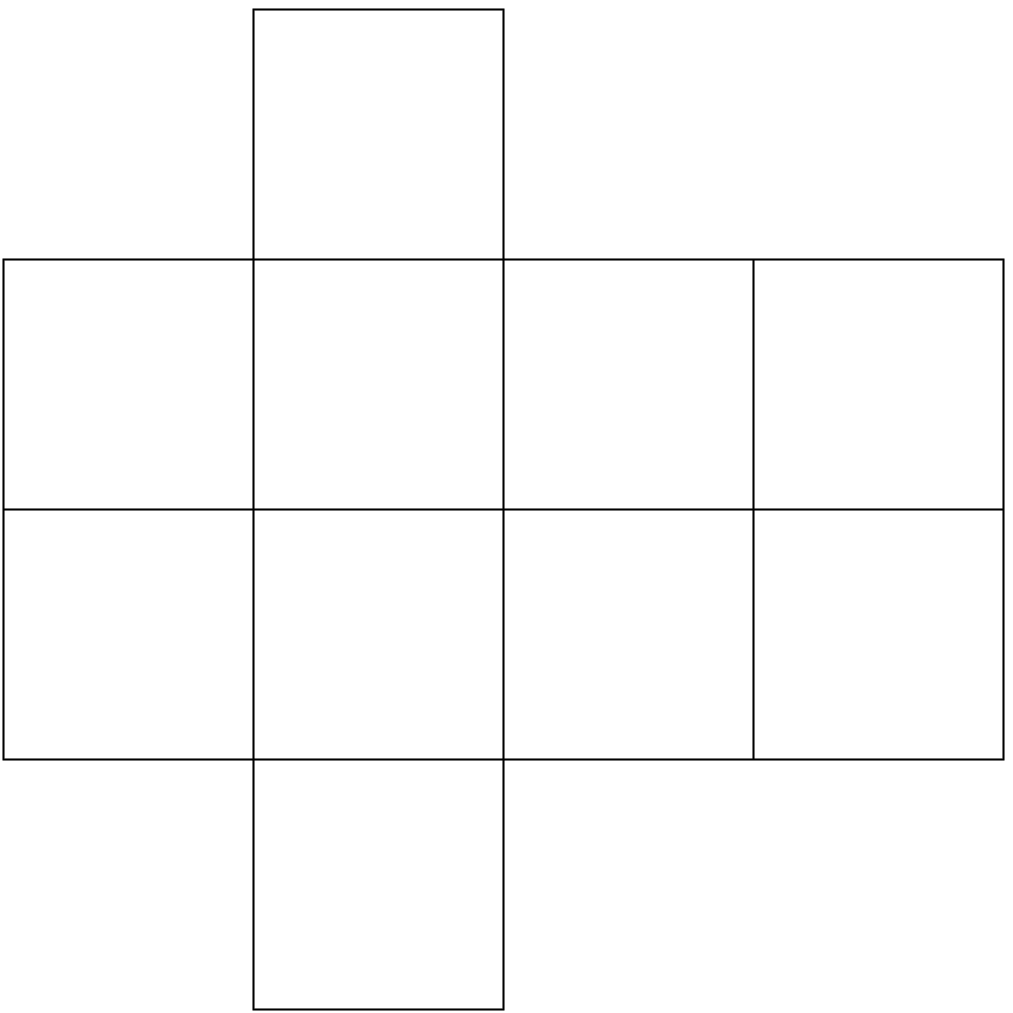}}&%
\epsfxsize=1.5in\epsfbox{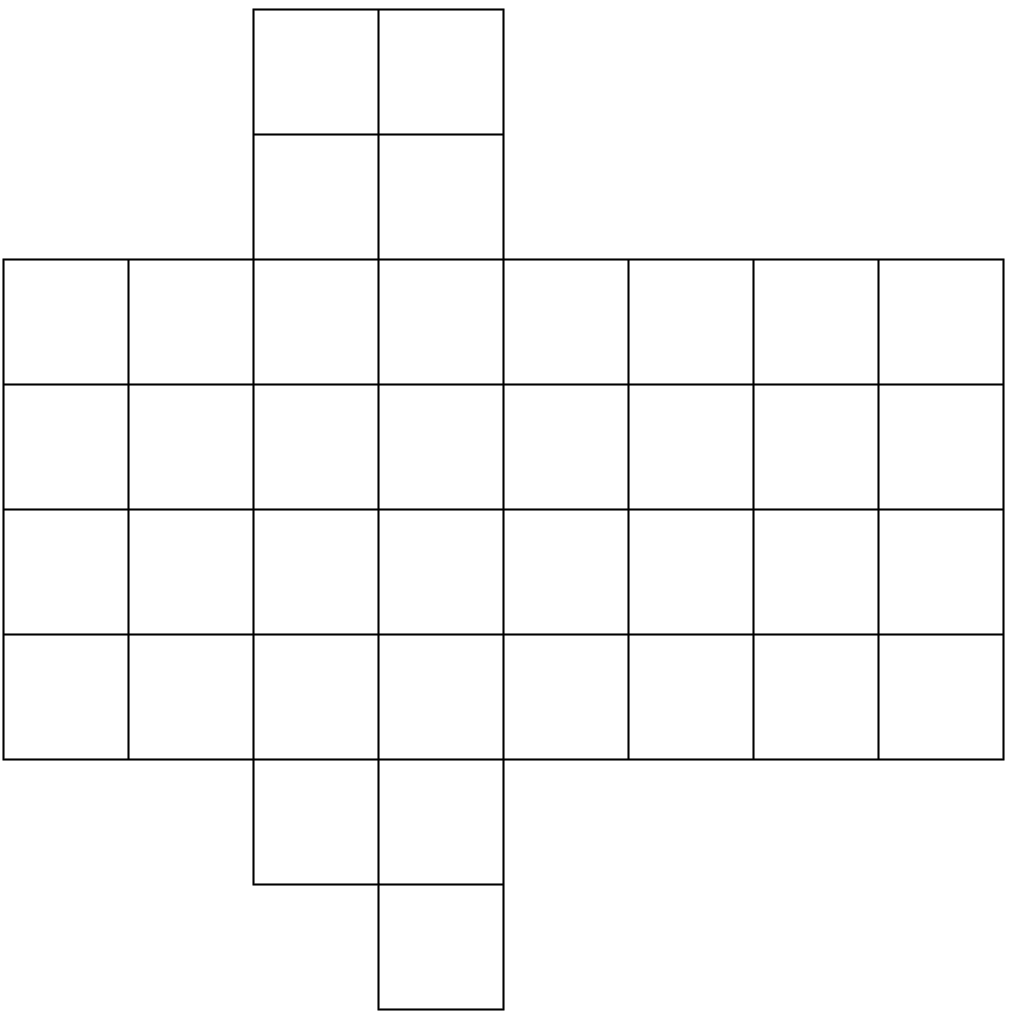}&%
\epsfxsize=1.5in\epsfbox{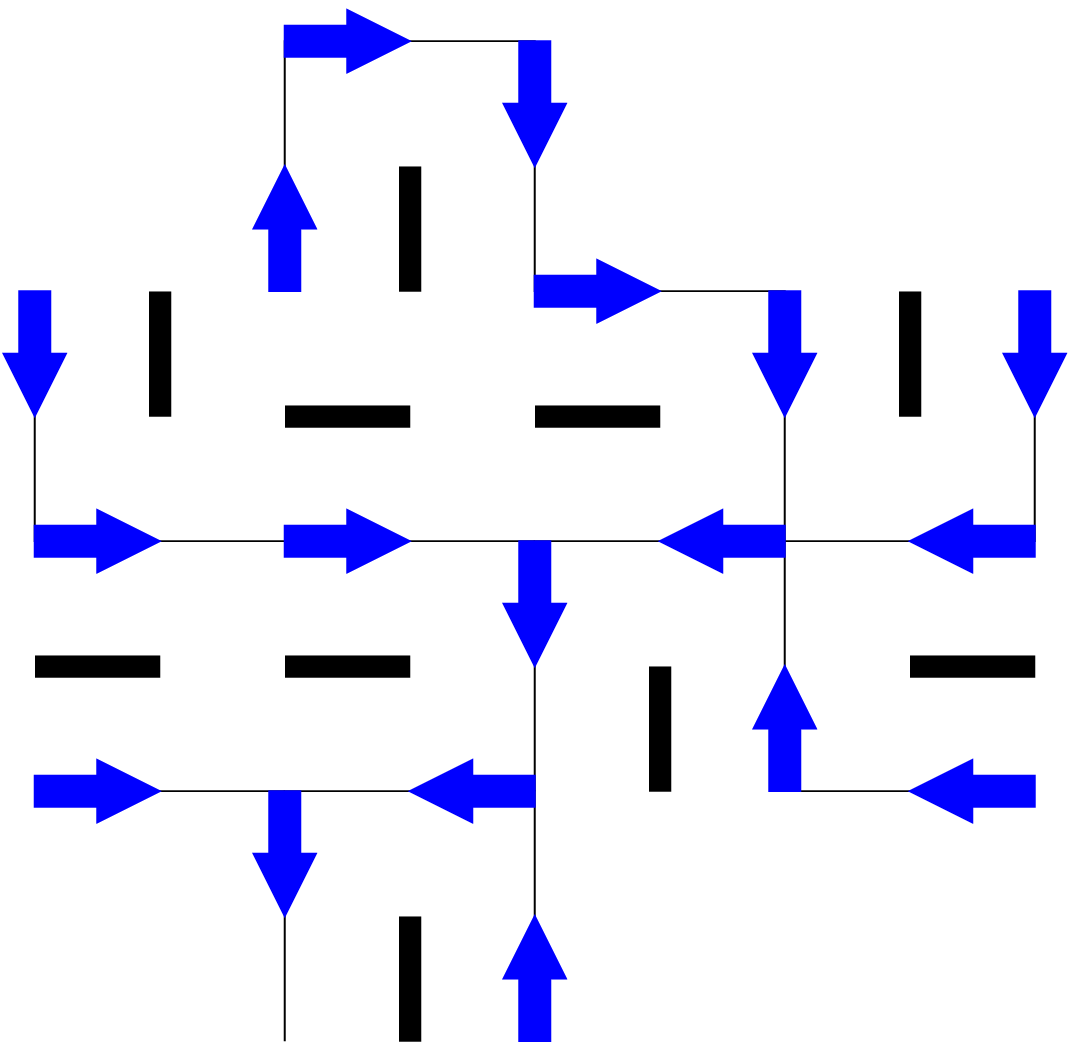}\cr
\noalign{\smallskip}
$G$&$\check G$&matching and tree\cr
}
\medskip
\narrower\noindent
\figlabel{temperley} Temperley's bijection. On the right,
the arrows are edges in the matching that containing vertices of $\Z^2$,
the solid segments are other edges in the matching, and the thin
lines are edges in the tree.
\endinsert

Temperley's bijection works also in more general situtations.
There is a simple modification to make it work for the wired
graph associated to the domain $D$.
Also, given a perfect matching on all of $(1/2)\Z^2$,
there is an associated (oriented) spanning forest of $\Z^2$.
The collection of all domino tilings of $\Z^2$ has a natural
probability measure (of maximal entropy), and for a.e.\ domino
tiling the corresponding spanning forest is a spanning tree.
Temperley's map from perfect matchings on $(1/2)\Z^2$ to
spanning forests of $\Z^2$ maps the cannonical probability
measure on the set of domino tilings to the law of the UST
of $\Z^2$.

Let $G$ and $\check G$ be as above, and let $\hat G$ be
the graph of the domino tiling, that is, the union
of the squares of edge length $1/2$ with centers at the
vertices of $\check G$, thought of as a subgraph
of the grid $(1/4,1/4)+(1/2)\Z^2$.  Associated to a domino
tiling of $\hat G$ is a {\bf height function} $h$ defined
on the vertices of $\hat G$.  
Here is the definition of $h$.  Pick some 
vertex $v_0\in\verts(\hat G)$ and some $a_0\in\R$,
and set $h(v_0)=a_0$.
Color a square face of the grid $(1/4,1/4)+(1/2)\Z^2$
{\bf white} if its center is 
a vertex of $\Z^2$ or if it is contained in a face of $\Z^2$,
and {\bf black} otherwise.
If $[u,v]\in\edges(\hat G)$ is on the boundary of a
domino tile in the tiling, then we require that $h(v)-h(u)=1$
if the square to the right of the directed edge $[u,v]$
is white and $h(v)-h(u)=-1$ if the square to the right of 
$[u,v]$ is black.  These constraints uniquely specify the height
function $h$ (except that the choices of $v_0$ and $a_0$ are
arbitrary).

We will work in the upper half plane $\H:=\{z\in\C\st \im z>0\}$.
Let $\Gd:=\G^W(\H,\delta)$, the wired graph of mesh $\delta$
associated with the domain $\H$, and
 $\hGd:=\H\cap\bigl((1/4,1/4)+(1/2)\Z^2\bigr)$.
The discussion above carries through for the grid
$\dZ$, in place of $\Z^2$.  (Although the distance between
adjacent vertices in the graph $\hat G$ is $\delta/2$ when
$G\subset\dZ$, we still work with the height function
where the height difference along an edge on the boundary
of a tile is $\pm1$.)
Temperley's bijection induces a measure preserving
transformation between domino tilings of the grid
$\hGd$ and the UST of $\Gd$. (If we keep the orientation,
then the UST is directed towards $\bd\H$.)   

We normalize the height function associated to a domino tiling
of $\hGd$ by requiring that $h\bigl((1/4,1/4)\bigr)=1/2$.
Then $h\bigl(((2k+1)/4,1/4)\bigr)=(-1)^k/2$, for $k\in\Z$.
If $v$ is some vertex in $\Gd$, let $h(v)$ be the average
of the value of $h$ on the vertices of $\hGd$ closest to $v$.

\procl l.windheight
Let $T$ be a spanning tree of $\Gd$, and let $h$ be the associated
height function.  Let $v\in\verts(\Gd)$ be a vertex different from
the wired vertex $\bd \H$, and let $a$ be
the real part of $v$.  Let $\ray_v'$ be the path from
$v$ to $\bd\H$ in $T$, considered as a path in the plane,
and let $\ray_v$ be the union of $\ray_v'$ with
the line segment joining the intersection $\ray_v'\cap\bd\H$ to $a$.
Then $-\pi h(v)/2=\winding(\ray_v,v)$, the winding number of
$\ray_v'$ around $v$.
\endprocl

This lemma is a special case of a more general observation
made by Kenyon.  
(Since $\ray_v$ is a path with $v$ as an endpoint,
we define
$\winding(\ray_v,v):=\lim_{r\to0}\winding\bigl(\ray_v-B(v,r),v\bigr)$;
which is the same as $\winding\bigl(\ray_v-B(v,\delta/2),v\bigr)$.)

\proof Use induction on the length of the path $\ray'_v$.
\Qed 

Symmetry implies that $\E[h(v)]=0$ for all $v\in\verts(\Gd)$.
Kenyon has shown \ref b.Kenyon:conf/ that
$$
\E[h(v)^2] =8\pi^{-2}\log(1/\delta)+O(1)
\label e.kenyon
$$
(provided that $v$ stays in a compact subset of $\H$).
Hence $\E[\winding(\ray_v,v)^2]=2\log(1/\delta)+O(1)$,
which seems very close to a proof of \ref e.wnv/.  However, to
make it into a proof of \ref e.wnv/ requires some effort (it seems).

The advantage of the height function over the winding number
is that the height difference between two vertices can be computed
along any path joining them.  On the other hand, to compute 
$\winding(\ray_v,v)$, one might think that it is necessary to follow
$\ray_v$, which is a random path. 
It is immediate that $h(v)$ is $\sum_j \lambda_j \chi_j$,
where $\chi_j$ is the indicator of the event that a certain
domino tile is present in the tiling, and $\lambda_j$ are some explicit
easy to compute (non-random) weights.  This means that to
calculate $\E[h(v)^2]$ one needs to have a good
estimate for the behavior of the correlations $\E[\chi_i\chi_j]$
for small $\delta$.  That's how Kenyon proves \ref e.kenyon/.

Recall that $A(p,r_1,r_2)$ denotes an annulus with center
$p$, inner radius $r_1$, and outer radius $r_2$.
The following result is an immediate consequence from \ref b.ABNW/.

\procl l.expo
Let $D\subset\C$ be a domain, and let $v_0\in\dZ\cap D$ be
some vertex.  Consider $T$, the UST on $\dZ\cap D$, with
free or wired boundary.  Let $r_2>2r_1>2\delta$,
and suppose that $r_2$ is smaller than the distance from
$v_0$ to $\partial D$.  Let $A$ be the annulus $A:= A(v_0,r_1,r_2)$,
and let $k(A)$ be the maximum number of disjoint paths in $T$ each
of which intersects both boundary components of $A$.
Then for each $k\in\N$
$$
\Pbig{k(A)\geq k}
\leq \nco z (r_2/r_1)^{\nco e (k-1)}
\,,
$$
where $\co z,\co e>0$ are universal constants. \Qed
\endprocl

\procl l.windtail 
Let $D\subset\C$ be some domain,
and consider $T$, the UST in $\G^W(D,\delta)$ (wired boundary).
Let $p_0\in D\cap \dZ$.
Given a set $K\subset D$, let
$X(K)$ denote the maximum winding number around $p_0$ of a path in $T$
with endpoints in $K$.
Let $r\in\bigl(0,d(p_0,\bd D)/2\bigr)$.
Then for each $h>0$ and $s\geq 2$,
$$
\P\Bigr[X\bigl(A(p_0,r/s,r)\bigr)>h\Bigl]\leq
\nco a \exp\bigl(-\nco b h/\log s \bigr)\log s
\,,
$$
where $\co a$ and $\co b$ are absolute constants.
\endprocl

\proof 
Set $A:= A(p_0,r/s,r)$, $X:= X(A)$.
To begin, assume that $s=2$.  Let $B_1,B_2,\dots,B_N$
be a covering of $A$ with balls of radius $r/10$, where $N\leq \nco N$,
with $\co N$ some universal constant.
Let $v_0\in A\cap\dZ$ be some vertex.
For any vertex $u$, let $W(u)$ be the signed winding 
number around $p_0$ of the path in $T$ from $v_0$ to $u$.
Consider neighbors $v,w$ in $A$, and let $\gamma$
be the UST path joining $v$ and $w$ in $T$.
If $\gamma\neq [v,w]$, then $\gamma\cup [v,w]$ is a simple
closed path, and hence has winding number at most
$2\pi$ around $p_0$. This implies that
$|W(v)-W(u)|$ is at most $2\pi$ plus the absolute
value of the winding number
of the edge $[v,u]$ around $p_0$, and therefore
$|W(v)-W(u)|\leq 3\pi$.   From this it follows that we may find
vertices $v_1,v_2,\dots,v_n$ in $A$ such that
$|W(v_j)-W(v_k)|\geq 3\pi$ for $j\neq k$ and $n\geq X/6\pi$.
Since $N\leq \co N$, we may find some ball $B_m$ from the
above collection, satisfying
$|B_m\cap \{v_1,\dots,v_n\}|\geq X/6\pi\co N$.
However, if $v,u\in B_m$ and $|W(v)-W(u)|\geq 2\pi$, then the
path in $T$ joining $v$ and $u$ must go around
$p_0$.  In particular, it must cross twice the annulus
$A_m:=A(c_m,r/10,r/5)$, where $c_m$ is the center of $B_m$.
It follows that the number of disjoint crossings of $A_m$
in $T$ is at least $|B_m\cap \{v_1,\dots,v_n\}|$.
Hence, by \ref l.expo/, for any fixed $m$ the probability that
$|B_m\cap \{v_1,\dots,v_n\}|>b$ is at most
$O(1)\exp(-\nco q b)$, for some constant $\co q>0$.
Consequently, 
$\P[X>h]\leq O(1)\co N\exp(-\co q h/6\pi\co N)$,
which completes the proof in the case $s=2$.

If $s>2$, then we may cover the annulus $A$ with
at most $2\log s+1$ disjoint concentric annuli with radii ratio $2$.
In order that $X(A)$ be at least $h$, there must be one
of these smaller annuli $A'$ with $X(A')\geq h/(2\log s +1)$.
The lemma follows.
\Qed

\procl l.rick \procname{{\rm\ref b.Kenyon:conf/}}
Let $p,q\in\H$ be any two points.
Consider a uniform domino tiling of the grid $\hGd\subset\H$ of mesh $\delta$,
and let $p',q'\in\verts(\hGd)$ be vertices closest to $p$ and $q$,
respectively.  Then
$$
\lim_{\delta\to0}\E [h(p')h(q')] =
8\pi^{-2}
\log\left|{\overline p-q\over p-q}\right|
\,.
\Qed
$$
\endprocl

It may be noted
that the right hand side is invariant under conformal automorphisms
of $\H$, since it is the log of the square root of a cross ratio of
$p,\overline p,q,\overline q$.

Let $\delta>0$, let $v_0$ be a vertex of the grid
$\delta\Z^2$ which is closest to $i=\sqrt{-1}$.
Fix some $r\in(0,1/2)$.
Let $\gamma_r^\delta$ be the connected component of
$\ray_{v_0}-B(v_0,r)$ that intersects $\partial\H$,
and let $W_\delta(r)$ be the winding number of $\gamma_r^\delta$ around
$v_0$.

\procl p.windw
Assuming $r<1/4$,
$$
\limsup_{\delta\to0}\Bigl|
\E\bigl[W_\delta(r)^2\bigr] - 2\log(1/r)\Bigr|
\leq \nco e \sqrt{\log(1/r)}
\,,
\label e.pwindw
$$
where $\co e$ is an absolute constant.
\endprocl 

\proof
Let $v_1$ be a vertex in $\delta\Z^2$ such that
$v_1-v_0-r/3\in[0,\delta)$.
Let $V_0$ be the set of vertices of $\delta\Z^2$ whose Euclidean distance to
$\{v_1,v_0\}$ is in the range $(r/9,2r)$.
Then, assuming that $\delta<r/9$,
$V_0$ separates $v_0$ from $v_1$ and $\{v_0,v_1\}$ from
$\partial\H$ in the grid $\delta\Z^2$.
Given a vertex $v$, let $\rayz_v$ be the union of
$\ray_v'$ with the line segment joining the intersection
$\ray_v'\cap\bd\H$ to $0$.
Set $\ray:=\bigcup\bigl\{\rayz_v\st v\in V_0\}$.

Set $X_j:= \winding(\ray_{v_j},v_j)$, $j=0,1$.
By \ref l.windheight/, have $X_j=-\pi h(v_j)/2$, and consequently,
\ref l.rick/ gives,
$$
\lim_{\delta\to 0}\E[X_0X_1] = 2\log(1/r) + O(1)
\,.
\label e.r
$$
Since $V_0$ separates $v_0$ from $v_1$, when conditioning on
$\ray$, $X_0$ becomes independent from $X_1$.  Therefore,
$$
\E[X_0X_1]=\E\bigl[\E[X_0X_1\mid \ray]\bigr]=
\E\Bigl[\E[X_0\mid \ray]\cdot\E[X_1\mid \ray]\Bigr]
\,.
\label e.inde
$$ 
We shall show that for small $\delta>0$,
$$
\E\Bigl[\bigl(W_\delta(r)-\E[X_j\mid \ray]\bigr)^2\Bigr]=O(1), \qquad j=0,1\,.
\label e.y
$$
Using \ref e.inde/, this implies
$$
\E[W_\delta(r)^2]-\E[X_0X_1]=O(1)\sqrt{\E[W_\delta(r)^2]}+O(1)
\,.
$$
Consequently, by \ref e.r/,
$$
\limsup_{\delta\to 0}\Bigl|\E[W_\delta(r)^2]-2\log(1/r)+O(1)\sqrt{\E[W_\delta(r)^2]}\Bigr|
\leq O(1)
\,,
$$
which implies \ref e.pwindw/. 
It therefore suffices to prove
\ref e.y/. 

Let $s:=\min\bigl\{|v-v_0|\st v\in \ray\bigr\}$.
Given $\ray$, let $T'$ be a UST of $\G^W(B(v_0,s/2),\delta)$,
and let $T''$ be the UST of the graph obtained from
$\Gd$ by identifying the vertices of $\ray$.  Note that
(by Wilson's algorithm, say) $T$, the UST on $\Gd$,
has the same law as $T''\cup\ray$ (as a set of edges).
By the domination principle, given $\ray$, we may couple $T''$
and $T'$ so that $\edges(T'')\supset \edges(T')$.
Given $\ray$, we couple $T$ and $T'$ so that $T\supset T'$.

Let $\alpha$ be the path in $T'$ from $v_0$ to 
$\bd B(v_0,s)$, and let $a$ be the point where $\alpha$
hits $\bd B(v_0,s)$.
Then $\E[\winding(\alpha,v_0)\mid\ray]=0$, by symmetry, because
given $\ray$, $T'$ is just ordinary UST on $\G^W(B(v_0,s/2),\delta)$.
But $\alpha$ is also a path in $T$.
Let $\beta$ be the path in $T$ from $a$ to the endpoint of $\gamma_r^\delta$
near $\bd B(v_0,r)$.
Then $X_0=\winding(\alpha,v_0)+\winding(\beta,v_0)+W_\delta(r)$,
and therefore,
$$
W_\delta(r)-\E[X_0\mid \ray]
=
-\E[\winding(\alpha,v_0)\mid\ray] -\E[\winding(\beta,v_0)\mid\ray] 
=
-\E[\winding(\beta,v_0)\mid\ray]
\,.
\label e.break
$$

Let $t>\delta$.
Note that if $s<t$, then there are at least two
disjoint crossings in $T$ of the annulus $A(v_0,t,r/10)$.
Therefore, \ref l.expo/ gives
$$
\P[s<t]\leq O(1) (t/r)^{\nco d}
\,.
\label e.ss
$$ 
Fix some $y\geq 2$.  By \ref l.windtail/, we have
$$
\P\Bigl[\bigl|\winding(\beta,v_0)\bigr|>t,\ s\geq r/y\Bigr]
\leq
O(1) \exp\bigl(-\co b t/\log y\bigr)\log y
\,.
$$
Hence, using \ref e.ss/,
\begineqalno
\P\Bigl[\bigl|\winding(\beta,v_0)\bigr|>t\Bigr]
&
\leq 
\P[ s< r/y\Bigr]+
O(1) \exp\bigl(-\co b t/\log y\bigr)\log y
\cr &
\leq
O(1) y^{-\co d} +
O(1) \exp\bigl(-\co b t/\log y\bigr)\log y
\,.
\endeqalno
Assuming that $t\geq1$, we may choose $y=\exp\sqrt t$,
and then get
$$
\P\Bigl[\bigl|\winding(\beta,v_0)\bigr|>t\Bigr]
\leq
O(1)\sqrt t \exp\bigl(-\nco z \sqrt t\bigr)
\,,
$$
for some constant $\co z>0$.
This gives $\E\bigl[\winding(\beta,v_0)^2\bigr]=O(1)$.
But for every random variable $Y$, we have
$\E[Y^2]\geq \E\bigl[\E[Y\mid\ray]^2\bigr]$.
Therefore, \ref e.break/ implies \ref e.y/ for $j=0$.
The proof of~\ref e.y/ for $j=1$ is entirely the same.
This completes the proof of the proposition.
\Qed

\procl p.theconst
Assuming \ref g.confinv/, $\lecl=\xplec$, where
$\lecl$ is the constant 
such that SLE with parameter $\lecl$ is the scaling limit of LERW.
\endprocl

\proof
Recall the definition of $W_\delta(r)$, which appears
above \ref p.windw/. 
Set $r_0=1/2$, and let $r_1>0$ be very small.
Let $\psi:\H\to\U$ be the conformal map satisfying
$\psi(i)=0$ and $\psi'(i)>0$.
Let $C_j$ be the circle of radius $r_j$ about $i$, $j=0,1$,
and set $C_j':=\psi(C_j)$.
Note that $Z_\delta(r_1)=W_\delta(r_1)-W_\delta(r_0)$ is the
winding number around $i$ of some arc on $\gamma^\delta$,
the LERW from a vertex near $i$ to $\partial\H$ in
$\delta\Z^2\cap\H$, and the arc has one endpoint near the
circle $\partial B(i,r_0)$ and the other endpoint
near the circle $\partial B(i,r_1)$.
It follows that $Z_\delta(r_1)$ converges weakly to a winding number
$Z(r_1)$ of an arc $\beta$ of $\psi^{-1}(\lep_\lecl)$
with endpoints on $C_0$ and $C_1$,
as $\delta\to0$ along some sequence, where
$\lep_\lecl$ is the SLE curve with parameter $\lecl$.
Moreover, since we have good tail estimates on $Z_\delta(r_1)$
(\ref l.windtail/), from the dominated convergence theorem
it follows that
$$
\E[Z(r_1)^2]=\lim_{\delta\to 0} \E[Z_\delta(r_1)^2]
\,.
$$
Hence, \ref p.windw/ gives,
$
\E\bigl[Z(r_1)^2\bigr] = 2\log(1/r_1) + O(1) \sqrt{\log(1/r_1)}
$.
Observe that for any path $\alpha$ in $\H-\{i\}$, the
winding number of $\alpha$ around $i$ minus the winding
number of $\psi(\alpha)$ around $0$ is bounded by some constant.
Consequently, the winding number $W'$ of $\beta':=\psi(\beta)$
around $0$ also satisfies
$$
\E\bigl[W'^2\bigr] = 2\log(1/r_1) + O(1) \sqrt{\log(1/r_1)}
\,.
\label e.www
$$ 
Set $t_j=\log r_j$, $j=0,1$, let $\tilde \beta$
be the arc $\lep_\lecl(t):[t_1,t_0]\to\U$,
and let $\tilde W$ be the winding number of $\tilde \beta$
around $0$.
By \ref t.wind/, with high probability,
the log of the absolute value of the
endpoints of $\tilde \beta$ is not far from the log of the
absolute value of the endpoints of $\beta'$.
Therefore, it is easy to conclude with the help of  \ref l.windtail/, 
that 
$$
\Ebig{(\tilde W- W')^2}=O(1)
\,.
\label e.wwp
$$ 
We know from \ref t.wind/ again that
$$
\Ebig{{\tilde W}^2}=\lecl|t_1| + O(1)\sqrt{|t_1|}\,.
$$
Combining this with \ref e.www/ and \ref e.wwp/ gives
$$
 (2-\lecl)\left|t_1\right| = O(1)\sqrt{\left|t_1\right|}\,.
$$ 
Letting $t_1\to-\infty$ now completes the proof.
\Qed 

\proofof t.lesl
Immediate from \ref c.brows/ and \ref p.theconst/. \Qed

\bsection{The critical value for the SLE}{s.crit}
\resetConsts

\procl t.kk
$\sup\KK\leq \crit$, where
$\KK$ is as in \ref d.lep/.
\endprocl

\proof
Fix some $\lecv\in\KK$, and let $f_t$ be the solution of the 
\lowner{} equation with parameter $\zeta(t)=\B(-\lecv t)$,
where $\B:[0,\infty)\to\bd\U$ is Brownian motion starting from
a uniform point in $\bd\U$. 
Note the for every $t<0$ the map
$f_t^{-1}$ is well defined and injective on $\partial \U-\{\zeta(t)\}$,
since $f_t$ is a Riemann map onto a slit domain, and the
slit hits $\partial\U$ at $\zeta(t)$.
Set $b(t)=-i\log f_t^{-1}(1)$, with $b(0)-\hat\B(0)\in[0,2\pi)$ and
$b(t)$ continuous in $t$.  Then $b(t)$ is real. As in \ref e.ab/,
we have
$$
b'(t) =
{\sin\bigl(\hat\B(-\lecv t)-b(t)\bigr)\over 1-\cos\bigl(\hat\B(-\lecv t)-b(t)\bigr)} 
= \cot{1\over 2}\bigl(\hat\B(-\lecv t)-b(t)\bigr)
\,.
\label e.deq
$$

Let $p(s)=b(-s)-\hat\B(\lecv s)$, and let $\eps>0$.
Set $\tau_\eps=\inf\{s\geq 0\st p(s)=\eps\}$ and
$\tau_\pi =\inf\{s\geq 0\st p(s)=\pi\}$.
For $x\in[\eps,\pi]$, let $g_\eps(x)$ be the probability that
$\tau_\pi<\tau_\eps$, conditioned on $p(0)=x$.  Also set 
$g_\eps(x)=0$ for $x<\eps$ and $g_\eps(x)=1$ for $x>\pi$.
We now show that $g_\eps$ satisfies 
$$
{\lecv\over 2} g_\eps''(x) + g_\eps'(x) \cot (x/2)=0
\label e.gde
$$
inside $(\eps,\pi)$, using It\^o's formula.
(The reader unfamiliar with stochastic calculus can 
have a look at \ref b.Durrett:BMM/, for example,
or try to derive \ref e.gde/ directly.  The latter is a bit
tricky, but can be done.)
Observe that $g_\eps\bigl(p(s^*)\bigr)$ is a martingale, where
$s^*=\min\{s,\tau_\eps,\tau_\pi\}$. 
By \ref e.deq/, we have
$$
dp(s) = -b'(-s) ds -d\hat\B(\lecv s)
= \cot\bigl(p(s)/2\bigr) ds - d\hat\B(\lecv s)
\,,
$$
and therefore, by It\^o's Formula
(assuming, for the moment, that $g_\eps$ is $C^2$), 
\begineqalno
d g_\eps\bigl(p(s)\bigr) 
&= g_\eps'\bigl(p(s)\bigr)\Bigl(\cot\bigl(p(s)/2\bigr)\,ds-d\hat\B(\lecv s)\Bigr)
+(1/2)g_\eps''\bigl(p(s)\bigr)\,d\left<\hat\B(\lecv s)\right>
\cr &
=\Bigl( \cot\bigl(p(s)/2\bigr)g_\eps'\bigl(p(s)\bigr)
+(\lecv/2)g_\eps''\bigl(p(s)\bigr)\Bigr)
\,ds
 -g_\eps'\bigl(p(s)\bigr)\,d\hat\B(\lecv s)
\endeqalno
for $s<\min\{\tau_\eps,\tau_\pi\}$. 
Since $g_\eps\bigl(p(s^*)\bigr)$ is a martingale,
the $ds$ term must vanish, and so \ref e.gde/ holds 
inside $(\eps,\pi)$.  Consequently,  in that range,
$$
g_\eps'(x) = c_\eps \bigl(\sin (x/2)\bigr)^{-4/\lecv}
\,,
$$ 
where $c_\eps$ is some constant depending on $\eps$. 
Since $g_\eps(\eps)=0$ and
$g_\eps(\pi)=1$, we have $\int_{\eps}^\pi g_\eps'(x)\,dx=1$, which gives
$$
c_\eps^{-1} = \int_{\eps}^{\pi}
\bigl(\sin (x/2)\bigr)^{-4/\lecv}\,dx
\,.
$$
We know that a.s.\ $p(s)\neq 0$ for all $s$,
which is equivalent to $\lim_{\eps\to 0} g_\eps(x)= 1$ on $(0,\pi)$.
This gives $\lim_{\eps\to0} g_\eps'(x)= 0$; that is, $\lim_{\eps\to 0}c_\eps=0$.
Therefore, $\lecv\leq 4$.

This completes the proof, except that we have not shown
that $g_\eps$ is $C^2$ (there should be a reference implying
this, but we have not located one).
To deal with this, the above procedure is reversed.
{\bf Define} $g_\eps$ as the solution of \ref e.gde/
satisfying $g_\eps(\eps)=0$ and $g_\eps(1)=1$.
Then the above application of  It\^o's Formula
shows that $g_\eps\bigl(p(s^*)\bigr)$
is a martingale.
By the Optional Sampling Theorem,
this implies that $g_\eps(x)$ is the probability that
$\tau_\pi<\tau_\eps$, conditioned on $p(0)=x$, and completes the proof.
\Qed 

\procl g.kk  
$\KK=[0,\crit]$.
\endprocl

\bsection{Properties of UST subsequential scaling limits in two dimensions}{s.ustsl}
\resetConsts

Before we go into the study of the UST scaling limit, let us remark that the
definition we have adopted for the scaling limit is by no means 
the only reasonable one.  There are several other reasonable variations,
and choosing one is partly a matter of convenience and taste.

We now recall some definitions.
Again, we think of
$\dZ$ as a subset of the sphere $\Sp2=\R^2\cup\{\infty\}$.
Recall that $\hT_\delta$ denotes the UST on $\dZ$, with the
point $\infty$ added, to make it compact.
Given two points $a,b\in\hT_\delta$,  $a\neq b$,
$\path_{a,b}=\path_{a,b}^\delta$
denotes the unique path in $\hT_\delta$ with endpoints
$a$ and $b$.
For the case $a=b$, we set $\path_{a,a}=\{a\}$.
Let $\ST_\delta$ be the collection of all triplets,
$(a,b,\path_{a,b})$, where $a,b\in\hT_\delta$.
$\ST_\delta$ will be called the {\bf paths ensemble} of $\hT_\delta$.
Let $\ST$ denote a random variable in
$\Ha\bigl(\Sp2\times\Sp2\times\Ha(\Sp2)\bigr)$
whose law is a
weak subsequential limit of the law of $\ST_\delta$
as $\delta\to0$.
The trunk is defined by
$$
\Tr=\Tr(\ST):=\bigcup_{(a,b,\path)\in\ST} ( \path-\{a,b\})
\,.
\label e.trunckdef
$$

Let $\eps\in(0,1]$, and $a,b\in\hT_\delta$. 
We define $\path_{a,b}(\eps)$ as follows.
Let $a'$ be the first point along the path $\path_{a,b}$ (which
is oriented from $a$ to $b$)
where $d(a,a')=\eps$, and let $b'$ be the last
point along the path where $d(b,b')=\eps$, provided
that such points exist.
If $a'$ and $b'$ exist, and $a'$ appears on the path
before $b'$, then let $\path_{a,b}(\eps)$ be the (closed)
subarc of $\path_{a,b}$ from $a'$ to $b'$; and otherwise
set $\path_{a,b}(\eps)=\emptyset$.
Let $\ST_\delta(\eps)$ denote the set of all triplets
$\bigl(a,b,\path_{a,b}(\eps)\bigr)$  such that
$a,b\in\hT_\delta$ and $\path_{a,b}(\eps)\neq\emptyset$.
Note that if $d(a,b)>2\eps$, then $\path_{a,b}(\eps)\neq\emptyset$.
We define
$$
\Tr_\delta(\eps)
:=\bigcup\Bigl\{\path\st (a,b,\path)\in\ST_\delta(\eps)\Bigr\}\,.
$$
Then $\Tr_\delta(\eps)$ is a compact subset of $\hT_\delta$,
which we call the {\bf $\eps$-trunk} of $\hT_\delta$.  
By compactness, for every $\eps\in(0,1]$ there is a subsequential scaling
limit of the law of $\Tr_\delta(\eps)$.
By passing to a subsequence, if necessary, we assume that
for all $n\in\Np:=\{1,2,\dots\}$ the weak limit
$\Tr_0(1/n)$ of $\Tr_\delta(1/n)$, as $\delta\to 0$, exists.

Recall that the {\bf dual} $\d T$ of a spanning tree $T\subset \dZ$ is the
spanning subgraph of the
dual graph $\d{(\dZ)}:=(\delta/2,\delta/2)+\dZ$
containing all edges that do not
intersect edges in $T$.  If $T$ is the UST on $\dZ$, then $\d T$
has the law of the UST on $\d{(\dZ)}$.  (See, e.g., \ref b.BLPS:usf/.)
Let $\d\hT_\delta$ be $\d T\cup\{\infty\}$, where $T$ is the
UST on $\dZ$.  $\d\Tr$ and $\d\Tr_0(\eps)$
are defined for $\d\hT_\delta$ as $\Tr$ and $\Tr_0(\eps)$
were defined for $\hT_\delta$.

We may think of the random variables $\Tr$,
$\d\Tr$, $\Tr_0(1/n)$, and $\d\Tr_0(1/n)$  ($n\in\Np$) 
as defined on the same probability space, by taking a subsequential
limit of the joint distribution of $\ST_\delta$, $\d \ST_\delta$,
$\bigl<\Tr_\delta(1/n)\st n\in\Np\bigr>$,
and $\bigl<\d\Tr_\delta(1/n)\st n\in\Np\bigr>$.
It is immediate to verify that a.s. 
$$
\Tr=\bigcup_{n\in\Np} \Tr_0(1/n)\,,
$$
and $\Tr_0\bigl(1/(n+1)\bigr)\supset \Tr_0(1/n)$ for $n\in\Np$.

We shall prove that $\Tr$ is a.s.\ a topological tree, in the sense
of the following definition.

\procl d.toptree \procname{Trees}
An {\bf arc} joining two points $x,y$ in a metric space
$X$ is a set $J\subset X$ 
such that there is a homeomorphism $\phi:[0,1]\to J$ with
$\phi(0)=x$ and $\phi(1)=y$.
A metric space $X$ will be called a {\bf topological tree}
if it is uniquely arcwise connected (that is,
given $x\neq y$ in $X$ there is a unique arc in $X$ joining $x$ and $y$)
and locally arcwise connected (that is, whenever $x\in U$ and $U$ is
an open subset of $X$ there is an open $W\subset U$ with $x\in W$ and
$W$ is arcwise connected).
A {\bf finite topological tree} is a topological space which is homeomorphic
to a finite, connected, simply connected, $1$-dimensional simplicial complex.
\endprocl

Note that a connected subset of a topological tree is
a topological tree \ref b.Bowditch:treelike/.

Although we shall not need this fact, it is instructive to note
that a metric space  which is a topological tree
is homeomorphic to an $\R$-tree%
\Ftnote{An $\R$-tree is a
metric space $(T,d)$ such that for every two
distinct points $x,y\in T$ there is a unique isometry $\phi$
from $[0,d(x,y)]$ onto a subset of $T$ satisfying
$\phi(0)=x$ and $\phi\bigl(d(x,y)\bigr)=y$.}
\ref b.MO:treechar/ 
(see also \ref b.MMOT:treechar/, for
a slightly less general but simpler proof). 

The next theorem establishes a finiteness property of the $\eps$-trunks,
which is the first step in the proof of \ref t.ustlim/.

\procl t.fin \procname{Finiteness}
For every $\eps>0$ there is a $\hat\delta>0$ with the following property.
Suppose that $0<\delta<\hat\delta$.
Let $\hat V$ be a set of vertices of $\dZ$ such that
every point in $\Sp2$ is within distance $\hat\delta$ of some vertex
in $\hat V$.
Let $\QQ:=\QQ_\delta(\hat V)$
be the subtree of $\hT_\delta$ that is spanned by
$\hat V$; that is, the minimal connected subset of $\hT_\delta$
containing $\hat V$. 
Then with probability at least $1-\eps$ we
have $\QQ\supset\Tr_\delta(\eps)$.
\endprocl

\proof
Fix some small $\hat \delta>0$, and suppose that $\delta\in(0,\hat \delta)$.
Let $V_0:=\hat V$, and for each $j\in\Np$
let $V_j$ be a set of vertices containing $V_{j-1}$
such that every vertex of $\dZ$ is within spherical distance
$\delta_j:=2^{-j}\hat\delta$ of some vertex in $V_j$, and $V_j$ is a
minimal set satisfying these properties.
Note that the number of vertices in $V_j-V_{j-1}$
is bounded by $O(1)\delta_j^{-2}$.
Let $\QQ_j$ be the subtree of $\hT_\delta$ spanned by $V_j$.

We now estimate the probability that there is some component
of $\QQ_{j+1}-\QQ_{j}$ whose diameter is large.
Let $v$ be some vertex in $\dZ$, 
let $\QQ(v,j)$ be the arc of $\hT_\delta$ that connects $v$ to $\QQ_j$,
and let $\ev D(v,j,a)$ be the event the diameter of $\QQ(v,j)$
is at least $a\delta_j$.
By Wilson's algorithm, we may obtain $\QQ(v,j)$ 
by conditioning on $\QQ_j$ and loop-erasing a simple random
walk from $v$ that stops when $\QQ_j$ is hit.
Every vertex $w\in \delta\Z^2$ is within distance
$\delta_j$ from a vertex in $\QQ_j$.  Since
$\QQ_j$ is connected and has diameter at least $1$,
\ref l.harmm/ shows that
there is a universal constant $\nco 1>0$
so that the probability that a random walk from $w$ gets to distance
$\co 1\delta_j$ from $w$ before hitting $\QQ_j$ is at most $1/2$. 
Consequently, $\ev D(v,j,a)$ has probability at most
$O(1)\exp({-\nco 3 a})$, where $\co 3>0$ is an absolute constant.
We choose $a_j:=j^2 (\log\hat \delta)^2/\co 3$.
Since there are at most $O(1)\delta_j^{-2}$ vertices in
$V_{j+1}-V_j$, we find that the probability of
$$
\ev D := \bigcup_{j=1}^{\infty}\bigcup_{v\in V_{j+1}}\ev D(v,j,a_j)
$$
is bounded by
$$
O(1)\hat\delta^{-2}\sum_{j=1}^\infty  2^{2j}\exp\bigl(-j^2(\log\hat\delta)^2\bigr)
\,,
$$
which goes to zero as $\hat\delta\to 0$.

Let $v\in \delta\Z^2$.
There is a sequence $v_1,v_2,\dots,v_{n}$ with
$v_j\in V_j$ such that $\QQ(v,1)\subset \bigcup_{j=2}^n \QQ(v_j,j-1)$,
and the latter union is connected.  If we are in the complement
of $\ev D$, it follows that the diameter of $\QQ(v,1)$ is
at most $s:=\sum_{j=1}^\infty a_j\delta_j$. 
Since $s\to 0$ as $\delta_0\to 0$, this establishes the theorem.
\Qed

Several corollaries follow from this theorem.

\procl c.fintree 
For each $n\in\Np$, a.s.\ $\Tr_0(1/n)$ is a finite topological tree.
\endprocl

\proof
Let $W\subset\R^2$ be finite.  For each $w\in W$ and $\delta>0$,
let $w_\delta\in\delta\Z^2$ be closest to $w$, with ties broken
arbitrarily, and set $W_\delta=\{w_\delta\st w\in W\}$.
Let $\QQ_\delta(W)$ be the subtree of $\hT_\delta$ spanned by $W_\delta$.
The theorem shows that we may choose a finite $W\subset\R^2$
such that
$\Tr_\delta(1/n)\subset \QQ_\delta(W)$ with probability at least $1-\eps$,
for every sufficiently small $\delta>0$.
Consequently, we may couple a subsequential scaling limit $\QQ(W)$
of $\QQ_\delta(W)$ as $\delta\to0$ so that 
$\Tr_0(1/n)\subset \QQ(W)$ with probability at least $1-\eps$.
Because $\Tr_0(1/n)$ is connected and $\eps$ is an arbitrary
positive number, it suffices to prove that $\QQ(W)$ is a.s.\ a finite
tree.  The latter is easily proved by induction on $|W|$ using
\ref t.noloop/, Wilson's algorithm, and the following easy fact:
the tree spanned by a subset of the points in $W$ is unlikely
to pass close by to the other points. (See \ref r.closetopt/.)
\Qed

\procl c.ustdim 
The Hausdorff dimension of $\Tr$ is in $(1,2)$.
Moreover,
if $I=(s_0,s_1]$ is an interval such that a.s.\ the Hausdorff dimension
of any scaling limit of LERW is in $I$, then the Hausdorff dimension
of $\Tr(\ST)$ is in $I$.
\endprocl

\proof
The second statement follows immediately from \ref t.fin/.
The first is now a consequence of the result of
\ref b.ABNW/, showing that there are $s_0,s_1\in(0,1)$ such
that a.s.\ the Hausdorff dimension of LERW scaling limit
is in $[s_0,s_1]$.%
\Ftnote{From \ref r.closetopt/ follows the weaker result that the
area measure of any subsequential scaling limit of LERW is zero,
hence that the area of $\Tr$ is zero. It is likely that with a bit
more effort the proof of
\ref r.closetopt/ is sufficient for the stronger claim that
the Hausdorff dimension is smaller than $2$.}
\Qed

\procl r.dimr 
The above-mentioned lower bound in \ref b.ABNW/ is based on the
ideas of \ref b.BJPP/.
Kenyon~\ref b.Kenyon:inprep/ can prove that we may take $s_1=5/4$.
In earlier work~\ref b.Kenyon:5o4/ he showed that 
$n^{5/4}$ times the expected number
of edges in a LERW from $(0,0)$ to the boundary of the
square $[-n,n]^2$ tends to a finite positive constant
as $n\to \infty$.
This supports the conjecture that the Hausdorff dimension of
the scaling limit of LERW is a.s.\ $5/4$,
and the same would apply to $\Tr(\ST)$.
\endprocl

The {\bf degree} of a point $p$ in a topological tree $T$ is
the number of connected components of $T-\{p\}$.
The following corollary is a strong form of the statement that the maximum
degree of points in $\Tr_0(1/n)$ is $3$.
{} From this
and the fact that $\Tr$ is a tree (which we prove further below)
it immediately
follows that the maximum degree in $\Tr$ is $3$,
because every finite subset of $\Tr$ is contained in some $\Tr_0(1/n)$.

Given a point $p\in\Sp2$ and two numbers $0<r_1<r_2<1$, let
$\Rsp(p,r_1,r_2)$ denote the annulus with center $p$, inner radius
$r_1$, and outer radius $r_2$, in the spherical metric.

\procl c.maxdeg
Given every $\eps\in(0,1)$, there is an $r\in(0,\eps)$ with the
following property.  For every sufficiently small $\delta>0$,
the probability that there is a point $p\in \Sp2$
such that there are $4$ disjoint crossings in $\Tr_\delta(\eps)$ of
the annulus $\Rsp(p,r,\eps)$ is at most $\eps$.
\endprocl

By having $4$ disjoint crossings in $\Tr_\delta(\eps)$ of an annulus $\An$,
we mean that there are $4$ disjoint connected subsets of $\Tr_\delta(\eps)$
that intersect both boundary components of $\An$.
Below, \ref c.stmaxdeg/ gives a strengthening of \ref c.maxdeg/.

\proof
By \ref t.fin/, it is enough to prove the statement with
$\QQ_\delta(W)$ replacing $\Tr_\delta(\eps)$, where $W\subset\R^2$ is 
a set of bounded size, provided that the value of $r$ does not depend
on $\delta$.  Again, induction on $|W|$ can be used together with
Wilson's algorithm.
One needs note the following easy facts. 
The tree $T_{k-1}$ spanned by $k-1$ points of $W$ is unlikely to pass close to
the other points of $W$, and when adding a further point, it
is unlikely that the attachment point of the new branch on $T_{k-1}$
will be close to another branch point.  Also, once a random
walk from the new point gets close to $T_{k-1}$ it will hit
$T_{k-1}$ close by, with high likelyhood.
The easy details are left to the reader.
\Qed

We now turn to the central issue in the proof of \ref t.ustlim/,
which is,

\procl t.noint
In any subsequential scaling limit of UST in $\Sp2$, a.s.\
the trunk and dual trunk do not intersect.
\endprocl

\procl l.3d2
Given $\delta,\eps>0$, let $T_\delta^3(\eps)$ be the set of points of
degree $3$ in $\Tr_\delta(\eps)$.  Let 
$\d{\Tr}_\delta(\eps)$ be the $\eps$-trunk of the dual tree
$\d{\hT_\delta}$.  Let $D$ be the spherical distance from
$T_\delta^3(\eps)$ to $\d{\Tr}_\delta(\eps)$; that is, the
least spherical distance between a point in $T_\delta^3(\eps)$ to
a point in $\d{\Tr}_\delta(\eps)$.  Then
$\lim_{t\to0}\P[D<t]\to 0$ uniformly in $\delta$.
\endprocl

The following simple observation is used in the proof.
Suppose that we condition on a set of edges $S$ to appear in
the UST tree in a planar graph.  For the dual tree,
this is the same as deleting the edges dual to the edges in $S$.
Consequently, one can perform a variation on Wilson's algorithm
for a planar graph, where one switches back and forth from building
the tree by adding LERW branches and building the dual tree.
When building the tree, the LERW acts with the constructed tree
as a wired absorbing boundary and the constructed dual tree as a free boundary,
and conversely when building the dual tree.

\proof We first choose a large but finite collection of
points $\QQ$ in $\delta\Z^2$ so that with high probability the
subtree $T$ of $\hT_\delta$ spanned by $\QQ$ contains $\Tr_\delta(\eps)$
(and $|\QQ|$ does not depend on $\delta$).
This can be done, by \ref t.fin/.  
Let $\d \QQ$ be a set of vertices of the dual graph $\d{(\delta\Z^2)}$,
such that with high probability the subtree $\d T$ spanned
by $\d \QQ$ in the dual graph contains $\d{\Tr}_\delta(\eps)$.
Let $a_1,a_2,a_3\in \QQ$ and $\d a_1,\d a_2\in\d \QQ$ be distinct points.
It suffices to show that the probability that the arc $\beta$
joining $\d a_1$ and $\d a_2$ in $\d T$ comes within distance
$t$ of the meeting point $m$ of $a_1,a_2$ and $a_3$ in $T$
goes to zero as $t\to 0$, uniformly in $\delta$.
This is easy.  We condition on the subtree $T_0$ of $T$ spanned
by $a_1,a_2,a_3$.  Let $\d z$ be a dual vertex close to $\d a_1$.
Then with high probability the dual tree path $\beta$ from $\d z$ to $\d a_1$
has diameter not much larger than the distance from $\d z$ to $\d a_1$.
In particular, it does not go close to $T_0$.  By the next lemma, conditioned
on $\beta$ and $T_0$, the probability that a simple random walk
starting at $\d a_2$, with $T_0$ acting as a reflecting boundary,
will get to within distance $t$ of $m$ before hitting $\beta$ is as small
as we wish.  Consequently, the same is true for the loop-erasure
of this walk, which can be taken as the path joining $\beta$
and $\d a_2$ in the dual tree.
\Qed

\procl l.harm
Let $D$ be a domain in $\Sp2$ with two boundary components,
$B_1,B_2$, and assume that both are not single points.
Consider a sequence $\delta_j$, $j\in\N$, of positive numbers tending to
zero.
Suppose that to each $j\in\N$ there are two connected subgraphs
$B_1^j,B_2^j$ of the grid $\delta_j\Z^2$, and that $B_1^j\to B_1$
and $B_2^j\to B_2$ in the Hausdorff metric on compact subsets of $\Sp2$. 
Let $m\in B_2$ be some point, and for each $t>0$ and $z\in\delta_j\Z^2-B_2^j$,
let $h_j(z,t)$ be the probability that simple random walk on
$\delta_j\Z^2-B_2^j$ starting at $z$
(with reflecting boundary conditions on $B_2^j$),
will get to within distance $t$ of $m$ before hitting $B_1^j$.
Let $K\subset\Sp2-B_2$ be compact.  Then
$$
\lim_{t\to 0\atop j\to\infty}
\sup\bigl\{h_j(z,t)\st z\in K\cap \delta_j\Z^2\bigr\}= 0
\,.
$$
\endprocl

\proof
Set $G_j:=\delta_j\Z^2-B_2^j$, 
let $S(t)$ be the vertices of $G_j$ that are within distance $t$ from $m$, and 
let $V_j(t)$ be the set of vertices of $G_j-B_1^j-S(t)$.
Then $h_j(z,t)$ is discrete-harmonic in $V_j(t)$. 
Recall that the Dirichlet energy of $h_j(z,t)$ is 
$\sum \bigl(h_j(z,t)-h_j(z',t)\bigr)^2$,
with the sum extending over all edges $[z,z']$ in $G_j$.
Let $\lambda=\lambda_j$ be the minimum of $h_j(z,t)$ on $K$,
and let $z$ be where the minimum is achieved.
Then there is a path $\beta$ from $z$ to $S(t)$
such that $h_j(z,t)\geq\lambda$ on $\beta$, by the maximum principle for
discrete harmonic functions.
Note that one can find a collection of $1/O(\delta_j)$
disjoint paths in $G_j$ which join $B_1^j$ and $\beta$
and each path in the collection has combinatorial length bounded by
$O(1)/\delta_j$.
The Dirichlet energy of $h_j(z,t)$ restricted to each such
path is at least $\lambda/O(\delta_j)$, and therefore
the Dirichlet energy of $h_j(z,t)$ is
at least $C\lambda$, where $C>0$ is a constant  
depending only on $K$ and $D$.

Let $d_j$ be the distance from $m$ to $B_1^j$.
Since $h_j(z,t)$ is harmonic in $V_j(t)$, it minimizes the Dirichlet energy
among functions on $G_j$ that are $1$ on $S(t)$ and $0$
on $B_1^j$.  Therefore, the Dirichlet energy of $h_j(z,t)$
is at most the Dirichlet energy of the function
$f:G_j\to\R$, which is $1$ on $S(t)$, $0$ outside of $S(d_j)$,
and equal to $\log \bigl(d_j/|z-m|\bigr)/\log(d_j/t)$ elsewhere,
which is $O(1)/\log(d_j/t)$, as $j\to\infty$.  This
gives, $\lambda_j=O(1)/\log(d_j/t)$, and the lemma follows.
\Qed

\proofof t.noint
Before we go into the actual details, the overall plan of the proof
will be given (in a somewhat imprecise manner).
Let $t_0>0$.
It is not hard to reduce the theorem to the claim that with probability
close to $1$ the
path $\gamma\subset \hT_\delta$ which joins two fixed points $a_1,a_2$
does not have points $p$ close to it such that 
the path $\alpha_p\subset\hT_\delta$ joining $p$ to $\gamma$ 
is not contained in a small neighborhood of $\gamma$.
Let $Z$ be the set of points $p$ such that $\alpha_p$ does
not stay close to $\gamma$.  When we condition on $\gamma$,
The probability that $p\in Z$ goes to zero as $p$ tends to a point in
$\gamma$, by a simple harmonic measure estimate.  However, this is
not enough, since there are many different $p$'s close to $\gamma$.
We fix some collection $L_1$ of points close to $\gamma$, and
take a thick collection of points $L_2$ which are much closer to $\gamma$.
What we show is that conditioned on $\gamma$ and on $Z\cap L_2\neq\emptyset$,
the expectation of $N:=|Z\cap L_1|$ is much larger than
$\E[N\mid\gamma]$.  This is established by observing that
when $p\in L_1\cap Z$ is appropriately chosen,
the expected number of points $p'\in L_2$ such that
$\alpha_{p'}$ is contained in $\alpha_p$, except for a
small initial segment of $\alpha_{p'}$,
is quite large.  
It follows that
$$
\P[Z\cap L_2\neq\emptyset\mid\gamma]\leq
{\E[N\mid\gamma]\over \E[N\mid\gamma,\ Z\cap L_1\neq\emptyset]}
$$
is small, which suffices to prove the theorem.

We now give the details.
Fix four distinct points $a_1,a_2,b_1,b_2\in\R^2$.
Given $\delta>0$, let $a_1'$ and $a_2'$ be points
of $\delta\Z^2$ that are closest to $a_1$ and $a_2$,
respetively, and let $b_1'$ and $b'_2$ be 
vertices of the grid dual to $\delta\Z^2$ that are
closest to $b_1$ and $b_2$, respectively.
Let $\gamma$ be the path in $\hT_\delta$
that joins $a_1'$ and $a_2'$, and
given any $p\in\delta\Z^2$, let $\alpha_p$ denote the
path in $\hT_\delta$ from $p$ to $\gamma$.
Let $\beta$ be the path of  $\d{\hT_\delta}$ that joins $b_1'$ and $b_2'$.

Since $\Tr=\bigcup_n\Tr_0(1/n)$, and 
$\d\Tr= \bigcup_n\d\Tr_0(1/n)$,
\ref t.fin/ shows that it suffices to prove that 
the probability that the distance between $\gamma$ and
$\beta$ is less than $t$ goes to zero, as $t$ goes down
to zero, uniformly in $\delta$.
We know that with probability close to one,
$\beta$ does not come close to $\{a_1,a_2\}$,
and $\gamma$ does not come close to $\{b_1,b_2\}$ \ref r.closetopt/.
Therefore, we need only consider the situation where there
is a point $q$ on $\gamma$, which is close to $\beta$,
but not close to $\{a_1,a_2,b_1,b_2\}$. 
Since $\beta$ and $\gamma$ cannot cross, and since $\beta$ locally
separates the sphere near every point of $\beta-\{b_1',b_2'\}$,
such a situation implies that there is a point $q'$ in
$\delta\Z^2$, which is near $q$, but in order to get
to $\gamma$ from $q'$ one must either cross $\beta$,
or go ``around'' it.  Consequently, $\diam\left(\alpha_{q'}\right)$
must be bounded away from zero, as $\alpha_{q'}$ cannot
cross $\beta$.
It therefore suffices to rule out the existence of
a point $q'\in\delta\Z^2$ close to $\gamma$ but with
$\diam\left(\alpha_{q'}\right)$ bounded away from zero. 
More precisely, let $K$ be a compact set
disjoint from $\{a_1,a_2\}$, and let
$\eps_1\in(0,1)$.  Let $\tilde h$ be the least distance
from $\gamma$ to some point $q'\in K\cap \delta\Z^2$
such that $\diam\left (\alpha_{q'}\right)\geq\eps_1$.
It suffices to show that
$$
\inf_{h_0>0}\limsup_{\delta\to0}\P[\tilde h<h_0]=0
\,.
\label e.targg
$$
Given any $p\in\delta\Z^2$, let $h(p)$ be the distance
from $p$ to $\gamma$ and let $k(p)$ be the maximal distance
from a point on $\alpha_p$ to $\gamma$.
By \ref c.maxdeg/, the probability that there is an
arc $\alpha$ in $\hT_\delta$, which is disjoint from $\gamma$,
satisfies $\diam(\alpha)\geq\eps_1$, and every point of $\alpha$ is
within distance $t$ of $\gamma$, goes to zero as $t\to0$, uniformly
in $\delta$.
Hence, to prove \ref e.targg/, it suffices to establish that
$$
\forall t>0\qquad
\inf_{h_0>0}\limsup_{\delta\to 0}
\PBig{\exists p\in K\cap\delta\Z^2\,\,\,h(p)<h_0,\,k(p)\geq t}=0
\,.
\label e.targa
$$

Since the proof is somewhat involved, we consider
first the simpler situation in which
$$
\gamma=[0,1]\times \{0\},\hbox{ and }
K=[1/3,2/3]\times[0,1]
\label e.simp
$$
(notwithstanding that this is an unrealistic situation, of extremely low
probability).
Obviously, it suffices to prove \ref e.targa/ for small $t>0$.

In the following arguments, several small positive quantities appear. 
Their dependence differs from the natural flow of the proof.
In order to make it clear that the proof is logically
sound, we state now that the dependence order is
as follows: 
$$
\eps_0,t_0,r_0,t_1,h_1,r_1,h_2,\delta\,;
$$
that is, each of these quantities may depend only on those appearing
before it in the list, and should be thought of as much smaller
than it predecessors.

Set $h_1':=\max\{k\delta\st k\in\Z,\ k\delta\leq h_1\}$, and
let $L_1:=\bigl\{(k\delta,h_1')\st k\in\Z,\ 1/6<k\delta<5/6\bigr\}$.
Given $\gamma$ and $p\in\delta\Z^2$, we may choose $\alpha_p$ by
loop-erasing a simple random walk from $p$ to $\gamma$.
Consequently, an easy harmonic measure estimate shows that 
$\P[k(p)>t]=O(h_1/t)$ for all $p\in L_1$.

Set $h_2':=\max\{k\delta\st k\in\Z,\ k\delta\leq h_2\}$,
and $L_2:=\bigl\{(k\delta,h_2')\st k\in\Z,\ 1/4<k\delta<3/4\bigr\}$.
Again, for all $p\in L_2$, $\P[k(p)\geq t_0]=O(h_2/t_0)$,
so we may assume that the event $\ev Q$ that the leftmost point
in $L_2$ satisfies $k(p)< t_0$ has probability at least $1-\eps_0$.
Let $\ev K$ be the event that there is some $p\in L_2$ with $k(p)\geq t_0$,
and let $\ev K':=  \ev K\cap\ev Q$.
For proving \ref e.targa/ in the simpler situation \ref e.simp/,
it suffices to show that $\P[\ev K']=O(\eps_0)$
for all sufficiently small $\delta>0$.

Consider the following procedure for generating $\hT_\delta$ given
$\gamma$.  Perform Wilson's algorithm starting with the
vertices in $L_2$, in left-to-right order.  If we encounter in
this procedure some vertex $p\in L_2$ such that $k(p)\geq t_0$,
we stop, and let $p_0$ denote that vertex.  Let $T_0$ be the tree
constructed up to that point (including $\alpha_{p_0}$).
On the event $\ev K'$, let $p_1$ be the first point on $\alpha_{p_0}$
whose distance to $\gamma$ is at least $t_0$, and let
$\alpha$ be the arc of $\alpha_{p_0}$ from $p_0$ to $p_1$.
Let $\ev A_1$ be the event that $\alpha$ is not contained
in the rectangle $[1/5,4/5]\times [0,t_0]$.
Note that $\ev A_1$ implies that there is an arc in
$\alpha\cap ([1/5,4/5]\times [0,t_0])$ with diameter
at least $1/15$.  By considering this
arc and $\gamma$, \ref c.maxdeg/ shows that
$\P[\ev A_1\cap\ev K']<\eps_0$, assuming that
$t_0$ is sufficiently small.

On the event $\ev K'-\ev A_1$, let
$$
x_1:=\max\bigl\{x\in \R\st (x,h_1/2)\in\alpha\bigr\}
\,,
$$ 
let $U$ be the component of
$\R^2-\Bigl(\alpha\cup\bigl([x_1,\infty)\times\{h_1/2\}\bigr)\cup
\bigl(\R\times\{t_1\}\bigr)\Bigr)$ that contains $(x_1,\infty)\times\{h_1/2\}$,
and let $U'$ be the
set of point in $U$ that are within distance $t_1$ from $\alpha$.
See \ref f.uprime/.
Let $\ev A_2$ be the event that $\ev K'-\ev A_1$
occurs and $U'$ intersects $T_0$.  

\midinsert
\SetLabels
\B(.1*.19)$T_0$\\
\L\T(.26*.09)$p_0$\\
\R(.13*.5)$\alpha$\\
\E\L(1.01*.35)$h_1/2$\\
\E\L(1.01*.02)$0$\\
\E\L(1.01*.89)$t_1$\\
\T(.47*-.01)$x_1$\\
\T(.47*.7)$U'$\\
\endSetLabels
\centerline{%
\AffixLabels{%
\epsfysize=2.2in\epsfbox{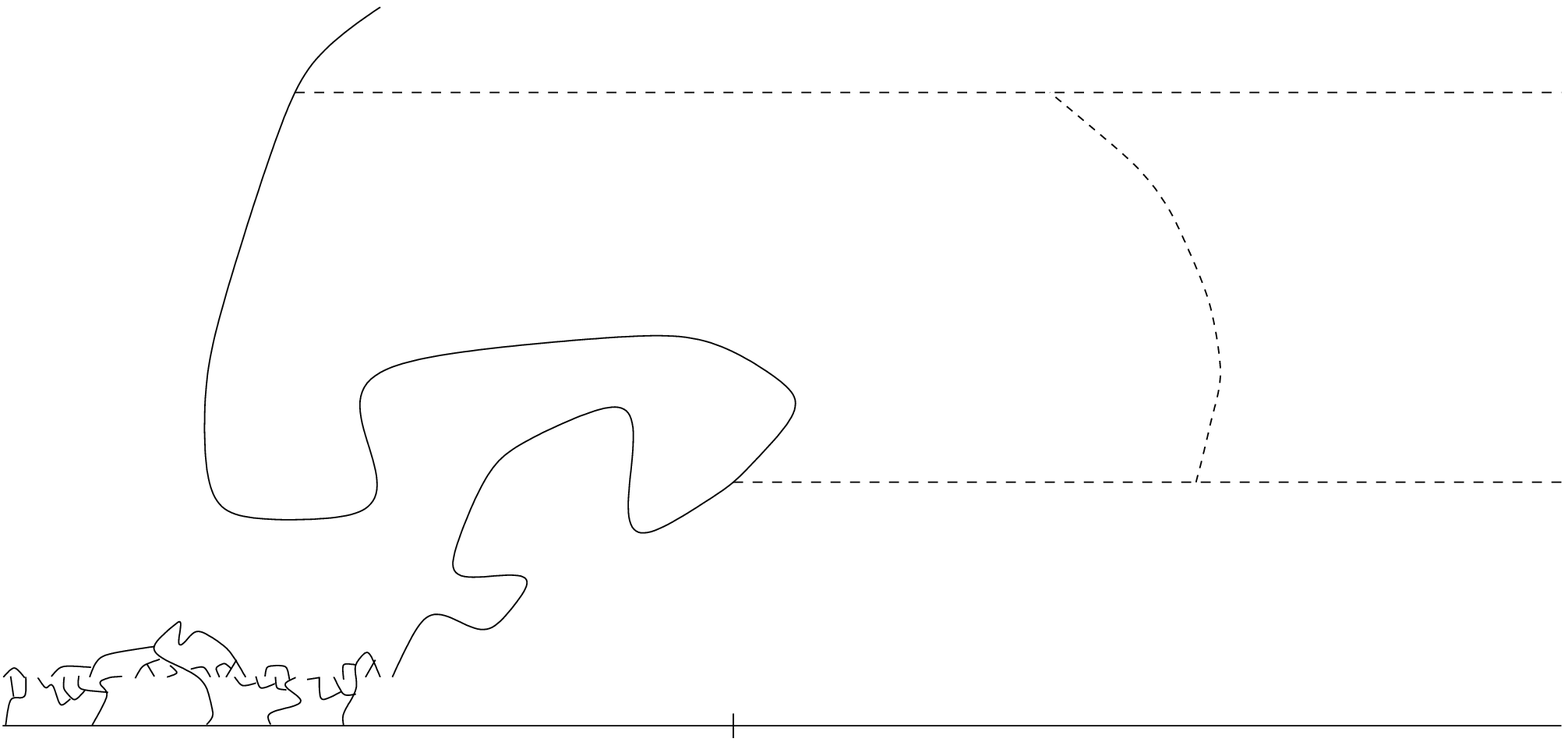}}}
\medskip
\caption{\figlabel{uprime}
}
\endinsert

We now prove that $\P[\ev K-\ev A_1-\ev A_2]=O(\eps_0)$.
Let $N$ be the number of points $p\in L_1$ such that
$k(p)\geq t_0$.  For a given $p\in L_1$, the probability
of $k(p)\geq t_0$ (given \ref e.simp/, but otherwise unconditioned)
is $O(h_1/t_0)$.  Therefore,
$$
\E[N]\leq O(1)h_1/(\delta t_0)
,.
\label e.enup
$$
On the other hand, condition on the event $\ev K'-\ev A_1-\ev A_2$
and on $T_0$.  Let $L_1'$ be the set
of $p\in L_1\cap U$ such that the distance from
$p$ to $\alpha$ is at most $t_1/2$. 
Note that conditioned on $p\in L_1'$, the probability that
$\alpha_p$ joins with $\alpha_{p_0}$ within distance
$2t_1$ from $p$ is at least 
$$
O(1)^{-1}{h_1\over \dist(p,\alpha)+h_1}
\,,
\label e.lowest
$$
since after generating $T_0$, we may continue by running Wilson's
algorithm starting at $p$, and the probability that
the random walk starting at $p$ will hit
$\alpha$ before the ray $(x_1,\infty)\times\{h_1/2\}$ is at least
\ref e.lowest/.
It therefore follows that conditioned on 
$\ev K-\ev A_1-\ev A_2$, we have 
\begineqalno
\E[N\mid \ev K-\ev A_1-\ev A_2]
& \geq
O(1)^{-1}\sum\bigl\{h_1/(k\delta)\st k\in\Z,\ h_1\leq k\delta
\leq t_1/2\bigr\}
\cr &
\geq 
O(1)^{-1}\delta^{-1}h_1 \log\bigl(t_1/(2 h_1)\bigr)
\,.
\endeqalno
Combining this with \ref e.enup/ gives
$$
\P[ \ev K-\ev A_1-\ev A_2] \leq 
{\E[N]\over \E[N\mid \ev K-\ev A_1-\ev A_2]}
\leq
 O(1)\Bigl( t_0 \log\bigl(t_1/(2 h_1)\bigr)\Bigr)^{-1}
\leq \eps_0
\,,
$$
provided that $h_1$ is sufficiently small.

It remains to establish that $\P[\ev A_2]\leq O(\eps_0)$.
First consider the case that there is some $p\in L_2$, to the
left of $p_0$, such that $\alpha_p\cap U\neq\emptyset$.
Then this must be the case for $p$ being the
left neighbor of $p_0$, namely $p=p_0':=p_0-(\delta,0)$,
because when $p\in L_2$ is to
the left of $p_0'$, the path $\alpha_p$ cannot cross
$\alpha_{p_0}\cup[p_0,p_0']\cup\alpha_{p_0'}$, and
$\alpha_p$ does not get to $\R\times\{t_0\}$.
If $\alpha_{p_0'}$ intersects $U$, then $\alpha_{p_0'}$
must first get to some point $z$ in $[x_1,\infty)\times\{h_1/2\}$.
Near $z$ there must be two points $\d z_1,\d z_2$, which are vertices
of the dual grid $\d{(\delta\Z^2)}$, and are locally
separated from each other by $\alpha_{p_0}$.
The path in the dual tree that joins  $\d z_1$ and $\d z_2$
has to contain the edge dual to the edge $[p_0',p_0]$.
Now consider another dual vertex $\d z_3$ just left
of $\alpha$ near the point $(x_1,h_1/2)$.
Let $\d m$ be the meeting point of
$\d z_1,\d z_2$ and $\d z_3$ in the dual tree.
If $\d m$ is not within distance $r_1$ of
$p_0$, we get in the $r_1/10$ trunk of the
dual tree at least four disjoint crossings of 
the annulus $\An(p_0,2h_2,r_1/2)$.
We may assume that this
has probability $\leq \eps_0$, by \ref c.maxdeg/.
Similarly, \ref l.3d2/, with the role  of the tree
and dual tree reversed, shows that we may take
the event that $\d m$ is within distance $r_1$ of
$p_0$ to have probability $\leq \eps_0$,
provided that $r_1$ is sufficiently small when compared with $h_1$.

To establish that $\P[\ev A_2]\leq O(\eps_0)$, it now
suffices to prove that
on the event $\ev K'-\ev A_1$, the probability that
$\alpha_{p_0}$ intersects $U'$ is $O(\eps_0)$.  The argument is
similar here, but occurs on a larger scale.
Suppose that $w\in\alpha_{p_0}\cap U'$.  Then $w$ must
be on the segment of $\alpha_{p_0}$ from $p_1$ to the 
point $p_2$ in $\alpha_{p_0}\cap\gamma$.
Because $w\in\alpha_{p_0}-\alpha$ and $w$ is within
distance $t_1$ to $\alpha$, there is a point near $w$
that is in the $(t_0/2)$-trunk of the dual tree; namely,
some point on the path connecting dual two vertices on opposite
sides of $\alpha_{p_0}$ near $p_1$.
Consequently, by \ref l.3d2/, we may rule out the possibility that
$p_2$ is within distance $r_0$ of $w$ as having small probability,
since $p_2$ is a point of degree $3$ of the $(t_0/2)$-trunk.
But if the distance between $p_2$ and $w$ is more than $r_0$,
then there are in the $(t_0/2)$-trunk at least four disjoint
crossings of the annulus $\An(w,2t_1,r_0)$:
two on $\alpha_{p_0}$ and two on $\gamma$.
An appeal to \ref c.maxdeg/ now establishes $\P[\ev A_2]\leq O(\eps_0)$.
This completes the proof in the situation \ref e.simp/.

We now explain how to modify the above proof to deal with the general
case.  First note that the restriction on $K$ is entirely
inconsequential; we could in the same way deal with
any compact set disjoint from the endpoints of $\gamma$.
More significant is the special selection of $\gamma$.
Observe that under the assumption of conformal invariance
of the LERW scaling limit, the general case can be reduced to
the case where $\gamma=[0,1]\times \{0\}$, because after
$\gamma$ is generated, the rest of the UST is just
unconditioned UST on the complement of $\gamma$
with wired boundary conditions.  We may then transform
$\gamma$ by a conformal homeomorphism to $[0,1]\times\{0\}$,
and refer to the above result.

Although we do not assume conformal invariance of the scaling
limit, it turns out that the proof above is itself conformally
invariant.  With some care, one can apply the conformal
map to the proof, in a manner of speaking.  This is actually
not very surprising, because the proof is ultimately
based on a simple (discrete) harmonic measure estimate, which is
conformally invariant.

Let us turn to the details.
We may couple the UST for a subsequence of $\delta$
tending to zero so that $\gamma$ tends to some path
$\gamma_0$ as $\delta\to 0$ along that subsequence (see the
discussion of the Prohorov metric in \ref s.back/).
Let $f_\delta:\Sp2-([0,1]\times\{0\})\to\Sp2-\gamma$
be the conformal map
normalized to take the endpoints of $[0,1]\times\{0\}$ to
the endpoints of $\gamma$ and so that $f_\delta(\infty)$ is on the
line which is the set of points at equal distance from both
endpoints of $\gamma$, say.
(The latter normalization is necessary to make $f_\delta$ unique,
but otherwise, it is quite arbitrary.)
It follows that $f_\delta$ tends to the similarly normalized
conformal map $f:\Sp2-([0,1]\times\{0\})\to\Sp2-\gamma_0$.
We may assume that $\delta$ is so small that $f$ and
$f_\delta$ are very close on compact subsets disjoint from
$[0,1]\times\{0\}$. 

For each $p\in L_1$, where $L_1$ is as before,
we let $\hat p$ denote a point
in $\delta\Z^2$ that is closest to $f(p)$.
For the general case, we consider $\hat N$,
the number of $p\in L_1$ such that $k(\hat p)$ is not small,
in place of $N$.
Let $\hat L_2$ denote the set of points in $\delta\Z^2$
that are within distance $3\delta$ of
$f\bigl([1/4,3/4]\times \{h_2\}\bigr)$.
The proof for the general case uses 
$\hat L_2$ in place of $L_2$.
For traversing $\hat L_1$,
there is no clear notion of the left-right order.
But any ordering that starts near $f\bigl((1/4,h_2)\bigr)$
and later does not visit any vertex before visiting an immediate
neighbor, will do.  Instead of the left neighbor $p_0'$ of
a vertex $p_0\in L_1$, we use for $\hat p_0\in \hat L_1$
that neighbor of $\hat p_0$ in $\hat L_1$ that is ``most counterclockwise'',
in the appropriate sense.
The rest of the proof proceeds with essentially no modifications,
except that the coordinate system used is transformed by $f$.
\Qed

\procl r.usfconj
In \ref b.BLPS:usf/ it has been asked whether the free
USF on every planar proper bounded degree graph is a tree.
The proof of \ref t.noint/ seems to be relevant.  It is plausible
that with a similar argument one can prove that for proper planar graphs
with bounded degree and a bounded number of sides per face, the
free USF is a tree.
\endprocl

We may now strengthen \ref c.maxdeg/, as follows

\procl c.stmaxdeg 
Given every $\eps\in(0,1)$, there is an $r\in(0,\eps)$ with the
following property.  For every sufficiently small $\delta>0$,
the probability that there is a point $p\in \Sp2$
such that there are $4$ disjoint crossings in $\hT_\delta$ of
the annulus $\Rsp(p,r,\eps)$ is at most $\eps$.
\endprocl

\proof 
Consider an annulus $\An=\Rsp(p,r,\eps)$, and suppose that there
are four disjoint paths $\alpha_0,\alpha_1,\alpha_2,\alpha_3$
in $\hT_\delta$ that cross it.  
Let $B_1$ be the component of $\Sp2-\An$ inside the inner boundary
component of $\An$, and let $B_2$ be the outside component.
Without loss of generality, we suppose that $\alpha_0\cup\alpha_2$
separate $\alpha_1$ from $\alpha_3$ inside $\An$; that is, the
circular order of these paths around $\An$ agrees with the order
of the indices.

Assume first that there are no paths in $\hT_\delta\cap \An$
that join two of the paths $\alpha_j$, $j=0,1,2,3$.
Then there must be paths in the dual tree
$\beta_0,\beta_1,\beta_2,\beta_3$, such that $\beta_j$ is
between $\alpha_{j-1}$ and $\alpha_j$ (indices mod $4$),
for each $j=0,1,2,3$.
If $\alpha_0,\alpha_1,\alpha_2$ and $\alpha_3$ can
all be connected to each other by paths in $\An\cup B_1$,
it follows that there are four crossings of the
annulus $\Rsp\left(p,r,\eps/3\right)$
in the $\eps/3$ trunk of $\hT_\delta$, and we know that has
small probability to happen anywhere, if $r$ is small,
by \ref c.maxdeg/.  If neither of the paths $\alpha_j$
connects to another in $\An\cup B_1$, the same argument applies
to the dual tree, because the paths $\beta_j$ must
all connect inside $\An\cup B_1$.
However, if two of the paths $\alpha_j$ connect in
$\An\cup B_1$, and one of the others does not connect to them,
then also two of the paths $\beta_j$ connect.  This implies
that the $\eps/3$ trunk gets within distance of $r$ from
the dual trunk, and again this can be discarded as having
small likelyhood.

We are left to deal with the situation where there is
a simple path $\gamma$ in $\An\cap\hT_\delta$ that  connects
two of the paths $\alpha_j$.  Note that for each pair
of paths $\alpha_j$ there can be at most one such $\gamma$
connecting them.  Also note that any path connecting
$\alpha_j$ and $\alpha_{j+2}$ (indices mod $4$)
must cross either $\alpha_{j+1}$ or $\alpha_{j+3}$.
Consequently,  if we consider any four concentric annuli
$\An_j:=\Rsp(p,r_j,r_{j+1})$, $j=0,1,2,3,4$, with
$r_j<r_{j+1}$ for each $j=0,1,2,3,4$, $r_0=r$, and $r_5=\eps$, at least
one of them will have the property that inside it there
is no path joining any two paths among the $\alpha_j$'s.
This allows a  reduction to the previous case, and completes the
proof.
\Qed

\procl t.trpaths
A.s., every simple path $\phi:[0,1)\to\Tr$
has a limit $\lim_{t\to1}\phi(t)$ in $\Sp2$,
and for every point $z\in\Sp2$ there is a simple path
$\phi:[0,1)\to\Tr$ such that $\lim_{t\to1}\phi(t)=z$.
\endprocl

\proof
Suppose that there are two distinct accumulation points,
$x$ and $y$, of $\phi(t)$ as $t\to 1$, and let
$m\in\Np$ satisfy $\dsp(x,y)>9/m$.
Then for each $t\in (0,1)$ the (spherical) diameter of
$\phi\bigl([t,1)\bigr)$ is 
greater than $9/m$.  Let $a\in(0,1)$ be such that the diameter of
$\phi\bigl([0,a])$ is at least $5/m$.  It easily follows that
$\phi\bigl([a,1)\bigr)\subset\Tr_0(1/m)$.
But since $\Tr_0(1/m)$ is a compact finite tree (\ref c.fintree/),
and the restriction of $\phi$ to $[a,1)$ is in $\Tr_0(1/m)$,
it follows that $\lim_{t\to 1}\phi(t)$ exists.  Contradiction.

Let $z\in\Sp2$ and $z'\in\Tr_0(1)$, $z'\neq z$.
We want to produce a simple path $\gamma\subset\Tr$ starting at $z'$ and
tending to $z$.
If $z\in\Tr$, then $z\in\Tr_0(1/n)$ for some $n\in\Np$,
and the existence of $\gamma$ is clear.  So suppose that $z\notin\Tr$.
For each $n\in\Np$, there is a point $z_n\in\Tr_0(1/n)$
which is within distance $2/n$ from $z$.
Let $\beta_n$ be the arc from $z'$ to $z_{n}$ in $\Tr_0(1/n)$.
For each $n$ and $m$, the intersection $\gamma_n\cap\Tr_0(1/m)$ is a
simple path.  Since $\Tr_0(1/m)$ is a compact finite topological
tree, there is a subsequence $\gamma_{n_j}$ such that
for each $m$ the Hausdorff limit $\lim_j\bigl(\gamma_{n_j}\cap\Tr_0(1/m)\bigr)$
exists, and is a simple path.  Because 
$\gamma_n-B(z,3/m)\subset \Tr_0(1/m)$ when $n\geq m$, it now
follows that $\gamma:=\lim_j\gamma_{n_j}$ is a simple
path.  Moreover, it is clear that $z\in\gamma$ and $\gamma-\{z\}\subset\Tr$.
\Qed

\proofof t.ustlim
We first prove that $\Tr$ is a topological tree.
Clearly, the trunk is
arcwise connected, since $\Tr:=\bigcup_n\Tr_0(1/n)$,
and each $\Tr_0(1/n)$ is arcwise connected.
It is also clear that the trunk is dense in $\Sp2$.
Let $x,y\in\Tr$.  Then there
is some $n\in\N$ such that $x,y\in\Tr_0(1/n)$, and there is a
unique arc $\gamma_1$ joining $x$ and $y$ in $\Tr_0(1/n)$.
Let $\gamma_2$ be an arc joining $x$ and $y$ in $\Tr$.
Since the dual trunk is dense, it must intersect all connected components
of $\Sp2-(\gamma_1\cup\gamma_2)$.  
Since the dual trunk is disjoint from the trunk, it does
not intersect $\gamma_1\cup\gamma_2$.  Because
the dual trunk is connected,
it now follows that $\Sp2-(\gamma_1\cup\gamma_2)$ is connected.
Consequently, $\gamma_1=\gamma_2$, and the trunk is uniquely
arcwise connected.

Let $n\in\N$ and
let $t$ be the spherical distance between
$\Tr_0(1/n)$ and $\d {\Tr}_0(1/n)$.  Since
these are compact and disjoint, $t>0$.
If $x\in\Tr_0(1/n)$ and there is a $y\in\Tr$
such that there is no path in $\Tr\cap B(x,3/n)$ joining $x$
and $y$, then there must be a path in $\d{\Tr}_0(1/n)$
separating $x$ and $y$ in $B(x,1/n)$.
(Indeed, if $\gamma$ is the path joining
$x$ and $y$ and $p\in\gamma- B(x,3/n)$,
then the path $\beta\subset\d\Tr$ connecting
two points $p'$ and $p''$ that are near $p$ and
are separated from each other by $\gamma$ near $p$,
will have a subarc in $\d{\Tr}_0(1/n)$
separating $x$ and $y$ in $B(x,1/n)$.)
  Consequently, for every point $x\in\Tr_0(1/n)$ and every
$y\in\Tr\cap B(x,t)\cap B(x,1/n)$ there is a path in
$\Tr\cap B(x,3/n)$ joining $x$ and $y$.  Hence,
the union of all arcs that contain
$x$ and are contained in $\Tr\cap B(x,3/n)$
is an arcwise connected subset of 
$\Tr\cap B(x,3/n)$ which contains 
$\Tr\cap B(x,t)\cap B(x,1/n)$.
This implies that $\Tr$ is locally arcwise connected,
and so it is a topological tree.  It is obviously dense
in $\Sp2$, and the proof of part (iii) is complete.

It is clear that for every $a,b\in\Sp2$ there is
some $\omega$ such that $(a,b,\omega)\in\ST$.

Let $\ST(\eps)$ be a subsequential scaling limit of $\ST_\delta(\eps)$.
We prove that a.s.\ every $\path$ such that $(a,b,\path)\in\ST(\eps)$
for some $a,b\in\Sp2$ is a simple path.
Let $\eps_1>0$.  It suffices to prove that the above statement
holds with probability at least $1-\eps_1$.
Let $V_1\subset\Sp2$ be a finite set of points,
and for each $\delta>0$ let $V_1^\delta$ be a collection
of vertices of $\dZ$, each close to one point of $V_1$,
and with $|V_1|=|V_1^\delta|$.
By \ref t.fin/, $V_1$ may be chosen
so that with probability at least $1-\eps_1$
the subtree of $\hT_\delta$ spanned by $V_1^\delta$
contains $\Tr_\delta(\eps)$, for all sufficiently small $\delta$.
This implies that each $\path_{a,b}(\eps)$ is a subarc
of $\path_{v,u}$ for some $v,u\in V_1^\delta$.
Because for every pair of points $v,u\in V_1$ the
scaling limit of the LERW from $v$ to $u$ is a simple path,
it follows that with probability at least $1-\eps_1$
for each $(a,b,\path)\in\ST(\eps)$, $\path$ is a simple path.

We may now conclude that a.s.\ for every $(a,b,\path)\in\ST$,
the set $\path-\bigl(\br B(a,\eps)\cup\br B(b,\eps)\bigr)$ is a 
$1$-manifold; that is, a disjoint union of simple paths.
Therefore, $\path':=\path-\{a,b\}$ is a $1$-manifold.
This means that each component of $\path'$
is an arc with endpoints in $\{a,b\}$.
It is clear that $\path$ may be oriented as a path
from $a$ to $b$. 
Suppose that $\path$, visits $a$ more than once.
If $a\neq b$, it then follows that $a\in\Tr$, and there is a simple
closed path in $\Tr$ containing $a$.  This is impossible,
since $\Tr$ is a topological tree.  Hence $a$,
and similarly $b$, are each visited only once in $\path$, which implies
that $\path$ is a simple path if $a\neq b$.
If $a=b$, the only possibility is that $\path=\{a\}$ or
that $\path$ is a simple closed path.
This proves the first and second statements in (ii).

Observe that if there is
simple curve $\alpha\subset\Tr$ such that $\br\alpha=\alpha\cup\{a\}$,
then $a$ must be in the dual trunk, for the dual trunk
is connected, intersects both components of 
$\Sp2-\br\alpha$, and is disjoint from $\Tr$.
This is a rare event, by \ref r.closetopt/ (or \ref c.ustdim/).
This proves (ii).

\ref c.stmaxdeg/ proves (iv).

The first claim in (i) is obvious.
Suppose that $a\neq b$ are such that there are two sets $\path$ and 
$\path'$ with $(a,b,\path),(a,b,\path')\in\ST$.  We know that
$\path$ and $\path'$ are simple paths.  If $\path\neq\path'$,
then there is a simple closed path, say $\gamma$,
contained in $\br\path\cup\br\path'$.
But as above, $\gamma$ must intersect the dual trunk,
since the dual trunk is connected and dense.
This implies that $a$ or $b$ are in the dual trunk.
This completes the proof of (i), and of the theorem.
\Qed 

\procl r.uniqp \procname{Uniqueness of paths}
We have seen in the above proof that the path in $\Tr$
from $a$ to $b$ is unique when $\{a,b\}\cap\d\Tr=\emptyset$.
The converse is also easily established.
\endprocl

\procl r.constructdual
\procname{Reconstructing $\d\Tr$}
It can be shown that the scaling limit dual trunk can be
reconstructed from the trunk.  
This can be seen from \ref r.uniqp/.  
Another description of the dual trunk from the trunk is
as follows.  Given distinct
$x,y\in\Sp2-\Tr$, let $\gamma(x,y)$ be the (unique) arc in $\Sp2-\Tr$
with endpoints $x,y$.  Then
$$
\d\Tr = \bigcup\Bigl\{\gamma(x,y)-\{x,y\}\st x\neq y,\, x,y\in\Sp2-\Tr\Bigr\}
\,.
$$
To prove this, it suffices to establish that $\gamma(x,y)$ is unique,
which follows from the fact that $\Tr$ is connected and dense.
\endprocl

\procl r.overtree
Consider the metric $d^*$ on $\Tr$, where $d^*(x,y)$ 
is the spherical diameter of the unique (possibly degenerate) arc
joining $x$ and $y$ in $\Tr$.  Since $\Tr$ is locally arcwise connected,
this new metric on $\Tr$ is compatible with the topology of $\Tr$ as
a subset of $\Sp2$.  Let $\Tr_*$ denote the completion of
this metric.  Then $\Tr_*$ is a compact topological tree,
and is naturally homeomorphic
with the ends compactification of $\Tr$.  Since $d^*$
majorizes the spherical metric, there is a natural projection
$\pi:\Tr_*\to\Sp2$, whose restriction to $\Tr$ is the identity.
It is easy to see that every point $p\in\Sp2-\d\Tr$ has
a unique preimage under $\pi$, and for points $p\in\d\Tr$,
the degree of $p$ in $\d\Tr$ is equal to $|\pi^{-1}(p)|$.

Consider some $o\in\Sp2$, and let $\ST^o$ be the appropriate
``slice'' of $\ST$; that is,
$\ST^o:=\bigl\{(b,\path)\st (o,b,\path)\in\ST\bigr\}$.
One can show that if $o\notin\d\Tr$, then
$\ST^o$ is homeomorphic with $\Tr_*$, and
$\pi:\ST^o\to\Sp2$ is the projection onto the first  coordinate,
when $\ST^o$ is identified with $\Tr_*$  through  this
homeomorphism.
\endprocl

\bsection{Free and wired trunks and conformal invariance}{s.ustci}
\resetConsts

We now want to give a precise formulation to a 
conformal invariance conjecture for the UST scaling limit,
and prove that it follows from the conjectured conformal invariance
of the LERW scaling limit.
(Such conformal invariance conjectures seem to be floating in the
air these days, with roots in the physics community.)
The conformal automorphisms of $\Sp2$ are M\"obius transformations.
We conjecture that the different notions of scaling limits of UST
in $\Sp2$, which where introduced in the previous section, exist
(without a need to pass to a subsequence) and are invariant under
M\"obius transformations.  
Moreover, the scaling limits in subdomains $D\subset \Sp2$ should be invariant
under conformal homeomorphisms $f: D\to D'\subset\Sp2$.  This
is a significantly stronger statement, since the M\"obius transformations
of $\Sp2$ form a $6$-dimensional group, while the space of conformal
homeomorphisms from the unit disk onto subdomains of $\Sp2$ is
infinite dimensional.

To formulate more precisely the invariance under conformal homeomorphisms
of subdomains $f:D\to D'$, we need first to discuss UST scaling limits in 
subdomains of $\Sp2$.  This will be now explained.  

For simplicity, we restrict our attention to simply connected domains
$D\subsetneqq\R^2$ whose boundary $\bd D$ is a simple closed path.
Let $p_0\in D$ be some basepoint.
Let $\FT_\delta^D$ be the uniform spanning tree of $G(D,\delta)$,
with free boundary conditions, and let $\WT_\delta^D$ be the
uniform spanning tree of $G(D,\delta)$ with
wired boundary conditions.
Let $\SpD$ be the metric space obtained from $\Sp2$ by contracting
$\Sp2-D$ to a single point.  Then we may think of $\WT_\delta^D$
as a random point in $\Ha(\SpD)$, which is a.s.\ a tree.
The tree $\WT_\delta^D$ may be thought of as a random point in
$\Ha(\br D)$, which is a.s.\ a tree.

Let $\WST_\delta^D:=\ST(\WT_\delta^D)$ and
$\FST_\delta^D:=\ST(\FT_\delta^D)$; that is,
the {\bf wired paths ensemble} $\WST_\delta$ is defined from
the wired tree $\WT_\delta$ in exactly the same way that the
ordinary paths ensemble $\ST_\delta$ was defined from $\hT_\delta$,
and similarly for $\FST_\delta$.
Note that $\WST_\delta\in\Ha\bigl(\SpD\times\SpD\times\Ha(\SpD)\bigr)$
and
$\FST_\delta\in\Ha\bigl(\br D\times\br D\times\Ha(\br D)\bigr)$.
Also the definitions of the scaling limits and the trunk
are the same as in the previous section. 

\procl t.dom
Let $D\subset\R^2$ be a domain whose boundary is a $C^1$-smooth simple
closed curve.  
\beginitems
\itemrm{(i)} \ref t.fin/, with $\Sp2$ replaced by $D$,
  holds for the free and wired spanning trees in $D$.
\itemrm{(ii)} The free scaling limit trunk in $D$ is disjoint
  from $\partial D$, in every (subsequential) scaling limit.
\itemrm{(iii)} The free scaling limit trunk in $D$ is disjoint
  from the scaling limit trunk of the dual tree (which is wired),
  in every (subsequential) scaling limit.
\endprocl

There are simply connected domains where $\bd D$ is not
a simple closed curve and (ii) fails:
the domain $(0,2)\times(0,2) - \bigcup_{n=1}^\infty (0,1]\times\{1/n\}$
is an example.

\procl l.freefin
There is an absolute constant $C>0$ such that the following holds
true.
Let $D$ be as in \ref t.dom/.  Then there
is a $\delta_0=\delta_0(D)>0$ with the following property.
Suppose that $\delta$ and $\delta_1$
are numbers satisfying $0<\delta\leq\delta_1\leq \delta_0$ and
$A$ is a connected subgraph of $\G(D,\delta)$ with diameter
at least $\delta_1$. Further suppose that
$p\in\verts\bigl(\G(D,\delta)\bigr)$ has distance $\delta_1$ to $A$.
Then the probability that a random walk on $\G(D,\delta)$
starting at $p$ will get to distance $C\delta_1$ from $p$ before
hitting $A$ is less than $1/2$.
\endprocl

\proof The proof is similar to the proof of \ref l.harm/.
Let $Z=Z(t)$ be the set of vertices with distance at least
$t$ from $p$, where we take $t>4\delta_1$
Let $A'$ be a component of $A\cap B(p,3\delta_1)$ containing some point at
distance $\delta_1$ to $p$ and having diameter at least $\delta_1$.
For $v\in\verts\bigl(\G(D,\delta)\bigr)$,
let $h(v)$ be the probability that a random walk starting
from $v$ will reach $Z$ before hitting $A'$.  
Then $h$ is discrete-harmonic, and minimizes Dirichlet energy
among functions that are $1$ on $Z$ and $0$ on $A'$.
As in the proof of \ref l.harm/, it follows that the
Dirichlet energy of $h$ is at most $O(1)/\log(t/\delta_1)$.
Let $B$ be the set of vertices $v$ at distance at most
$2\delta_1$ from $p$ such that $h(v)\geq h(p)$, and
let $B'$ be the component of $B$ containing $p$. 
Note that the diameter of $B'$ is at least $\delta_1$,
as $B'$ must neighbor with some vertex with distance
$2\delta_1$ from $p$, by the maximum principle for $h$.
As in the proof of \ref l.harm/, it can be shown that
when $\delta_1$ is sufficiently small 
(how small depends on the scale in which $\partial D$ appears smooth),
one can find $O(1)/\delta$ disjoint paths
in $\delta\Z^2\cap D$ connecting
$A'$ to $B'$, each of combinatorial length $O(1)/\delta$.
Because $h$ is zero on $A'$ and at least $h(p)$
on $B'$, it follows that the Dirichlet energy of $h$ is
at least $O(1)h(p)$.
We conclude that $h(p)\leq O(1)/\log(t/\delta_1)$,
which proves the lemma.
\Qed

\proofof t.dom
The proof for (i) in the wired case is the same as the
proof of \ref t.fin/ (and we don't need to assume anything about 
$D$).  The free case is the same, except that
one needs to appeal to \ref l.freefin/.
Assuming (ii), the proof of (iii) is identical to the proof of
\ref t.noint/.  The proof of (ii) is also the same as the
proof of \briefref t.noint/, except that one needs to find
the appropriate substitutes for \ref l.3d2/ and \ref c.maxdeg/;
namely, for every $\eps>0$ there is an $\eps_0>0$ such that
for all $\delta\in(0,\eps_0)$ 
with probability at least $1-\eps$
all points of degree three in the $\eps$-trunk of $\FT_\delta^D$ have distance
at least $\eps_0$ from $\partial D$,
and the probability that there is a point $p\in\partial D$
such that there are two disjoint crossings of
the annulus $\An(p,\eps_0,\eps)$ in the
$\eps$-trunk of $\WT_\delta^D$ is at most $\eps$.
The latter statement follows from the proof of \ref c.maxdeg/.

It remains to prove the appropriate substitute for
\briefref l.3d2/.  Consider three distinct points in
$D$, $p_1,p_2,p_3$, and for $\delta>0$ let $p_1',p_2',p_3'$
be a triple of points in $\delta\Z^2$ which is close to $p_1,p_2,p_3$,
respectively.  Let $m$ be the meeting point of $p_1',p_2',p_3'$ in
$\FT_\delta^D$, and let $B$ be any disk whose center is in $\partial D$
and which does not intersect $\{p_1',p_2',p_3'\}$.
Let $B'$ and $B''$ be disks concentric with $B$ of $1/2$ and $1/4$ of
its size, respectively.
By part (i), it suffices to prove that with probability going to $1$
as $\eps_0\to 0$ and $\delta\to 0$,
$m$ is not within distance $\eps_0$ from $\partial D\cap B''$.

Suppose that $m\in B$.  Let $\gamma_1,\gamma_2,\gamma_3$
be the arcs of $\FT_\delta^D$ that join $m$ to $p_1',p_2',p_3'$,
respectively, and let $\gamma_j'$ be the largest initial
segment of $\gamma_j$ that is contained in $B$, $j=1,2,3$.
There is a unique $k\in\{1,2,3\}$ such that
$\bigcup_{j\neq k}\gamma_j'$ 
separates $\gamma_k'$ from $\partial D$ in $B$.
By symmetry, it suffices to estimate the probability that
$m$ is close to $B\cap\partial D$ and $k=3$.  
We generate $m$ in the following way. 
Let $\gamma$ be a LERW from $p_1'$ to $p_2'$ in $\delta\Z^2\cap D$. 
Let $X$ be a random walk on $\delta\Z^2\cap D$ starting at $p_3'$,
let $\tau_\gamma$ be the first time $t$ where $X(t)\in\gamma$,
and let $\tau_1$ be the first time $t$ when $X(t)$ is incident
with an edge intersecting $\partial D\cap B'$.
Then we may take $m=X(\tau_\gamma)$.
Consider the event $\ev A$ where $m\in B$, $k=3$ and $\tau_1\geq\tau_\gamma$.
With high probability, the points $X(t)$, $t\leq\tau_1$,
which are close to $B'\cap\partial D$, are also close
to $X(\tau_1)$.  This is just a property of simple random walk
absorbed at $B'\cap\partial D$.
Consequently, on $\ev A$, with high probability,
if $m$ is close to $B\cap\partial D$, then $\gamma$
passes close to $X(\tau_1)$. 
Since the probability that $\gamma$ passes near any point,
which is not too close to $p_1'$ and $p_2'$, is small
(\ref r.closetopt/ applies here), and $X(\tau_1)$ is independent from
$\gamma$, we see that $\ev A$ has arbitrarily small probability.

We now need to consider the case $k=3$ and $\tau_1<\tau_\gamma$.
Let $\tau_1'$ be the first $t\geq\tau_1$ such that
$X(t)\notin B$, and let $\tau_2$ be the first $t\geq\tau_1'$
such that $X(t)$ is incident with an edge intersecting $\partial D\cap B'$.
Inductively, let $\tau_n'$ be the first $t\geq \tau_n$
such that $X(t)\notin B$ and let $\tau_{n+1}$ be the first
$t\geq\tau_n'$ such that $X(t)$ is incident with an edge
intersecting $\partial D\cap B'$.  Note that if $k=3$ and
$\tau_\gamma>\tau_1$, then $\tau_\gamma>\tau_1'$, for the random
walk $X$ must go around $\gamma_1'\cup\gamma_2'$ before hitting
$\gamma$.
Similarly, if $k=3$ and $\tau_\gamma>\tau_1$, then
$\tau_\gamma\in(\tau_n',\tau_{n+1})$ for some $n\in\N$.
If we fix a finite $k\in\N$, then the same argument as above
shows that with high probability, $\gamma$ does not pass close to the set
$\bigl\{X(\tau_n)\st n=1,2,\dots,k\bigr\}$.
Because $\P[\tau_k'<\gamma]\to 0$ as $k\to\infty$, uniformly in $\delta$,
the required result follows.
\Qed

Suppose that $D$ and $D'$ are two domains in $\R^2$ such that
the boundaries $\bd D,\bd D'$ are simple closed paths in $\R^2$.
Then there is a conformal homeomorphism $f:D\to D'$.
Moreover, $f$ extends continuously to a homeomorphism
of $\br D$ onto $\br D'$, which we will also denote by $f$.
It follows that $f$ induces maps
\begineqalno
f_W
&:\Ha\bigl(\SpD\times\SpD\times\Ha(\SpD)\bigr)\to
\Ha\bigl(\SS^2_{D'}\times\SS^2_{D'}\times\Ha(\SS^2_{D'})\bigr) 
\,,
\cr
f_F
&:\Ha\bigl(\br D\times\br D\times\Ha(\br D)\bigr)\to
\Ha\bigl(\br D\times\br D\times\Ha(\br D)\bigr)
\,.
\cr
\endeqalno

\procl t.confust
Let $D\subset\R^2$ be a domain whose boundary is a $C^1$-smooth simple
closed path.
Assuming \ref g.confinv/, the following is true.
\beginitems
\itemrm{(i)}
  The free and the wired UST scaling limits, $\FST,\WST$, in $D$ exist.
  (That is, do not depend on the sequence of $\delta$ tending to $0$.)
\itemrm{(ii)}
  If $f:D\to D'$ is a conformal homeomorphism between such domains,
  then $f_W$ is measure preserving from the law of $\WST$ in
$D$ to the law of $\WST$ in $D'$, and similarly for free boundary conditions.
\enditems
\endprocl

\proof
The proof for wired boundary conditions follows from Wilson's algorithm
and \ref t.dom/.
The easy details are left to the reader.  

For the free boundary conditions, observe that $\WT_\delta^D$
is dual to $\FT_\delta^D$ (on the dual grid). 
\ref r.constructdual/ is
also valid in the present setting, and shows that the free scaling
limit trunk can be reconstructed from the wired scaling limit trunk.
It is easy to see that the free scaling limit $\FST$ can be reconstructed
from the trunk.
Hence, conformal invariance of the wired UST implies conformal invariance
of the free.
\Qed

\bsection{Speculations about the Peano curve scaling limit}{s.peano}
This section will discuss the Peano curve winding
between the UST and its dual. 
{}From here on, the discussion will be somewhat speculative,
and we omit proofs, not because the proofs are particularly hard,
but because the paper is long enough as it is, and it is not clear
when another paper on this subject will be produced.

This Peano curve was briefly mentioned in \ref b.BLPS:usf/.
Consider the set of points $\Peano_\delta\subset\R^2$ which have the same
Euclidean distance from $\hT_\delta$ as from its dual $\d\hT_\delta$.
It is easy to verify that $\Peano_\delta$ is a simple path in
a square grid $\G_P(\delta)$, of mesh $\delta/2$, which visits all the vertices
in that grid.  Set $\hat \Peano_\delta:=\Peano_\delta\cup\{\infty\}$.
Then $\hat\Peano_\delta$ is a.s.\ a simple closed path in $\Sp2$
passing through $\infty$.

To consider the scaling limit of $\hat\Peano_\delta$, it is no
use to think of it as a set of points in $\Sp2$, because then the
scaling limit will be all of $\Sp2$.  Rather, one needs to
parameterize $\hat\Peano_\delta$ in some way.
One natural parameterization would be by the area of its
$\delta/2$-neighborhood, but there are several other plausible
parameterizations.  Another, more sophisticated approach,
would be to think of $\hat\Peano_\delta$ as defining a circular
order on the set $\hat\Peano_\delta$.
The circular order $R_\delta$ is a closed subset of $\bigl(\Sp2)^4$,
and $(a,b,c,d)\in R_\delta$ iff $a,b,c,d\in\hat\Peano_\delta$ and
$\{a,c\}$ separates $b$ from $d$ on $\hat\Peano_\delta$.
Then the (subsequential) scaling limit of $\hat\Peano_\delta$
may be taken as the weak limit of the law of $R_\delta$
in $\Ha\Bigl(\bigl(\Sp2\bigr)^4\Bigr)$.

Let $\Peano$ denote the Peano curve scaling limit, defined 
 as a path, or as a circular order, or some other reasonable definition.
Here is what we believe to be a description of $\Peano$,
in terms of the scaling limit of the UST.
Recall that in \ref r.overtree/, we have introduced a completion
$\Tr_*$ of the trunk, in the metric $d^*$, where the distance
between any two points of $\Tr$ is the diameter of the arc
connecting them, and that $\pi:\Tr_*\to\Sp2$ is the natural
projection.  Consider the joint distribution of
$\Tr_*$ and the dual $\d\Tr_*$, and let $\d\pi$
denote the projection $\d\pi:\Tr_*\to\Sp2$.
Let $\tilde\Peano$ be the set of points $(p,q)\in\Tr_*\times\d\Tr_*$
such that $\pi p=\d\pi q$.  Then $\tilde\Peano$ is a simple
closed path, and the map $\tilde\Peano\mapsto\Sp2$
defined by $(p,q)\mapsto \pi p$ gives the scaling limit $\Peano$.

Fix some $a,b\in\R^2$, $a\neq b$. 
Let $\path_{a,b}$ be the path such that $(a,b,\path_{a,b})\in\ST$
(this $\path_{a,b}$ is a.s.\ unique), and let
$\d\path_{a,b}$ be such that $(a,b,\d\path_{a,b})\in\d\ST$.
Then $\path_{a,b}$ and $\d\path_{a,b}$ are a.s.\ simple paths.
Let $D_1$ and $D_2$ be the two components of 
$\Sp2-(\path_{a,b}\cup\d\path_{a,b})$.  A.s.\ $\infty\in D_1\cup D_2$,
and without loss of generality take $\infty\in D_2$.
It is then clear that the part of $\Peano$ which is
between $a,b$ and does not contain $\infty$ is $\br D_2$,
and that the part which does contain $\infty$ is $\br D_1$.
Suppose that we {\bf condition} on $\path_{a,b}$ and on $\d\path_{a,b}$,
and look at some point $c\in D_2$.  We'd like to know the distribution
of the part of $\Peano$ between $c$ and $a$ which does not include $b$, say.
Recall that on a finite planar graph, we may generate the
UST and the dual UST by a modification of Wilson's algorithm,
where at each step in which we start from a vertex in the graph,
the dual tree built up to that point acts as a free boundary component,
and the tree built up to that point acts as an absorbing wired boundary
component, and at steps in which we start from a dual vertex, the
tree built up to that point acts as a free boundary component,
and the dual tree built up to that point acts as an absorbing wired
boundary component.  Consequently, we let $\alpha$ be the scaling
limit of LERW  on $\dZ-\d\path_{a,b}$
starting at $c$ that stops when it hits $\path_{a,b}$.
Then we let $\d\alpha$ be the scaling limit of LERW on
$\dZ-\bigl(\alpha\cup\path_{a,b}\bigr)$ starting at $c$ that stops
when it hits $\d\path{a,b}$.  (Since $c\in\alpha$, to define this
requires taking a limit as the starting point tends to $c$.)
Then $\alpha\cup\d\alpha$ separates $D_2$ into two regions,
say $D_2'$ and $D_2''$, and if $b\notin\br D_2'$, then $D_2'$
is the part of $\Peano$ ``separated'' from $b$ by $\{a,c\}$.

In the above construction, the domain $D_2$ was considered with
mixed boundary conditions.  One arc of $\bd D_2-\{a,b\}$ was taken
as wired, while the other was free.  The resulting Peano path
scaling limit $\Peano$ is a path joining $a$ and $b$ in $\br D_2$.
{} From \ref g.confinv/ should follow a conformal invariance
result for UST in such domains with mixed boundary conditions.
Therefore, having an understanding of the law of the
Peano curve for one triplet $(D,a,b)$, where 
$a$ and $b$ are distinct points in $\bd D$, which is a simple closed curve,
suffices for any other such triplet.  This suggests that we should take
the simplest possible such configuration; that is,
$D=\H$, the upper half plane, $a=0,b=\infty$.
Suppose that we then take a point $c\in\H$ and condition on the part
of $\Peano$ between $0$ and $c$ and separated from $\infty$.
The effect of that on $\Peano$ in the remaining subdomain $D_c$ of $\H$
is all in the boundary $\bd D_c$.  This is a kind of Markovian property
for the Peano curve, similar to the property given by \ref l.dprod/.
By taking the conformal map from $D_c$ to $\H$, which fixes $\infty$,
takes $c$ to $0$, and is appropriately normalized at $\infty$, we
may return to base one.  This suggests that,
as we have claimed in the introduction, a representation of
the Peano curve scaling limit similar to the SLE representation
of the LERW from $0$ to $\bd U$ which we have introduced.
The analogue of the \lowner{} differential equation for this situation
is \ref e.loewh/.
Due to the Markovian nature of the Peano curve, the corresponding parameter
$\zeta$ in \ref e.loewh/ should have the form $\zeta=\B(\lec t)$,
where $\B$ is Brownian motion on $\R$ starting at $0$,
and $\lec$ is some constant.  One can, in fact, show that $\lec=8$,
by deriving an appropriate analogue of Cardy's \ref b.Cardy/ conjectured
formula,
using the representation \ref e.loewh/ and the techniques of \ref s.crit/.
The details will appear elsewhere.


\bibsty{mybibstyle}
\bibfile{\jobname.bib}

\startbib{MMOT92}

\bibitem[Aiz]{Aizenman:web}
{\sc M.~Aizenman}.
\newblock Continuum limits for critical percolation and other stochastic
  geometric models.
\newblock
Preprint.
\newblock \urlref{http://xxx.lanl.gov/abs/math-ph/9806004}.

\bibitem[ABNW]{ABNW}
{\sc M.~Aizenman, A.~Burchard, C.~M. Newman, and D.~B. Wilson}.
\newblock Scaling limits for minimal and random spanning trees in two
  dimensions.
\newblock
Preprint.
\newblock \urlref{http://xxx.lanl.gov/abs/math/9809145}.

\bibitem[ADA]{ADA}
{\sc M.~Aizenman, B.~Duplantier, and A.~Aharony}.
\newblock Path crossing exponents and the external perimeter in 2{D}
  percolation.
\newblock
Preprint.
\newblock \urlref{http://xxx.lanl.gov/abs/cond-mat/9901018}.

\bibitem[Ald90]{Aldous:ust}
{\sc D.~J. Aldous}.
\newblock The random walk construction of uniform spanning trees and uniform
  labelled trees.
\newblock {\em SIAM J. Discrete Math. 3}, 4 (1990), pages 450--465.

\bibitem[Ben]{Benjamini:ustsl}
{\sc I.~Benjamini}.
\newblock Large scale degrees and the number of spanning clusters for the
  uniform spanning tree.
\newblock In M.~Bramson and R.~Durrett, editors, {\em Perplexing Probability
  Problems: Papers in Honor of Harry Kesten}, Boston. Birkh{\"a}user.
\newblock To appear.

\bibitem[BLPS98]{BLPS:usf}
{\sc I.~Benjamini, R.~Lyons, Y.~Peres, and O.~Schramm}.
\newblock Uniform spanning forests.
\newblock
Preprint.
\newblock
  \urlref{http://www.wisdom.weizmann.ac.il/\string~schramm/papers/usf/}.

\bibitem[BJPP97]{BJPP}
{\sc C.~J. Bishop, P.~W. Jones, R.~Pemantle, and Y.~Peres}.
\newblock The dimension of the {B}rownian frontier is greater than $1$.
\newblock {\em J. Funct. Anal. 143}, 2 (1997), pages 309--336.

\bibitem[Bow]{Bowditch:treelike}
{\sc B.~H. Bowditch}.
\newblock Treelike structures arising from continua and convergence groups.
\newblock {\em Mem. Amer. Math. Soc.\/}.
\newblock To appear.

\bibitem[Bro89]{Broder:ust}
{\sc A.~Broder}.
\newblock Generating random spanning trees.
\newblock In {\em Foundations of Computer Science}, pages 442--447, 1989.

\bibitem[BP93]{BP:ust}
{\sc R.~Burton and R.~Pemantle}.
\newblock Local characteristics, entropy and limit theorems for spanning trees
  and domino tilings via transfer-impedances.
\newblock {\em Ann. Probab. 21}, 3 (1993), pages 1329--1371.

\bibitem[Car92]{Cardy}
{\sc J.~L. Cardy}.
\newblock Critical percolation in finite geometries.
\newblock {\em J. Phys. A 25}, 4 (1992), pages L201--L206.

\bibitem[DD88]{Duplantier:Peano}
{\sc B.~Duplantier and F.~David}.
\newblock Exact partition functions and correlation functions of multiple
  {H}amiltonian walks on the {M}anhattan lattice.
\newblock {\em J. Statist. Phys. 51}, 3-4 (1988), pages 327--434.

\bibitem[Dur83]{Duren:book}
{\sc P.~L. Duren}.
\newblock {\em Univalent functions}.
\newblock Springer-Verlag, New York, 1983.

\bibitem[Dur84]{Durrett:BMM}
{\sc R.~Durrett}.
\newblock {\em Brownian motion and martingales in analysis}.
\newblock Wadsworth International Group, Belmont, Calif., 1984.

\bibitem[Dur91]{Durrett:probability}
{\sc R.~Durrett}.
\newblock {\em Probability}.
\newblock Wadsworth \& Brooks/Cole Advanced Books \& Software, Pacific Grove,
  CA, 1991.

\bibitem[EK86]{EK:markov}
{\sc S.~N. Ethier and T.~G. Kurtz}.
\newblock {\em Markov processes}.
\newblock John Wiley \& Sons Inc., New York, 1986.

\bibitem[Gri89]{Grimmett:book}
{\sc G.~Grimmett}.
\newblock {\em Percolation}.
\newblock Springer-Verlag, New York, 1989.

\bibitem[H{\"a}g95]{Hag:ustlim}
{\sc O.~H{\"a}ggstr{\"o}m}.
\newblock Random-cluster measures and uniform spanning trees.
\newblock {\em Stochastic Process. Appl. 59}, 2 (1995), pages 267--275.

\bibitem[It{\^o}61]{Ito:lectures}
{\sc K.~It{\^o}}.
\newblock {\em Lectures on stochastic processes}.
\newblock Notes by K.~M.~Rao.
\newblock Tata Institute of Fundamental Research, Bombay, 1961.

\bibitem[Jan12]{Janiszewski}
{\sc Janiszewski}.
\newblock {\em J. de l'Ecole Polyt. 16\/} (1912), pages 76--170.

\bibitem[Ken98a]{Kenyon:conf}
{\sc R.~Kenyon}.
\newblock Conformal invariance of domino tiling.
\newblock
Preprint.
\newblock \urlref{http://topo.math.u-psud.fr/\string~kenyon/confinv.ps.Z}.

\bibitem[Ken98b]{Kenyon:5o4}
{\sc R.~Kenyon}.
\newblock The asymptotic determinant of the discrete laplacian.
\newblock
Preprint.
\newblock \urlref{http://topo.math.u-psud.fr/\string~kenyon/asymp.ps.Z}.

\bibitem[Ken99]{Kenyon:longrange}
{\sc R.~Kenyon}.
\newblock Long-range properties of spanning trees.
\newblock Preprint.

\bibitem[Ken]{Kenyon:inprep}
{\sc R.~Kenyon}.
\newblock In preparation.

\bibitem[Kes87]{Kesten:hitting}
{\sc H.~Kesten}.
\newblock Hitting probabilities of random walks on $\Z \sp d$.
\newblock {\em Stochastic Process. Appl. 25}, 2 (1987), pages 165--184.

\bibitem[Kuf47]{Kufarev}
{\sc P.~P. Kufarev}.
\newblock A remark on integrals of {L}\"owner's equation.
\newblock {\em Doklady Akad. Nauk SSSR (N.S.) 57\/} (1947), pages 655--656.

\bibitem[LPSA94]{Langlands:Bull}
{\sc R.~Langlands, P.~Pouliot, and Y.~Saint-Aubin}.
\newblock Conformal invariance in two-dimensional percolation.
\newblock {\em Bull. Amer. Math. Soc. (N.S.) 30}, 1 (1994), pages 1--61.

\bibitem[Law93]{Lawler:makarov}
{\sc G.~F. Lawler}.
\newblock A discrete analogue of a theorem of {M}akarov.
\newblock {\em Combin. Probab. Comput. 2}, 2 (1993), pages 181--199.

\bibitem[Law]{Lawler:survey}
{\sc G.~F. Lawler}.
\newblock Loop-erased random walk.
\newblock In M.~Bramson and R.~Durrett, editors, {\em Perplexing Probability
  Problems: Papers in Honor of Harry Kesten}, Boston. Birkh{\"a}user.
\newblock To appear.

\bibitem[L{\"o}w23]{Lowner}
{\sc K.~L{\"o}wner}.
\newblock Untersuchungen {\"u}ber schlichte konforme abbildungen des
  einheitskreises, {I}.
\newblock {\em Math. Ann. 89\/} (1923), pages 103--121.

\bibitem[Lyo98]{Lyons:usfsurvey}
{\sc R.~Lyons}.
\newblock A bird's-eye view of uniform spanning trees and forests.
\newblock In {\em Microsurveys in discrete probability (Princeton, NJ, 1997)},
  pages 135--162. Amer. Math. Soc., Providence, RI, 1998.

\bibitem[MR]{MR:lowner}
{\sc D.~E. Marshall and S.~Rohde}.
\newblock In preparation.

\bibitem[MMOT92]{MMOT:treechar}
{\sc J.~C. Mayer, L.~K. Mohler, L.~G. Oversteegen, and E.~D. Tymchatyn}.
\newblock Characterization of separable metric ${{\R}}$-trees.
\newblock {\em Proc. Amer. Math. Soc. 115}, 1 (1992), pages 257--264.

\bibitem[MO90]{MO:treechar}
{\sc J.~C. Mayer and L.~G. Oversteegen}.
\newblock A topological characterization of ${{\R}}$-trees.
\newblock {\em Trans. Amer. Math. Soc. 320}, 1 (1990), pages 395--415.

\bibitem[New92]{Newman:planesets}
{\sc M.~H.~A. Newman}.
\newblock {\em Elements of the topology of plane sets of points}.
\newblock Dover Publications Inc., New York, second edition, 1992.

\bibitem[Pem91]{Pemantle:dichot}
{\sc R.~Pemantle}.
\newblock Choosing a spanning tree for the integer lattice uniformly.
\newblock {\em Ann. Probab. 19}, 4 (1991), pages 1559--1574.

\bibitem[Pom66]{Pommerenke:Lowner}
{\sc C.~Pommerenke}.
\newblock On the {L}oewner differential equation.
\newblock {\em Michigan Math. J. 13\/} (1966), pages 435--443.

\bibitem[Rus78]{Russo}
{\sc L.~Russo}.
\newblock A note on percolation.
\newblock {\em Z. Wahrscheinlichkeitstheorie und Verw. Gebiete 43}, 1 (1978),
  pages 39--48.

\bibitem[Sch]{Schramm:inprep}
{\sc O.~Schramm}.
\newblock In preparation.

\bibitem[Sla94]{Slade:survey}
{\sc G.~Slade}.
\newblock Self-avoiding walks.
\newblock {\em Math. Intelligencer 16}, 1 (1994), pages 29--35.

\bibitem[SW78]{SW}
{\sc P.~D. Seymour and D.~J.~A. Welsh}.
\newblock Percolation probabilities on the square lattice.
\newblock {\em Ann. Discrete Math. 3\/} (1978), pages 227--245.
\newblock Advances in graph theory (Cambridge Combinatorial Conf., Trinity
  College, Cambridge, 1977).

\bibitem[TW98]{TW}
{\sc B.~T{\'o}th and W.~Werner}.
\newblock The true self-repelling motion.
\newblock {\em Probab. Theory Related Fields 111}, 3 (1998), pages 375--452.

\bibitem[Wil96]{Wilson:alg}
{\sc D.~B. Wilson}.
\newblock Generating random spanning trees more quickly than the cover time.
\newblock In {\em Proceedings of the Twenty-eighth Annual ACM Symposium on the
  Theory of Computing (Philadelphia, PA, 1996)}, pages 296--303, New York,
  1996. ACM.

\endbib

\endreferences

\filbreak
\begingroup
\eightpoint\sc
\parindent=0pt 

\htmlref{http://www.wisdom.weizmann.ac.il/} {Mathematics Department},
\htmlref{http://www.weizmann.ac.il/} {The Weizmann Institute of Science},
Rehovot 76100, Israel

\medskip
\emailwww{schramm@wisdom.weizmann.ac.il}
{http://www.wisdom.weizmann.ac.il/\string~schramm/}
\bye